\documentclass[12pt]{iopart}

\usepackage{iopams,amsthm,amssymb} 
\usepackage{amsfonts} 
\usepackage{cite}
\usepackage{amsmath}
\usepackage[dvips]{graphicx}
\usepackage[ruled]{algorithm}
\usepackage{algorithmic}
\usepackage{bm}
\eqnobysec	   

\newcommand{\myi}{\ensuremath \rmi}		

\newcommand{\myvec}[1]{\boldsymbol{#1}}		
\newcommand{\unitvec}[1]{\hat {\myvec{#1}}}	

\newcommand{\mymat}[1]{{\mathsf{#1}}}

\newcommand{\dbyd}[2]{\frac{\partial {#1}}{\partial {#2}}}

\newcommand{\operator}[1]{\ensuremath {\mathcal{#1}}}
\newcommand{\myspace}[1]{\ensuremath {\mathrm{#1}}} 
\newcommand{\set}[1]{\left\{#1\right\}}

\newcommand{\labelmu}[1]{\mu_{\mathrm{#1}}}	
\newcommand{\mua}{\labelmu{a}}			
\newcommand{\mus}{\labelmu{s}'}			
\newcommand{\must}{\labelmu{s}}			

\newcommand{\phidiff}{\ensuremath \Phi}		
\newcommand{\pos}{\ensuremath r} 		
\newcommand{\posvec}{\ensuremath \myvec{\pos}}	
\newcommand{\ang}{\ensuremath s} 		
\newcommand{\angvec}{\unitvec{\ang}}		
\newcommand{\angvecp}{\angvec\,'}		
\newcommand{\domain}{\ensuremath \Omega} 	
\newcommand{\bdomain}{\partial \domain}		
\newcommand{\void}{\overline{\domain}} 		
\newcommand{\bvoid}{\partial \void}		
\newcommand{\bnormal}{\ensuremath \nu} 		
\newcommand{\bnormvec}{\unitvec{\bnormal}}	
\newcommand{\refmm}{\ensuremath R}              
\newcommand{\bccons}{\ensuremath A}		
\newcommand{\fembasis}{u}
\newcommand{\mneumop}[1]{\neumop{#1}}

\newcommand{\pmdf}{\ensuremath \rho}

\newcommand{\angint}[1]{\int_{S^{n-1}} \hspace{-1em} #1 \,\mathrm{d}\angvec}
\newcommand{\angintp}[1]{\int_{S^{n-1}} \hspace{-1em}  #1 \,\mathrm{d}\angvecp}

\newcommand{\nboundint}{\int_{\bdomain} }

\newcommand{\LM}{Levenberg-Marquardt~}

\newcommand{\Frechet}{Fr\'{e}chet}
\newcommand{\KZ}{Kaczmarz}

\renewcommand{\mneumop}[1]{\bnormvec \cdot \nabla {#1}}
\newcommand{\difftensor}{\mymat{K}}

\newcommand{\bea}{\begin{eqnarray}}
\newcommand{\eea}{\end{eqnarray}}
\newcommand{\bwea}{\endtwocolumns \topline \begin{eqnarray}}
\newcommand{\ewea}{\end{eqnarray} \botline \twocolumns }

\newcommand{\vintvar}[3]{\int_{#1} {#2} \,\mathrm{d}{#3}}
\newcommand{\sintvar}[3]{\int_{#1} {#2} \,\mathrm{d}{#3}}

\newcommand{\boundintvar}[2]{\nboundint {#1} \,\mathrm{d}{#2}}

\newcommand{\gboundintvar}[2]{\int_{\bgdomain} {#1} \,\mathrm{d}{#2}}

\newcommand{\gdomintvar}[2]{\int_{\gdomain} {#1} \,\mathrm{d}{#2}}
\newcommand{\domintvar}[2]{\int_{\domain} {#1} \,\mathrm{d}{#2}}
\renewcommand{\myvec}[1]{\bi{#1}}
\newcommand{\mvk}[2]{\myvec{#1}_{#2}}
\newcommand{\tmvk}[2]{\tilde{\myvec{#1}}_{#2}}

\newcommand{\elem}{\triangle}
\newcommand{\logparam}{\xi}
\newcommand{\logy}{\zeta}
\newcommand{\Prob}{\pi}

\newcommand{\Lgn}{\mathcal{J}}
\newcommand{\Krylsp}{\myspace{K}}
\newcommand{\Krylset}{\myset{K}}
\newcommand{\Krylvec}{\myvec{v}}
\newcommand{\MRF}{W}
\newcommand{\regln}{\alpha}
\newcommand{\step}{\tau}
\newcommand{\stepk}[1]{\tau_{#1}}
\newcommand{\cgp}{\beta}
\newcommand{\cgpk}[1]{\cgp_{#1}}
\newcommand{\reglr}{\Psi}
\newcommand{\mystack}[2]{{#1}_{#2}}
\newcommand{\sds}{\myvec{s}}
\newcommand{\sdsk}[1]{\mvk{s}{#1}}
\newcommand{\tsdsk}[1]{\tmvk{s}{#1}}

\newcommand{\cgsk}[1]{\mvk{p}{#1}}
\newcommand{\tcgsk}[1]{\tmvk{p}{#1}}
\newcommand{\mA}{\mymat{A}}
\newcommand{\fv}{\myvec{\param}^\delta}
\newcommand{\fvk}[1]{\myvec{\param}^{\delta({#1})}}
\newcommand{\xv}{\myvec{\param}}
\newcommand{\xvk}[1]{\mvk{\param}{#1}}
\newcommand{\yyv}{\myvec{\y}}
\newcommand{\yv}{\myvec{\y}^\delta}

\newcommand{\meas}{\y^{\mathrm{meas}}}
\newcommand{\tmA}{\tilde{\mA}}
\newcommand{\tfv}{\tilde{\myvec{\param}}}
\newcommand{\tfvk}[1]{\tmvk{\param}{#1}}
\newcommand{\tyv}{\tilde{\myvec{\y}}}
\newcommand{\tKrylvec}{\tilde{\Krylvec}}
\newcommand{\cvx}{\cvmat^{-1}_{\param}}
\newcommand{\cvy}{\cvmat^{-1}_{e}}
\newcommand{\icvx}{\cvmat_{\param}}
\newcommand{\icvy}{\cvmat_{e}}
\newcommand{\lcvx}{\mymat{L}_{\param}}
\newcommand{\lcvy}{\mymat{L}_{e}}
\def\mtrx#1#2{
  \left(
    \begin{array}{#1}
      #2
    \end{array}
  \right)}
\newcommand{\R}{\mathbb{R}}

\newcommand{\cvmat}{\mymat{\Gamma}}
\def\tr{^{\rm T}}
\def\th{^{\rm th}}

\def\diag{{\rm diag\,}}

\newcommand{\Id}{\mymat{I}}
\newcommand{\LMp}{\gamma}
\newcommand{\Obj}{\mathcal{E}}

\newcommand{\field}{U}

\newcommand{\adjfield}{\field^{\ast}}

\newcommand{\lgnfield}{Z}

\newcommand{\ptlop}{\myop{V}}
\newcommand{\ef}{e\rightarrow f}
\newcommand{\PDE}{\mathcal{L}} 
\newcommand{\rPDE}{\mathcal{W}} 
\newcommand{\U}{\Phi}     
\newcommand{\aU}{\U^{\ast}}     
\newcommand{\Ue}{\U^{e}}     
\newcommand{\Uf}{\U^{f}}     
\newcommand{\aUf}{\U^{f\ast}}     
\newcommand{\Uef}{\U^{\ef}}     
\newcommand{\aUef}{\U^{\ef\ast}}     
\newcommand{\aPDE}{\mathcal{L}^{\ast}} 
\newcommand{\param}{x}

\newcommand{\diffusivity}{k}

\newcommand{\fx}{h}
\newcommand{\flph}{p}
\newcommand{\pparam}{\param^{\delta}} 
\renewcommand{\bccons}{\zeta}
\newcommand{\sbs}{b}

\renewcommand{\r}{\myvec{r}}
\newcommand{\rs}{\r_{\mathrm{s}}}
\newcommand{\rd}{\r_{\mathrm{m}}}
\newcommand{\gr}{\myvec{\xi}}
\newcommand{\grs}{\gr_{\mathrm{s}}}
\newcommand{\grd}{\gr_{\mathrm{m}}}
\newcommand{\sr}{\myvec{r}}

\newcommand{\gdomain}{\Xi}
\newcommand{\bgdomain}{\partial\Xi}
\newcommand{\surf}{\Sigma}
\newcommand{\rtebnd}{B}
\renewcommand{\myspace}[1]{\mathrm{#1}}
\newcommand{\myop}[1]{\mathcal{#1}}
\newcommand{\myset}[1]{\mathbb{#1}}
\newcommand{\Jin}{J^-}
\newcommand{\Jout}{J^+}
\newcommand{\eJin}{J^{e-}}
\newcommand{\eJout}{J^{e+}}
\newcommand{\fJin}{J^{f-}}
\newcommand{\fJout}{J^{f+}}
\newcommand{\efJout}{J^{\ef+}}
\newcommand{\exsurf}{\surf_{\mathrm{ext}}}
\newcommand{\sbc}{\myop{B}^-}
\newcommand{\mbc}{\myop{B}^+}
\newcommand{\ambc}{\myop{B}^{+\ast}}
\newcommand{\tfn}{\Lambda}
\newcommand{\fmap}{F}
\newcommand{\dffmap}{F^{'}}
\newcommand{\affmap}{F^{'\ast}}
\newcommand{\ddffmap}{F^{''}}

\newcommand{\mfmap}{\operator{F}}
\newcommand{\dfmfmap}{\mfmap^{'}}
\newcommand{\afmfmap}{\mfmap^{'\ast}}

\newcommand{\aafmfmap}{\mfmap^{''\ast}}
\newcommand{\y}{y} 
\newcommand{\dfy}{y^{'}} 
\newcommand{\apind}{j} 
\newcommand{\ap}{w_{\apind}} 
\newcommand{\sind}{i}
\newcommand{\bdind}{\ell}

\newcommand{\itind}{n}
\newcommand{\s}{\Jin_{\sind}}
\newcommand{\m}{\Jout_{\sind}}
\newcommand{\Q}{\myspace{Q}}
\newcommand{\Z}{\myspace{Z}}
\newcommand{\diffcoef}{D}
\newcommand{\submua}[1]{ \mu_{{\mathrm a}_{#1}}}	
\newcommand{\msset}{{\boldsymbol \lambda}}
\newcommand{\sspace}{\myspace{X}}
\newcommand{\dspace}{\myspace{Y}}
\newcommand{\absspace}{\sspace_{\mua}}

\newcommand{\kapspace}{\sspace_{\diffcoef}}
\newcommand{\absspacel}[1]{\sspace_{{\mua}}^{#1}}

\newcommand{\kapspacel}[1]{\sspace_{{\diffcoef}}^{#1}}

\newcommand{\fluspace}{\sspace_{\fx}}

\newcommand{\mfpespace}{\sspace_{\mathrm{\bf \flph}}}

\newcommand{\domj}{\domain_{\sind}}

\newcommand{\norgrad}[1]{\frac{\partial #1 }{\partial \nu}}

\newcommand{\shape}{S}

\renewcommand{\diag}{{\mbox{{\rm diag}}}}

\newcommand{\shaperecon}{\hat{\shape}}

\newcommand{\shapeiter}[1]{\shape^{(#1)}}

\newcommand{\I}{\ensuremath U}
\newcommand{\J}{\ensuremath J}
\newcommand{\Brho}{\boldsymbol\rho}

\renewcommand{\Im}{{\mathrm {Im}}}
\renewcommand{\Re}{{\mathrm {Re}}}

\newcommand{\lap}{\nabla^{2}}
\newcommand{\curl}{\nabla \times}
\renewcommand{\div}{\nabla \cdot }
\newcommand{\grad}{\nabla}

\newcommand{\Br}{{\bm r}}

\newcommand{\Bp}{{\bm p}}

\newcommand{\Bq}{{\bm q}}
\newcommand{\Bk}{{\bm k}}

\newcommand{\BR}{{\bm R}}
\newcommand{\BE}{{\bm E}}

\newcommand{\BJ}{{\bm J}}
\renewcommand{\BR}{{\bm R}}

\newcommand\n{{\hat{\boldsymbol\nu}}}
\newcommand\svec{{\bf \hat s}}

\newcommand\z{{\bm \hat z}}

\newcommand\K{{\mathcal{K}}}   
\newcommand{\dist}{\hspace{.5pt}{\rm dist}}

\newtheorem{theorem}{Theorem}[section]

\theoremstyle{remark}

\theoremstyle{corollary}

\newcommand\half{\frac{1}{2}}

\begin{document}

\title{Optical tomography: forward and inverse problems}

\author{Simon R. Arridge}
\address{Department of Computer Science, University College London, Gower Street, London WC1E6BT, UK}
\ead{s.arridge@cs.ucl.ac.uk}

\author{John C. Schotland}
  
\address{Department of Bioengineering and Graduate Group in Applied
Mathematics and Computational Science, University of Pennsylvania,
Philadelphia, PA 19104, USA} 
\ead{schotland@seas.upenn.edu}

\begin{abstract}
This paper is a review of recent mathematical and computational advances in optical tomography. We discuss the physical foundations of forward models for light propagation on microscopic, mesoscopic and macroscopic scales. We also consider direct and numerical approaches to the inverse problems which arise at each of these scales. Finally, we outline future directions and open problems in the field.
\end{abstract}

\section{Introduction}
Optical tomography is a biomedical imaging modality that uses scattered light as a probe of structural variations in the optical properties of tissue. In a typical experiment, a highly-scattering medium is illuminated by a narrow collimated beam and the light which propagates through the medium is collected by an array of detectors. There are many variants of this basic scenario. For instance, the source may be pulsed or time-harmonic,  coherent or incoherent, and the illumination may be spatially structured or multispectral. Likewise, the detector may be time- or frequency-resolved, polarization or phase sensitive,  located in the near- or far-field and so on. The inverse problem that is considered is to reconstruct the optical properties of the medium from boundary measurements. The mathematical formulation of the corresponding forward problem is dictated primarily by \emph{spatial scale}, ranging from the Maxwell equations at the microscale, to the radiative transport equation at the mesoscale, and to diffusion theory at the macroscale. In addition, experimental time scales vary from the femtosecond on which light pulses are generated, through the nanosecond on which diffuse waves propogate, to the millisecond scale on which biological activation takes place and still longer for pathophysiologic changes. In this paper, we review the mathematical structure and computational approaches to the forward and inverse problems which arise at each of these scales. 

The state of progress on the inverse problem of optical tomography was 
reviewed a decade ago in this journal~\cite{arridge99}. Since then, the field 
has grown out of all proportion. It is almost impossible to comprehensively 
summarize the numerous new measurement systems, applications, algorithms, 
physical approximations and theoretical results. Nevertheless, in this paper 
we attempt to give an overview of those advances that have proven to be most 
signiÞcant, and to enumerate those topics which remain unresolved and which 
represent the most fruitful areas for forthcoming research. 

We have structured the paper as a tutorial, emphasizing aspects of wave propagation in random media, the mathematical structure of related inverse problems, and computational tool for image reconstruction. The selection of topics is to some extent personal and reflects the tastes and biases of the authors. We do not attempt to summarize the considerable body of work that concerns instrumentation or clinical applications of optical tomography. For some reviews of these aspects we refer the reader to several articles\cite{yodh95,hebden97a,arridge97a,hawysz2000,boas2001a,gibson2005a,leff2008}.
Nevertheless, it is important to note that there have been significant advances on the experimental front that have motivated and impacted on the corresponding theoretical developments.

\begin{itemize}
\item The introduction of noncontact imaging systems wherein a scanned beam and a lens-coupled CCD is employed to replace the illumination and detection fiber-optics of conventional optical tomography experiments~\cite{Ripoll_2004_1,Ripoll_2003_1,Cuccia_2005_1}. Such systems generate data sets that are orders of magnitude larger than those acquired with fiber-based systems, leading to significant computational challenges for image reconstruction. 
\item The development of continuously tunable high-power pulsed light sources has led to increased attention to the multispectral aspects of optical tomography.
\item The substantial growth in the availability of targeted fluorescent probes for molecular imaging has stimulated the development of computationally efficient algorithms for fluorescence optical tomography.
\end{itemize}
In addition, from a theoretical point of view, the rapid progress in the analysis of optical phenomena at smaller scales has allowed the field to move from simply considering diffuse transport as the model of light propogation to include more general radiative transport, coherence and polarization effects.

The paper is organized as follows. In section 2 we introduce some general concepts and an abstract framework within which it is possible to describe the class of inverse problems that arise in optical tomography. In section 3 we present the essential physical ideas that are need to the propagation of light in a random medium. We also describe the scattering theory of diffuse waves within radiative transport theory. In section 4 we define the various imaging modalities that we will study using the framework established in section 2. In section 5 we introduce numerical methods for computing the Green's functions for the RTE and diffusion equation in the context of forward modeling. In section 6 we describe direct reconstruction methods and associated fast algorithms for several inverse problems in optical tomography. Such methods are well suited to reconstructing images from the large data sets that are available in noncontact optical tomography systems. In section 7 we focus on numerical inversion methods using the tools of optimization theory and Bayesian statistics. In section 8 we introduce the basic concepts of shape-based methods for image reconstruction. In section 9 we survey some remaining topics including reconstruction of anisotropic and time-varying optical parameters. 

There are a number of topics that we do not discuss. These include thermoacoustic and photoacoustic tomography, acousto-optic imaging, coherent imaging and nanoscale optical tomography. Lack of space or expertise, rather than interest, prohibit the authors from treating these important subjects.

\section{A General Framework}
\label{sect:general_framework}
In this section we introduce some general concepts and notation for the forward and inverse problems occuring in optical tomography.
Since the various models are described in terms of different physical quantities (space, time or frequency, direction, wavelength), we attempt some abstract definitions in terms of a 
\emph{generalised coordinate} that we denote by $\gr$.
\subsection{Forward problems\label{sect:general_forward}}
We consider a domain $\domain$ with boundary $\bdomain$. To allow for generalised coordinates we define the generalised domains $\gdomain$ with boundary $\bgdomain$. Experiments are performed by applying inward directed photon currents $
\Jin(\grs), \grs \in \bgdomain$ and measuring outgoing photon currents $\Jout(\grd), \grd \in \bgdomain$. In the ``non-contact" application, these functions are applied and/or measured by free-space propogation from $\bgdomain$ to an exterior surface $\exsurf$. The propagation of light inside the domain is governed by a differential operator $\PDE(\param)$ parameterised by functions $\param(\sr), \sr \in \domain$. Let $\sspace$ denote the function space of these parameters.
In most cases the inverse problem is a parameter estimation problem for $\param$.

Let $\myset{S}$ represent a (possibly infinite) set of source functions on $\bgdomain$ or $\exsurf$. Let $\Q$ represent the space spanned by the functions in $\myset{S}$. Then for a particular applied source current $\s \in \myset{S}$ the measureable current is found by solving
\begin{eqnarray}
	\PDE(\param) \field_{\sind} & =& 0 \label{eq:gpde}\\
	\sbc \field_{\sind} &=& \s \label{eq:gin}\\
\m &=& \mbc \field_{\sind} \label{eq:gout}
\end{eqnarray}
where $\sbc$, $\mbc$ are the boundary conditions. We assume that in all cases considered herein, these boundary conditions are specific enough to ensure uniqueness of $\field_{\sind}$. The measureable currents belong to a space $\Z$.

In the case where $\s$ is a $\delta$-function located at $\grs$ the solution $\field_{\sind}$ defines the \emph{Green's function} $G(\gr,\grs)$, and we define the linear Green's operator
\begin{equation}
\field_{\sind} = \myop{G}(\param) \s = \langle G(\param),\s\rangle_{\bgdomain} = \gboundintvar{ G(\gr,\grs)\s(\grs) }{\grs}
\end{equation}
where $\rmd \grs$ is the surface integration measure on $\bgdomain$.

The definitions \eref{eq:gpde}-\eref{eq:gout} allow the definition of the mappings
\begin{eqnarray}
	&\textnormal{Transfer function} :& \tfn : \Q\rightarrow \Z \ , \quad \m = \tfn(\param)\s = \mbc\myop{G}(\param) {\s} \label{eq:transfer_function}\\
	&\textnormal{Forward Map} :& \fmap : \sspace \rightarrow \Z \ ,  \quad   \m = \fmap_{\sind}(\param) \label{eq:gfmap} \label{eq:forward_map}
\end{eqnarray}

The (non-linear) forward map is specified for a particular source pattern. Let $\myset{F}_{\myset{S}} = \{F_{\sind}, \s \in \myset{S}\}$ be the (possibly infinite) set of all forward mappings for the set of source currents $\myset{S}$.

Measured data is obtained as the result of a measurement operator
\begin{equation}
\operator{M} : \Z \rightarrow \dspace
\end{equation}
where $\dspace$ is a (possibly infinite) vector space.
Combining the forward map and measurement operator defines a forward operator
\begin{equation}
\mfmap = \operator{M}\fmap : \sspace \rightarrow \dspace
\end{equation}
Let $\myset{M}$ represent a (possibly infinite) set of aperture functions on $\bgdomain$ or $\exsurf$. Then a particular datum is obtained by applying a given aperture $\ap \in \myset{M}$ to the measurable
\begin{equation}
	\y_{\apind,\sind} = \operator{M}_{\apind}\m = \left<\ap,\m\right>_{\bgdomain} = \left<\ap,\fmap_{\sind}(\param)\right>_{\bgdomain}= \left<\ap,\mbc\myop{G}(\param) {\s}\right>_{\bgdomain}
\end{equation}

For the particular case of a $\delta$-function source we have
\begin{equation}
	\y_{\apind,\delta(\grs)} = \left<\ap,\mbc G(\grd,\grs)\right>_{\bgdomain}
\end{equation}
Similarly, if the measurements are considered simply point samples we can consider
\begin{equation}
	\y_{\delta(\grd),\delta(\grs)} = \mbc G(\grd,\grs)
\end{equation}

\subsection{Linearisation\label{sect:general_linearisation}}
The key tool for the study of the inverse problem is the linearisation of the forward map $\fmap$.
Consider a perturbation in the parameters $\param(\sr) \rightarrow \param(\sr) + \pparam(\sr)$. The \Frechet\, derivative of the forward mapping is the linear mapping defined by
\begin{equation}
\dffmap_{\sind}(\param) \pparam = \fmap_{\sind}(\param+\pparam) - \fmap_{\sind}(\param) + o(\|\pparam\|^2)
\label{eq:gfrechet}
\end{equation}
Let $\myset{F}'_{\myset{S}} = \{F_{\sind}, \s \in \myset{S}\}$ be the (possibly infinite) set of all \Frechet\, derivatives for the set of source currents $\myset{S}$.

When the \Frechet\, derivative is projected onto an aperture function we define the change in measurement  (for a particular perturbation $\pparam$)
\begin{equation}
	\dfy_{\apind,\sind} = \left<\ap,\dffmap_{\sind}(\param)\pparam\right>_{\bgdomain} = \left<\affmap_{\sind}(\param)\ap,\pparam\right>_{\gdomain}
\label{eq:dfy}
\end{equation}
which serves as the definition of the \emph{adjoint} \Frechet\, derivative. The adjoint \Frechet\, derivative will involve the solution of the adjoint problem
\begin{eqnarray}
	\aPDE(\param) \adjfield_{\apind} & =& 0 \label{eq:agpde}\\
	\ambc \adjfield_{\apind} &=& \ap \label{eq:agin}
\end{eqnarray}

It is convenient to define another notation for the \Frechet\, derivative as a general linear integral operator 
\begin{equation}
\dffmap_{\sind} \pparam = \myop{K}_{\sind}\pparam = \gdomintvar{ K_{\sind}(\grd,\gr) \pparam (\sr)}{\gr}
\end{equation}
with kernel $K_{\sind}$ and $\rmd \gr$ the volume integral measure on $\gdomain$. The notation $\myset{K}_{\myset{S}}$ is used for the complete set of these linear operators. 

>From the \Frechet\, derivative of the forward mapping we are lead to a definition of the \Frechet\, derivative of the Green's operator
\begin{equation}
\myop{K}_{\sind}\pparam = \mbc \myop{G}' [\s \otimes \pparam] = 
\left< \left<\mbc G',\pparam\right>_{\gdomain} , \Jin \right>_{\bgdomain}
\end{equation}
The operator $\myop{G}'$ has kernel $ G'(\grd,\gr,\grs)$ 
whose form depends on the particular problem considered. We may define in
general terms a potential operator
\begin{equation}
\ptlop(\pparam)=\PDE(\param+\pparam) - \PDE(\param)
\label{eq:ptlop_def}
\end{equation}
which leads to
\begin{equation}
\myop{G}' = -\myop{G}\ptlop\myop{G}\,,\quad\textnormal{and}\,\quad
\myop{G}^{'\ast} = -\myop{G}^{\ast}\ptlop\myop{G}
\end{equation}
These definitions will be made more precise in \sref{sect:analytic_approaches}.

For the particular case of $\delta$-function source and sampling measurements $\ap = \delta(\grd)$ the change in measurement given by \eref{eq:dfy} is simply the inner product of 
$G'$ with the perturbation
\begin{equation}
\dfy_{\delta(\grd),\delta(\grs)} =  \left<\mbc G'(\grd,\gr,\grs),\pparam(\sr)\right>_{\gdomain} =  \left<\pmdf,\pparam\right>_{\domain}
\end{equation}
which defines the \emph{sensitivity function} or \emph{measurement density function} $\pmdf(\sr)$. More generally,
the photon measurement density functions are the projection of the three-point Green's function onto a given source and measurement aperture
\begin{equation}
\pmdf_{\apind,\sind}(\sr) = \gboundintvar{ \ap(\grd)\gboundintvar{ \s(\grs) G'(\grd,\gr,\grs) }{\grs}}{\grd}
\label{eq:pmdf1}
\end{equation}

For the case of finite source and aperture functions, the forward 
mapping $\myset{F}_{\myset{S}}$ and all its \Frechet\, derivatives 
are continuous-to-discrete, and the adjoint mappings are discrete-to-continuous; 
such mappings necessarily contain an infinite dimensional null-space. When 
considering an infinite set of source and aperture functions, the forward 
mapping does not contain a null-space but is bounded and compact, leading to an 
inverse problem that is unbounded and ill-posed.

For computational purposes, the parameters of the inverse problem are considered in a basis representation
\begin{equation}
\param(\sr) = \sum_k \param_k \sbs_k(\sr)
\label{eq:sbs}
\end{equation}
This allows the definition of a discrete-to-continuous linear mapping
\begin{equation}
\dffmap_{\sind} \pparam = \myop{K}_{\sind}\pparam = \sum_k \pparam_k \gdomintvar{ K_{\sind}(\grd,\gr) \sbs_k (\sr)}{\gr}
\end{equation}
Considering a finite set of measurement functions of dimension $n_{\sind}$ and a finite set of aperture functions of dimension
$n_{\apind}$ leads to a vector of measurements of size $n_{\sind} \times n_{\apind}$. Taking the basis representation \eref{eq:sbs} leads to a discrete matrix representation of the linearised problem
\begin{equation}
\yv = \mA \fv
\label{eq:lin_discrete}
\end{equation}
where $\mymat{A}$ has matrix elements 
\begin{eqnarray}
\domintvar {\pmdf_{\apind,\sind}(\sr) \sbs_k (\sr)}{\sr} = \nonumber \\
\hspace{0.5cm}
\gboundintvar{ 
    \ap(\grd) \gboundintvar{
          \s(\grs)\gdomintvar{
              G'(\grd,\gr,\grs) \sbs_k (\sr)}{\gr}}{\grs}}{\grd}
\end{eqnarray}
This matrix is typically referred to as the \emph{Jacobian}, \emph{weight matrix} or \emph{sensitivity matrix} for the linearised inverse problem.

The adjoint operation
\begin{equation}
\fv = \mA\tr \yv
\end{equation}
can be considered as the discretisation of the adjoint operator
\begin{eqnarray}
\param &=& \afmfmap_{\sind} \y^{\delta} = \myop{K}^{\ast}_{\sind}\myop{M}^{\ast}\y^{\delta} = \gboundintvar{ K_{\sind}(\grd,\gr) \myop{M}^{\ast}\y(\rd)}{\grd}\label{eq:adjFrechet}\\
&=& \myop{G}^{'\ast} \ambc [\s \otimes \y] = 
\left< \left< G^{'\ast},\ambc\y\right>_{\bgdomain} , \Jin \right>_{\bgdomain}
\end{eqnarray}

\subsection{Inverse problems
\label{sect:general_inverse}}

Our goal is to invert the forward map $\fmap$ and recover the parameters $x$.
To proceed, we note that
the adjoint \Frechet\,derivative operator defined in \eref{eq:adjFrechet} acts 
on an element in data space to give an element in parameter space representing 
a change in $\param$ that reduces the norm of the residual difference between 
the data and the forward operator acting on $\param$. In terms of optimisation 
theory this direction is the  steepest descent  update direction for the 
inverse problem in the sense that it gives the direction in which this residual 
norm decreases locally most rapidly. This operator is therefore the basis for
gradient based optimisation methods. We explore this idea in more 
detail in \sref{sect:linear_methods}.


Solution methods for inverse problems typically perform better if making used of second derivative information. The second \Frechet\,derivative is given by
\begin{equation}
\ddffmap_{\sind}(\param) (\pparam_1,\pparam_2) = \fmap_{\sind}(\param+\pparam_1 + \pparam_2) - \dffmap_{\sind}(\param) (\pparam_1+\pparam_2) - \fmap_{\sind}(\param)
\label{eq:gfrechet2}
\end{equation}
with corresponding second derivative Green's operator $\myop{G}''$ and kernel $G''$, given by
\begin{equation}
\myop{G}'' = \myop{G}\ptlop \myop{G}\ptlop\myop{G}
\end{equation}
and in general we can define the nonlinear mapping $\fmap$ as
\begin{equation}
\fmap(\param+\pparam) = \mbc\left[\myop{G}(\param) - \myop{G}(\param)\ptlop(\pparam)\myop{G}(\param) + \cdots + (-1)^n\myop{G}(\param)\left(\ptlop(\pparam) \myop{G}(\param)\right)^{n} + \cdots \right]
\label{eq:born_general}
\end{equation}
which is the Born series. We discuss this series in more detail in~\eref{Born_series} and \eref{pert_V} in section~\ref{sec:physics}.

\begin{figure}[t]
\hspace{1.5cm}
	\includegraphics[width=0.22\textwidth]{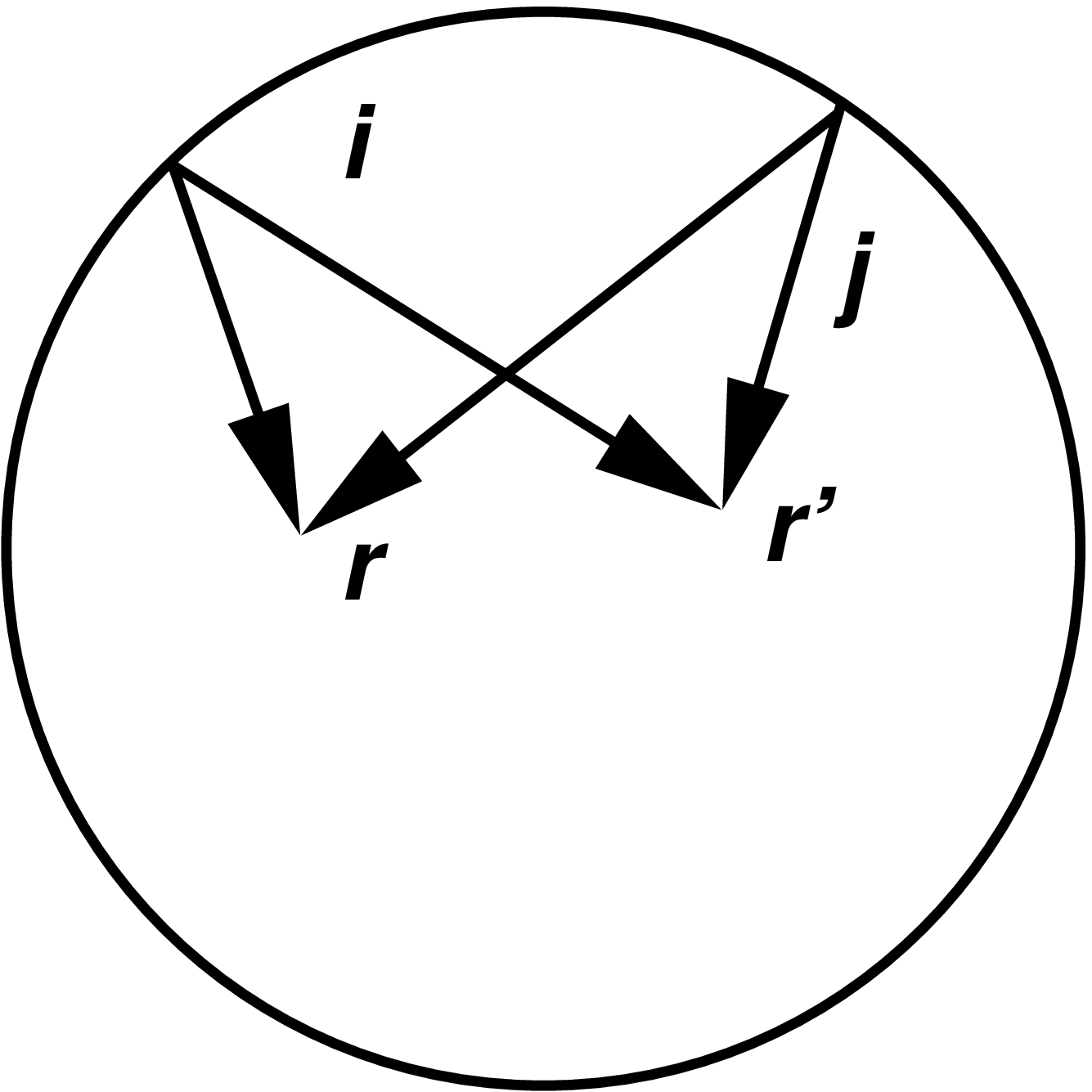}
\hspace{1em}
	\includegraphics[width=0.22\textwidth]{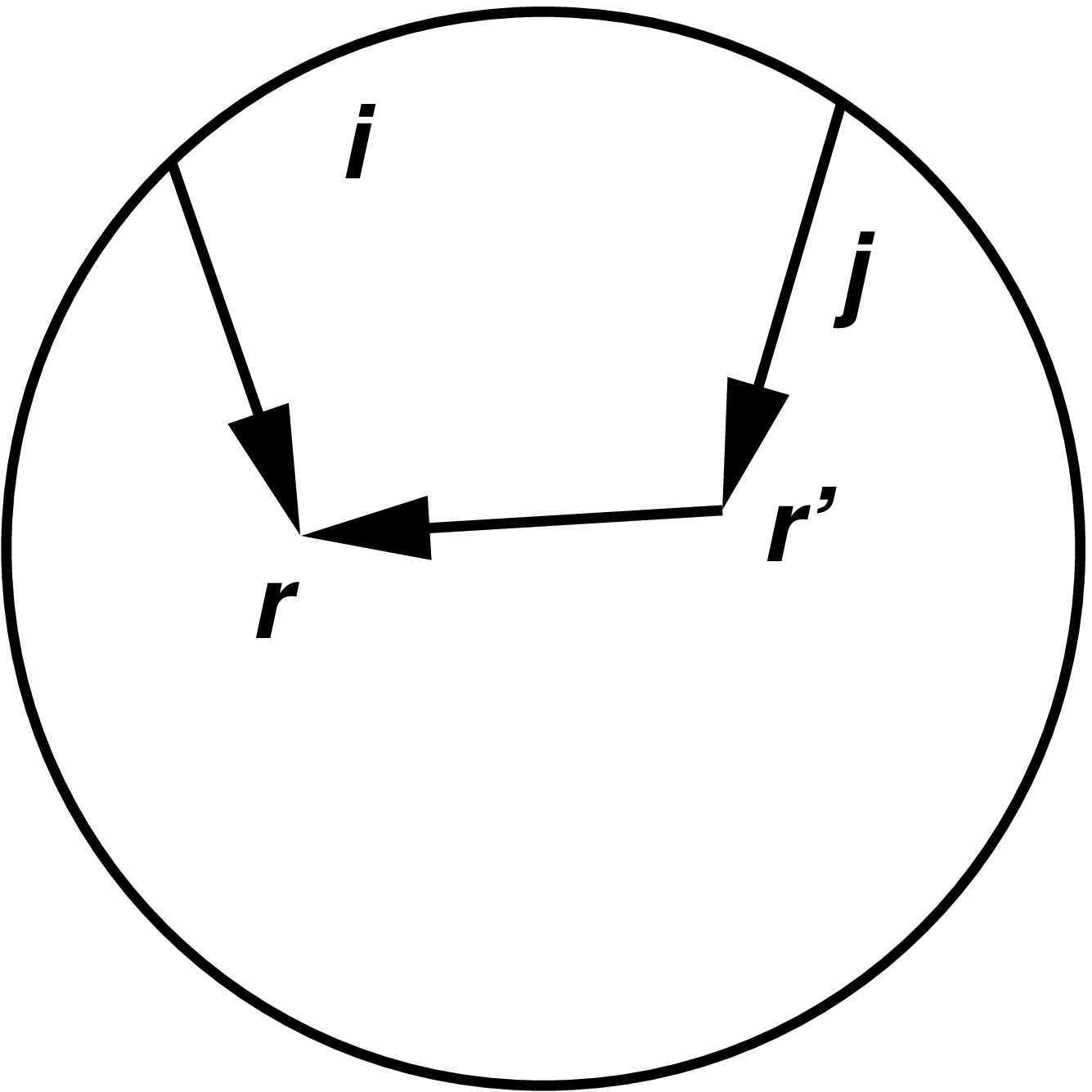}
\hspace{1em}
	\includegraphics[width=0.22\textwidth]{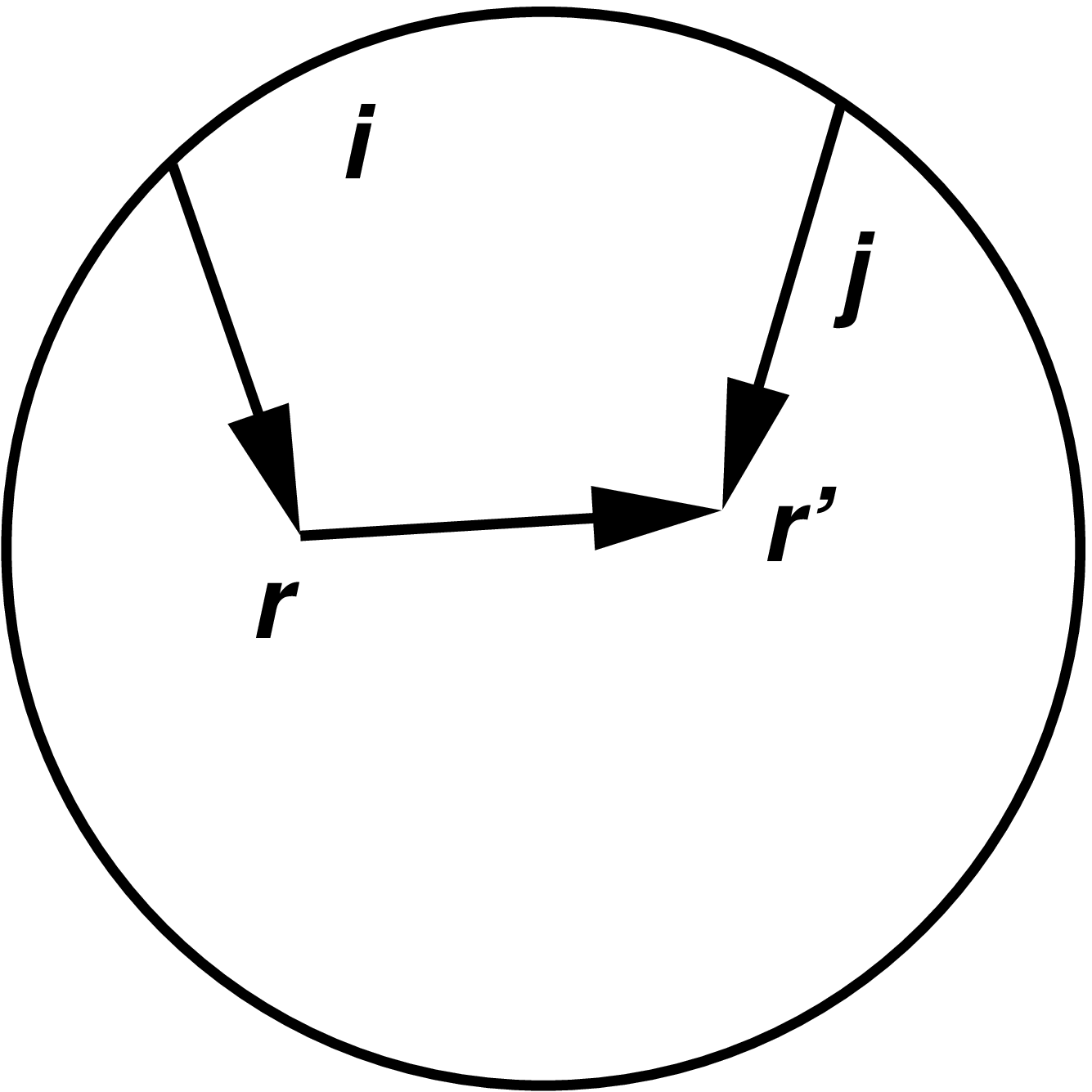}
\caption{Second order derivative operators. On the left is the kernel of
 $\myop{K}^{\ast}\myop{K}$
In the centre the term $ \myop{G}^{\ast}\ptlop \myop{G}\ptlop\myop{G} $. On the right the term $\myop{G}^{\ast}\ptlop \myop{G}^{\ast}\ptlop\myop{G}$.
\label{fig:fourgreen}
}
\end{figure}

For optimisation methods the second \Frechet\,derivative gives rise to a Hessian operator
\begin{equation}
\myop{H} = \myop{K}^\ast \myop{K} - \aafmfmap\y^{\delta} 
\end{equation}
where the adjoint operator $\aafmfmap$ is composed of two terms
\begin{equation}
\myop{G}^{''\ast} = \myop{G}^{\ast}\ptlop \myop{G}\ptlop\myop{G} + \myop{G}^{\ast}\ptlop \myop{G}^{\ast}\ptlop\myop{G}
\end{equation}
The second derivative Green's kernels are represented in~\fref{fig:fourgreen}.

Taking the basis representation as in \eref{eq:lin_discrete}
leads to a discrete Hessian representation
\begin{equation}
\mymat{H} = \mymat{A}\tr\mymat{A} - \mymat{F}^{''\mathrm{T}}  \yv
\end{equation}
which forms the basis of the Newton update scheme
\begin{equation}
\mymat{H} \fv = \mymat{A}\tr \yv
\end{equation}
This is discussed in further detail in \sref{sect:Optimisation}.

\section{Waves, Transport and Diffusion}
The mathematical description of light propagation in random media changes according to the length scale of interest~\cite{vanrossum_1999_1}. We begin with the Maxwell equations, which are valid on microscopic scales. The mesoscale, in which the characteristic scale is set by the scattering length, is described by the radiative transport equation (RTE). Finally, we discuss the macroscale, which is described by the diffusion approximation to the RTE.

\subsection{Electromagnetic waves}
The Maxwell equations together with the Lorentz force law govern the classical description of all electromagnetic phenomena. For simplicity, we restrict our attention to nonmagnetic media and to monochromatic fields. We employ the convention that physical quantities are given by the real part of the corresponding complex quantities. In the absence of sources, the electric field $\BE$ in an inhomogeneous medium with a position-dependent permittivity $\varepsilon$ satisfies the time-independent wave equation
\begin{equation}
\label{wave_E}
\curl\curl \BE(\Br) - k_0^2 \varepsilon(\Br) \BE(\Br) = 0 \ ,
\end{equation}
where $k_0$ is the free-space wavenumber.
In addition, the divergence condition 
\begin{equation}
\label{div_condition}
\div \varepsilon(\Br)\BE(\Br) = 0 
\end{equation}
must be obeyed. If $\varepsilon$ varies slowly on the scale of the wavelength, then $\BE$ 
satisfies the scalar wave equation
\begin{equation}
\label{scalar_wave}
\lap U(\Br) + k_0^2 \varepsilon(\Br) U(\Br) = 0 \ ,
\end{equation}
where $U$ is any of the components of $\BE$. Note that the components of $\BE$ are still coupled according to~(\ref{div_condition}).
If we restrict our attention to (\ref{scalar_wave}) we will say that we are working within the scalar theory of electromagnetic fields.

The conservation of energy is governed by the relation
\begin{equation}
\label{energy_cons}
\div \BJ + \frac{4\pi k_0}{c} \Im\left(\varepsilon\right) I = 0 \ ,
\end{equation}
where the energy current density $J$ is defined by 
\begin{equation}
\label{scalar_current}
\BJ = \frac{1}{2ik_0}\left(U^*\grad U - 
U\grad U^* \right) \ 
\end{equation}
and the intensity $I$ is given by
\begin{equation}
I = \frac{c}{4\pi}|U|^2 \ .
\end{equation}
Note that $\BJ$ plays the role of the Poynting vector in the scalar theory.

The solution to (\ref{scalar_wave}), which obeys the Sommerfeld radiation condition, is given by the Lippmann-Schwinger equation
\begin{equation}
\label{Lippmann-Schwinger}
U(\Br) = U_i(\Br) + k_0^2\int d\Br' G(\Br,\Br')U(\Br')\eta(\Br') \ ,
\end{equation}
where $\eta=(\varepsilon-1)/4\pi$ is the dielectric susceptibility.
Here $U_i$ is incident field, which obeys (\ref{scalar_wave}) with $\eta=0$ and $G$ is the Green's function which, in free space, is given by
\begin{equation}
\label{scalar_G_0}
G_{\infty}(\Br,\Br') = \frac{e^{ik_0|\Br-\Br'|}}{|\Br-\Br'|} \ .
\end{equation}
If we iterate (\ref{Lippmann-Schwinger}) starting from $U=U_i$,  we obtain an  infinite series for $U$ of the form
\begin{eqnarray}
\label{Born_series}
\nonumber
U(\Br)  = U_i(\Br) + k_0^2\int d\Br' G(\Br,\Br')\eta(\Br')U_i(\Br') \\
+ k_0^4 \int d\Br' d\Br'' G(\Br,\Br')\eta(\Br')G(\Br',\Br'')\eta(\Br'')U_i(\Br'') + \cdots  
\ . 
\end{eqnarray}
Eq.~(\ref{Born_series}) is known as the Born series and each of its terms corresponds to successively higher orders of scattering of the incident field. If only the first term in the series is retained, this is known as the Born approximation.  
 
\subsection{Radiative transport}

In radiative transport theory, the propagation of light through a material medium 
is formulated in terms of a conservation law that accounts for gains and losses of photons due to scattering and absorption~\cite{Ishimaru_book,Case_book}. The fundamental quantity of interest is the specific intensity $I(\Br,\svec,t)$, defined as the intensity at the position $\Br$ in the direction $\svec$ at time $t$. The specific intensity obeys the radiative transport equation (RTE):
\begin{eqnarray}
\label{RTE_general}
\frac{1}{c} \dfrac{\partial I}{\partial t} + \svec\cdot\grad I + (\mu_a+\mu_s)I= \mu_s\int p(\svec',\svec)I(\Br,\svec')d\svec' \ , \ \ \ \Br\in \Omega \ ,
\end{eqnarray}
where $\mu_a$ and $\mu_s$ are the absorption and scattering coefficients. The speciÞc intensity also satisÞes the half-range boundary condition 
\begin{equation}
\label{RTE_bc_ihg}
I({\bm r},\hat{\bf s}) 
=  I_{\rm inc}({\bm r},\hat{\bf s}) \ , \ \ \hat{\bf s} \cdot \n < 0 \ , \ \ \Br \in \partial \Omega \ ,
\end{equation}
where $\n$ is the outward unit normal to $\partial\Omega$ and
$I_{\rm inc}$ is the incident specific intensity at the boundary. The above choice of boundary condition guarantees the uniqueness of solutions to the RTE~\cite{Case_book}.
The phase function $p$ is symmetric with respect to interchange of its arguments and obeys the normalization condition
\begin{equation}
\label{normalization_p}
\int p(\svec,\svec') d\svec'  = 1 \ ,
\end{equation}
for all $\svec$. We will often assume that $p(\svec,\svec')$ depends only upon the angle between $\svec$
and $\svec'$, which holds for scattering by spherically-symmetric particles. Note that the choice $p=1/(4\pi)$ corresponds to isotropic scattering. For scattering that is strongly peaked in the forward direction ($\svec\cdot\svec'\approx 1$), an asymptotic expansion of the right hand side of~\eref{RTE_general} may be performed~\cite{Leakeas_2001_1}. This leads to the Fokker-Planck form of the RTE
\begin{equation}
\label{RTE_FP}
\frac{1}{c} \frac{\partial I}{\partial t} + \svec\cdot\grad I + (\mu_a+\mu_s)I= \mu_s L I \ , \ \ \ L=
\half (1-g) \Delta_{\svec} \ ,
\end{equation}
where $\Delta_{\svec}$ is the Laplacian on the two-dimensional unit sphere $S^2$ and $g\approx 1$ is the anisotropy of scattering, as given by \eref{eq:gdef} in \sref{sec:diffuse}. A rational approximation to the transport term in which the operator $L$ is given by
\begin{equation}
\label{LL}
L= a \Delta_{\svec}({\rm 1}-b \Delta_{\svec})^{-1} 
\end{equation}
has also been proposed~\cite{Leakeas_2001_1,Gonzalez-Rodriguez_2008_1}. Here the constants $a,b$ can be computed from angular moments of the phase function. In the case of the Henyey-Greenstein phase function
\begin{equation}
p(\svec\cdot\svec')= \frac{1}{4\pi}\sum_{l=0}^\infty (2l+1) g^l P_l(\svec\cdot\svec') \ ,
\end{equation}
it can be seen that 
\begin{equation}
a= \half (1-g)(1+2b)\ , \ \ \ b= \frac{2-3g+g^2}{6g(1-g)} \ .
\end{equation}
Note that the form of $L$ in~\eref{LL} is a bounded operator on $L^2(S^2)$, which is not the case for the spherical Laplacian that appears in~\eref{RTE_FP}.

The total power $P$ passing through a surface $\Sigma$ is related to the specific intensity by
\begin{equation}
P = \int_{\Sigma} d\Br \int d\svec I(\Br,\svec,t) \svec \cdot\n  \ .
\end{equation}
The energy density $\Phi$ is obtained by integrating out the angular dependence of the specific intensity:
\begin{equation}
\label{def_u}
\Phi(\Br,t) = \frac{1}{c}\int I(\Br,\svec,t) d\svec \ .
\end{equation}
We note that the RTE allows for the addition of intensities. As a result, it cannot explain certain wavelike phenomena.

\subsubsection{From waves to transport}
The RTE can be derived by considering the high-frequency asymptotics of wave propagation in a random medium. We briefly recall the main ideas in the context of monochromatic scalar waves. The general theory for vector electromagnetic waves is presented in~\cite{Ryzhik_1996_1}. We assume that the random medium is statistically homogeneous and that the susceptibility $\eta$ is a Gaussian random field such that
\begin{equation}
\langle \eta(\Br) \rangle = 0 \ , \ \ \
\langle \eta(\Br)\eta(\Br')\rangle = C(|\Br-\Br'|) \ ,
\end{equation}
where $C$ is the two-point correlation function and $\langle \cdots \rangle$ denotes statistical averaging. Let $L$ denote the propagation distance of the wave. At high frequencies, $L$ is large compared to the wavelength and we introduce a small parameter $\epsilon = 1/(k_0 L) \ll 1$. We suppose that the fluctuations in $\eta$ are weak so that $C$ is of the order $O(\epsilon)$. We then rescale the spatial variable according to $\Br\to \Br/\epsilon$ and define the scaled field $U_\epsilon(\Br)=U(\Br/\epsilon)$, so that (\ref{scalar_wave}) becomes
\begin{equation}
\label{scalar_wave_scaled}
\epsilon^2\lap U_{\epsilon}(\Br) + U_{\epsilon}(\Br) = -4\pi \sqrt{\epsilon} \eta\left(\Br/\epsilon\right) U_{\epsilon}(\Br)\ .
\end{equation}
Here we have introduced a rescaling of $\eta$ to be consistent with the assumption that the fluctuations are of strength $O(\epsilon)$.

Although (\ref{energy_cons}) gives some indication of how the intensity of the field is distributed in space, it does not prescribe how the intensity propagates. To overcome this difficulty, we introduce the Wigner distribution $W_\epsilon(\Br,\Bk)$, which is a function of the position $\Br$ and the wave vector $\Bk$:
\begin{equation}
\label{wigner}
W_\epsilon(\Br,\Bk)=\int d \BR e^{i \Bk \cdot \BR} U_\epsilon\left(\Br-\half \epsilon \BR\right)
U_\epsilon^*\left(\Br+\half \epsilon \BR\right)  \ .
\end{equation}
The Wigner distribution has several important properties. It is real-valued and is related to the intensity and energy current density by the formulas
\begin{eqnarray}
I = \frac{c}{4\pi}\int \frac{d\Bk}{(2\pi)^3} W_\epsilon(\Br,\Bk)  \ , \ \ \ \BJ = \int \frac{d\Bk }{(2\pi)^3} \Bk W_\epsilon(\Br,\Bk) \ .
\end{eqnarray}
Making use of (\ref{scalar_wave_scaled}), it can be seen that the Wigner distribution obeys the equation
\begin{equation}
\Bk \cdot \grad_{\Br} W_\epsilon + i \frac{2\pi}{\sqrt{\epsilon}} \int d\Bq e^{-i \Bq \cdot x/\epsilon} \tilde\eta(\Bq) \left(W_\epsilon(\Br,\Bk+\half \Bq) -W_\epsilon(\Br,\Bk-\half \Bq)  \right) = 0 \ ,
\end{equation}
where we have assumed that $\eta$ is real-valued and $\tilde\eta$ denotes the Fourier transform of $\eta$ which is defined by
\begin{equation}
\tilde\eta(\Bq) = \int d\Br e^{i \Bq \cdot \Br} \eta(\Br) \ .
\end{equation}

We now consider the asymptotics of the Wigner function in the homogenization limit $\epsilon \to 0$. This corresponds to the regime of high-frequencies and weak fluctuations. We proceed by introducing a two-scale expansion for $W_\epsilon$ of the form
\begin{equation}
W_\epsilon(\Br,\Br',\Bk) = W_0(\Br,\Bk) + \sqrt{\epsilon} W_1(\Br,\Br',\Bk) + \epsilon W_2(\Br,\Br',\Bk) + \cdots \ ,
\end{equation}
where $\Br'=\Br/\epsilon$ is a fast variable. By averaging over the fluctuations on the fast scale,
it is possible to show that $\langle W_0 \rangle$, which we denote by $W$, obeys the equation
\begin{equation}
\label{transport_eqn}
\Bk \cdot \grad_{\Br} W = \int d\Bk' \tilde C(\Bk-\Bk')\delta(k^2-k'^2)\left( W(\Br,\Bk') -W(\Br,\Bk)\right) \ .
\end{equation}
Evidently, (\ref{transport_eqn}) has the form of a time-independent transport equation. The role of the delta function is to conserve momentum, making it possible to view $W$ as a function of position and the direction $\Bk/|\Bk|$. We note that the phase function and scattering coefficient are related to statistical properties of the random medium, as reflected in the appearance of the correlation function $C$ in~\eref{transport_eqn}. If the medium is composed of spatially uncorrelated point particles with number density $\rho$, then
\begin{equation}
\mu_a = \rho\sigma_a \ , \quad \mu_s = \rho\sigma_s \ , \quad 
p= \frac{d\sigma_s}{d\Omega}\big/\sigma_s \ ,
\end{equation}
where $\sigma_a$ and $\sigma_s$ are the absorption and scattering cross sections of the 
particles and $d\sigma_s/d\Omega$ is the differential scattering cross section. Note that 
$\sigma_a$, $\sigma_s$ and $p$ are wavelength dependent quantities.

The above derivation of the RTE can be modified when the speed of light is not constant. Following Ref.~\cite{Bal_2006_1}, the statistically-averaged Wigner function can be seen to obey the equation
\begin{equation}
\label{RTE_variable_n}
n \Bk \cdot \grad_{\Br} W -k_0\grad n \cdot \left(I- \svec\otimes\svec)\right)\grad_{\svec} W = \int d\Bk' \tilde C(\Bk-\Bk')\delta(k^2-k'^2)\left( W(\Br,\Bk') -W(\Br,\Bk)\right) \ ,
\end{equation}
where $n$ is the index of refraction. Using the relation $I=c(k_0/n(\Br))^3 W$ and the above result, we obtain the required generalization of the RTE.

Radiative transport theory can be generalized to account for energy transport for different states of polarization~\cite{Clark_2009_1,Ryzhik_1996_1}. The RTE is then formulated in terms of the 
Stokes vector ${\bf I}=(I_\parallel, I_\perp, U, V)$. Here the Stokes parameters, which characterize the polarization of the field,  are defined as
\begin{eqnarray}
I_\parallel = \langle E_\parallel E_\parallel^* \rangle \ , \quad I_\perp = \langle E_\perp E_\perp^* \rangle \ , \\
U = \langle E_\parallel E_\perp^* + E_\parallel^* E_\perp^\rangle \ , \ \quad V = i(\langle E_\parallel E_\perp^* - E_\parallel^* E_\perp \rangle) \ ,
\end{eqnarray}
where $E_\parallel, E_\perp$ denote two orthogonal components of the electric field in the plane perpendicular to the direction of propagation. Note that the total intensity $I=I_\parallel + I_\perp$. The Stokes vector satisfies the generalized RTE 
\begin{equation}
\frac{1}{c} \frac{\partial {\bf I}}{\partial t} + \svec\cdot\grad {\bf I} + (\mu_a+\mu_s){\bf I}= \mu_s\int M(\svec',\svec) \cdot {\bf I}(\Br,\svec')d\svec' \ , 
\end{equation}
where $M$ denotes the Mueller matrix, which accounts for mixing of different states of polarization upon scattering.

\subsection{Collision expansion}

The RTE (\ref{RTE_general}), obeying the boundary condition (\ref{RTE_bc_ihg}), is equivalent to the integral equation
\begin{equation}
\label{RTE_int}
I({\bm r}, \hat{\bf s}) = I_0({\bm r}, \hat{\bf s}) + \int G_0({\bm r},\hat{\bf s}; {\bm r}^{\prime},\hat{\bf s}^{\prime}) \mu_s({\bm r}^{\prime}) p(\hat{\bf s}^{\prime},\hat{\bf s}^{\prime\prime})
I({\bm r}^{\prime}, \hat{\bf s}^{\prime\prime}) d\Br^{\prime} d
{\svec}^{\prime}d{\svec}^{\prime\prime} \ .
\end{equation}
Here $I_0$ is the unscattered (ballistic) specific
intensity, which satisfies the equation
\begin{equation}
\label{RTE_ballistic}
\left[ \hat{\bf s} \cdot \nabla + \mu_a + \mu_s \right] 
I_0 = 0 \ , 
\end{equation}
and $G_0$ is the ballistic Green's function
\begin{equation}
\label{G_b}
G_0({\bm r},\hat{\bf s}; {\bm r}^{\prime},\hat{\bf s}^{\prime}) = 
g({\bm r}, {\bm r}^{\prime})
\delta\left(\hat{\bf
    s}^{\prime} - \frac{\Br-\Br'}{|\Br-\Br'|}\right) \delta(\hat{\bf s} - \hat{\bf s}^{\prime}) \ , 
\end{equation}
where
\begin{equation}
\label{g}
g({\bm r}, {\bm r}^{\prime}) = 
\frac{1}{\vert {\bm r} - {\bm r}^{\prime}
\vert^2} \exp\left[ - \int_0^{\vert {\bm r} - {\bm r}^{\prime}\vert}
\mu_t\left({\bm r}^{\prime} + \ell \dfrac{\Br-\Br'}{|\Br-\Br'|}\right)d\ell \right] \ ,
\end{equation}
and the extinction coefficient $\mu_t = \mu_a + \mu_s$. 
Note that if a narrow collimated beam of intensity $I_{\rm inc}$ is incident on the medium at the point ${\bm r}_0$ in the direction $\hat{\bf s}_0$, then $I_0({\bm r},
\hat{\bf s})$ is given by
\begin{equation}
\label{I_b}
I_0({\bm r}, \hat{\bf s}) = I_{\rm inc} G_0({\bm r},\hat{\bf s}; {\bm r}_0,\hat{\bf s}_0) \ ,
\end{equation}

To derive the collision expansion, we iterate (\ref{RTE_int}) starting from $I^{(0)}=I_0$ and obtain
\begin{equation}
\label{series}
I({\bm r}, \hat{\bf s})= I^{(0)}({\bm r}, \hat{\bf s})+ I^{(1)}({\bf
r}, \hat{\bf s})+ I^{(2)}({\bm r}, \hat{\bf s})+\cdots \ ,
\end{equation}
where  each term of the series is given by
\begin{equation}
I^{(n)}({\bm r}, \hat{\bf s})= \int d\Br^{\prime} d
{\svec}^{\prime}d{\svec}^{\prime\prime} G_0({\bm r},\hat{\bf s}; {\bm
  r}^{\prime},\hat{\bf s}^{\prime}) \mu_s({\bm r}^{\prime}) p(\hat{\bf
  s}^{\prime},\hat{\bf s}^{\prime\prime}) I^{(n-1)}({\bm r}^{\prime},
\hat{\bf s}^{\prime\prime}) \ ,
\end{equation}
with $n=1,2,\dots$.
The above series is the analog of the Born series
for the RTE, since each term accounts for successively higher orders of scattering.

It is instructive to examine the expression for $I^{(1)}$, which is the contribution
to the specific intensity from single scattering: 
\begin{equation}
\label{I_s}
I^{(1)}({\bm r}, \hat{\bf s}) = \int d\Br^{\prime} d
{\svec}^{\prime}d{\svec}^{\prime\prime}G_0({\bm
  r},\hat{\bf s}; {\bm r}^{\prime},\hat{\bf s}^{\prime}) \mu_s({\bm
  r}^{\prime}) p(\hat{\bf s}^{\prime},\hat{\bf s}^{\prime\prime})
I_0({\bm r}^{\prime}, \hat{\bf s}^{\prime\prime})  \ .
\end{equation}
It follows from~\eref{I_b} that the change in intensity $\delta I=I-I^{(0)}$ measured by a point detector at $\Br_2$ in the direction $\svec_2$, due to a unit-amplitude point source at $\Br_1$ in the direction $\svec_1$ is given by
\begin{eqnarray}
\label{delta_integrated}
\nonumber
\delta I(\Br_1,\svec_1;\Br_2,\svec_2)&=&p(\svec_1,\svec_2)\int_0^\infty dR R^2 g(\Br_2,\Br_1 + R \svec_1)g(\Br_1 + R \svec_1,\Br_1) \\
&&\times\delta\left(\dfrac{\Br_2-\Br_1 - R\svec_1}{|\Br_2-\Br_1 - R\svec_1|}-\svec_2 \right)\mu_s(\Br_1+R\svec_1) \ .
\end{eqnarray}
Suppose that $\Br_1$ and $\Br_2$ are located on the boundary of a bounded domain and $\svec_1$ points into and $\svec_2$ points out of the domain. Then the rays in the directions $\svec_1$ and $\svec_2$ must intersect at a point $\BR$ that lies in the interior of the domain. In addition, the delta function in (\ref{delta_integrated}) implements the constraint that $\svec_1$, $\svec_2$ and $\Br_1-\Br_2$ all lie in the same plane. Using this fact and carrying out the integration in (\ref{delta_integrated}), we find that
\begin{eqnarray}
\label{delta_final}
\delta I(\Br_1,\svec_1;\Br_2,\svec_2)\propto\exp\left[-\int_0^{L_1} \mu_t(\Br_1+\ell \svec_1)d\ell-\int_0^{L_2} \mu_t(\BR +\ell\svec_2)d\ell \right] \ ,
\end{eqnarray}
where $L_{1}=|\BR-\Br_{1}|$, $L_2=|\BR-\Br_2|$ and we have omitted overall geometric prefactors. Note that the argument of the exponential corresponds to the integral of $\mu_t$ along the broken ray which begins at $\Br_1$, passes through $\BR$, and terminates at $\Br_2$. The significance of such broken rays will be discussed in section~\ref{sec:SSOT}.

The terms in the collision expansion can be classified by their smoothness. The lowest order term is the most singular. In the absence of scattering, according to~\eref{g}, this term leads to a Radon transform relationship between the absorption coefficient and the specific intensity, under that condition that the source and detector are collinear. Inversion of the Radon transform is the basis for optical projection tomography~\cite{Sharpe_2002_1,vinegoni_2008_1}.
The first order term is also singular, as is evident from the presence of a delta function in~\eref{delta_integrated}. Terms of higher order are of increasing smoothness. This observation has been exploited to prove uniqueness of the inverse transport problem and to study its stability. A comprehensive review is presented in~\cite{bal2009a}.

The above discussion has implicitly assumed that the angular dependence of the specific intensity is measurable. In practice, such measurements are extremely difficult to obtain. The experimentally measurable intensity is often an angular average of the specific intensity over the aperture of an optical system. The effect of averaging is to remove the singularities that are present in the specific intensity. The resulting inverse problem is then highly ill-posed~\cite{Bal_2008_1}.
To illustrate this point, we observe that it follows from~\eref{def_u} and~\eref{RTE_int} that the energy density $\Phi$ obeys the integral equation
\begin{eqnarray}
\label{averaged_I}
\Phi(\Br) = \Phi^{(0)}(\Br) + \frac{1}{4\pi}\int d\Br' \frac{e^{-\mu_t|\Br-\Br'|}}{|\Br-\Br'|^2}p(\Br')\Phi(\Br') \ ,
\end{eqnarray}
where the scattering is assumed to be isotropic and
\begin{equation}
\Phi^{(0)}(\Br) = \frac{1}{c}\int I^{(0)}(\Br,\svec)d\svec \ .
\end{equation}
Here the kernel appearing in~\eref{averaged_I} is smoothing with a Fourier transform that decays algebraically at high frequencies. 

\subsection{Transport regime\label{sect:Physics_transport_regime}}

We now develop the scattering theory for the RTE in an inhomogeneously absorbing medium. We proceed by decomposing $\mu_a$ into a constant part $\bar\mu_a$ and a spatially varying part $\delta\mu_a$:
\begin{equation}
\label{def_delta_mu_a}
\mu_a(\Br)=\bar\mu_a + \delta\mu_a(\Br) \ .
\end{equation}
The stationary form of the RTE~\eref{RTE_general} can be rewritten in the form
\begin{equation}
\label{RTE_delta}
\svec\cdot\grad I + \bar\mu_tI -\mu_s\int p(\svec',\svec)I(\Br,\svec') d\svec' = -\delta\mu_a(\Br)I \ ,
\end{equation}
where $\bar\mu_t=\bar\mu_a+\mu_s$.
The solution to (\ref{RTE_delta}) is given by
\begin{equation}
\label{LS}
I(\Br,\svec)=I_i(\Br,\svec) - \int d\Br' d\svec' G(\Br,\svec;\Br',\svec')\delta\mu_a(\Br')I(\Br',\svec') \ ,
\end{equation}
where $G$ denotes the Green's function for a homogeneous medium with absorption $\bar\mu_a$ and $I_i$ is the incident specific intensity. Eq.~(\ref{LS}) is the analog of the Lippmann-Schwinger equation for the RTE. It describes the ``multiple scattering'' of the incident specific intensity from inhomogeneities in $\delta\mu_a$. If only one absorption event is considered, then the intensity $I$ on the right hand side of (\ref{LS}) can be replaced by the incident intensity $I_i$. This result describes the linearization of the integral equation (\ref{LS}) with respect to $\delta\mu_a$. If the incident field is generated by a point source at $\Br_1$ pointing in the direction $\svec_1$, then the change in specific intensity due to spatial fluctuations in absorption $\delta I$ can be obtained from the relation
\begin{equation}
\label{phi}
\delta I(\Br_1,\svec_1;\Br_2,\svec_2)= I_0\int d\Br d\svec G(\Br_1,\svec_1;\Br,\svec)G(\Br,\svec;\Br_2,\svec_2)\delta\mu_a(\Br) \ .
\end{equation}
Here $I_0$ denotes the intensity of the source and $\Br_2,\svec_2$ are the position and orientation of a point detector.

To make further progress requires knowledge of the Green's function for the RTE. The Green's function for an infinite medium with isotropic scattering can be obtained explicitly~\cite{Case_book}. In three dimensions we have
\begin{eqnarray}
\nonumber
\fl
G(\Br,\svec;\Br',\svec')=G_0(\Br,\svec;\Br',\svec') 
+\frac{\mu_s}{4\pi} \int &&\frac{d\Bk}{(2\pi)^3}e^{i\Bk\cdot(\Br-\Br')}  \frac{1}{(\mu_t +
i \svec\cdot\Bk)(\mu_t + i \svec'\cdot\Bk)} \\
&& \times \frac{1}{1-\frac{\mu_s}{|\Bk|}
\tan^{-1}\left(\frac{|\Bk|}{\mu_t}\right)} \ ,
\end{eqnarray}
where $G_0$ is the ballistic Green's function which is defined
by~\eref{G_b}. More generally, numerical procedures such as the discrete ordinate method or the $P_N$ approximation may be employed to compute the Green's function in a bounded domain with anisotropic scattering, as described in section~\ref{sect:modelling_methods}.

The method of rotated reference frames is a spectral method for the computing the Green's function for the three-dimensional RTE in a homogeneous medium with anisotropic scattering and planar boundaries~\cite{Panasyuk_2006_1}. It is derived by considering the plane-wave modes for the RTE which are of the form
\begin{equation}
\label{modes}
I(\Br,\svec)=A_\Bk(\svec)e^{\Bk\cdot\Br} \ ,
\end{equation}
where the amplitude $A$ is to be determined. Evidently, the components of the wave vector $\Bk$ cannot be purely real; otherwise the above modes would have exponential growth in the $\hat\Bk$ direction. We thus consider evanescent modes with
\begin{equation}
\Bk = i\Bq \pm \sqrt{q^2 + 1/\lambda^2} \ \z \ ,
\end{equation}
where $\Bq\cdot\z=0$ and $\Bk\cdot\Bk\equiv 1/\lambda^2$. These modes are oscillatory in the transverse direction and decay in the $\pm z$-directions. By inserting (\ref{modes}) into the RTE (\ref{RTE_general}), we find that
$A_\Bk(\svec)$ satisfies the equation
\begin{equation}
\label{eigenproblem}
\left(\svec\cdot\Bk + \mu_a + \mu_s \right)A_\Bk(\svec) =\mu_s\int  p(\svec,\svec')A_\Bk(\svec')d\svec' \ .
\end{equation}

To solve the eigenproblem defined by (\ref{eigenproblem}) it will prove useful to expand $A_\Bk(\svec)$ into a basis of spherical functions defined in a  rotated reference frame whose $z$-axis coincides with the direction $\hat\Bk$.
We denote such functions by $Y_{lm}(\svec;\hat\Bk)$ and define them via the relation
\begin{equation}
Y_{lm}(\svec;\hat\Bk)=\sum_{m'=-l}^l D_{mm'}^l(\varphi,\theta,0)Y_{lm'}(\svec) \ ,
\end{equation}
where $Y_{lm}(\svec)$ are the spherical harmonics defined in the laboratory frame, $D_{mm'}^l$ is the Wigner $D$-function and $\varphi,\theta$ are the polar angles of $\hat\Bk$ in the laboratory frame. 
We thus expand $A_\Bk$ as
\begin{equation}
\label{A_expansion}
A_\Bk(\svec)=\sum_{l,m}C_{lm}Y_{lm}(\svec;\hat\Bk) \ ,
\end{equation}
where the coefficients $C_{lm}$ are to be determined.
Note that since the phase function $p(\svec,\svec')$ is invariant under simultaneous rotation of $\svec$ and $\svec'$, it may be expanded into rotated spherical functions according to
\begin{equation}
p(\svec,\svec')=\sum_{l,m}p_lY_{lm}(\svec;\hat\Bk)Y_{lm}^*(\svec';\hat\Bk) \ ,
\end{equation}
where the expansion coefficients $p_l$ are independent of $\hat\Bk$. An alternative approach that may be used to solve the eigenproblem~\eref{eigenproblem} is to employ the method of discrete ordinates~\cite{Kim_2003_1,Kim_2004_1,Kim_2006_1}.

Substituting (\ref{A_expansion}) into (\ref{eigenproblem}) and making use of the orthogonality properties of the spherical functions, we find that the coefficients $C_{lm}$ satisfy the equation
\begin{equation}
\label{eigen}
\sum_{l',m'}R_{l'm'}^{lm}C_{l'm'}=\lambda\sigma_lC_{lm} \ .
\end{equation}
Here the matrix $R$ is defined by
\begin{eqnarray}
R_{l'm'}^{lm}&=&\int d\svec \ \svec\cdot\hat\Bk Y_{lm}(\svec;\hat\Bk)Y_{l'm'}^*(\svec;\hat\Bk) \\
&=& \delta_{mm'}\left(b_{lm}\delta_{l',l-1} + b_{l+1,m}\delta_{l',l-1}\right) \ ,
\nonumber
\end{eqnarray}
where
\begin{equation}
b_{lm} = \sqrt{(l^2-m^2)/(4l^2-1)} \ ,
\end{equation}
and
\begin{equation}
\sigma_l=\mu_a + \mu_s(1-p_l) \ .
\end{equation}
Eq.~(\ref{eigen}) defines a generalized eigenproblem which can be transformed into a standard eigenproblem as follows. Define the diagonal matrix $S_{l'm'}^{lm}=\delta_{mm'}\delta_{ll'}\sqrt{\sigma_l}$. Note that $\sigma_l > 0$ since $p_l\le 1$ and thus $S$ is well defined. We then pre and post multiply $R$ by $S^{-1}$ and find that $W\psi=\lambda\psi$ where $W=S^{-1}RS^{-1}$ and $\psi=SC$. It can be shown that $W$ is symmetric and block tridiagonal with both a discrete and continuous spectrum of eigenvalues $\lambda_\mu$ and a corresponding complete orthonormal set of eigenvectors $\psi_\mu$, indexed by $\mu$~\cite{Panasyuk_2006_1}. We thus see that the modes (\ref{modes}), which are labeled by $\mu$, the transverse wave vector $\Bq$, and the direction of decay, are of the form
\begin{equation}
\label{modes_general}
I_{\Bq\mu}^\pm(\Br,\svec)=\sum_{l,m}\sum_{m'} \frac{1}{\sqrt\sigma_l}\psi_{lm}^\mu D_{mm'}^l(\varphi,\theta,0)Y_{lm'}(\svec)e^{i\Bq\cdot\Brho \mp Q_\mu(\Bq)z} \ ,
\end{equation}
where 
\begin{equation}
Q_\mu(\Bq)=\sqrt{q^2 + 1/\lambda_\mu^2} \ . 
\end{equation}

The Green's function for the RTE in the $z\ge 0$ half-space may be constructed as a superposition of the above modes:
\begin{equation}
\label{G_modes}
G(\Br,\svec;\Br',\svec')=\int \frac{d\Bq}{(2\pi)^2} \sum_\mu {\mathcal A}_{\Bq\mu}
I_{\Bq\mu}^\pm(\Br,\svec)I_{-\Bq\mu}^\mp(\Br',-\svec) \ ,
\end{equation}
where the upper sign is chosen if $z>z'$, the lower sign is chosen if $z<z'$ and the coefficients ${\mathcal A}_{\Bq\mu}$ are found from the boundary conditions. Using this result, we see that $G$ can be written as the plane-wave decomposition
\begin{equation}
\label{Green_RTE}
G(\Br,\svec;\Br',\svec')=\int \frac{d\Bq}{(2\pi)^2} \sum_{lm,l'm'}
g_{l'm'}^{lm}(z,z';\Bq)e^{i\Bq\cdot(\Brho-\Brho')}Y_{lm}(\svec)
Y_{l'm'}^*(\svec') \ ,
\end{equation}
where
\begin{eqnarray}
\label{g_RTE_def}
g_{l'm'}^{lm}(z,z';\Bq) &=& \frac{1}{\sqrt{\sigma_l\sigma_l'}}\sum_\mu \sum_{M,M'}{\mathcal A}_{\Bq\mu}\psi_{lm}^\mu \psi_{l'm'}^\mu  \\ &&\times D_{mM}^l(\varphi,\theta,0) D_{m'M'}^{l'}(\varphi,\theta,0)e^{-Q_\mu(\Bq)|z-z'|} 
\nonumber \\ 
\label{B_def}
&\equiv& \sum_\mu {\mathcal B}_{l'm'}^{lm}(\Bq,\mu)e^{-Q_\mu(\Bq)|z-z'|} \ ,
\end{eqnarray}
which defines ${\mathcal B}_{l'm'}^{lm}$.
It is important to note that the dependence of $G$ on the coordinates
$\Br,\Br'$ and directions $\svec,\svec'$ is \emph{explicit} and that the expansion is computable for any rotationally invariant phase function.

\subsection{Diffuse light\label{sec:diffuse}}

The diffusion approximation (DA) to the RTE is widely used in applications. It is valid in the regime where the scattering length $l_s = 1/\mu_s$ is small compared to the distance of propagation. The standard approach to the DA is through the $P_N$ approximation, in which the angular dependence of the specific intensity is expanded in spherical harmonics. The DA is obtained if the expansion is truncated at first order. The DA may also be derived using asymptotic methods~\cite{Larsen_1974_1}. The advantage of this approach is that it leads to error estimates and treats the problem of boundary conditions for the resulting diffusion equation in a natural way. 

The DA holds when the scattering coefficient is large, the absorption coefficient is small, the point of observation is far from the boundary of the medium and the time-scale is sufficiently long. Accordingly, we perform the rescaling
\begin{equation}
\mu_a \to \epsilon\mu_a \ , \ \ \ \mu_s \to \frac{1}{\epsilon}\mu_s \ , \ \ \ 
t \to t/\epsilon \ ,
\end{equation}
where $\epsilon \ll 1$. Thus the RTE~\eref{RTE_general} becomes
\begin{equation}
\label{RTE_rescaled}
\frac{\epsilon^2}{c} \frac{\partial I}{\partial t} + \epsilon \svec\cdot\grad I
+ \epsilon^2 \mu_a I + \mu_s I = \mu_s \int p(\svec,\svec')I(\Br,\svec') d\svec' \ .
\end{equation}
We then introduce the asymptotic expansion for the specific intensity
\begin{equation}
I(\Br,\svec) = I_0(\Br,\svec) + \epsilon I_1(\Br,\svec) + \epsilon^2 I_2(\Br,\svec) + \cdots \ 
\end{equation}
which we substitute into \eref{RTE_rescaled}. Upon collecting terms of $O(1)$, $O(\epsilon)$ and $O(\epsilon^2)$ we have
\begin{eqnarray}
\int  p(\svec,\svec') I_0(\Br,\svec') d\svec' = I_0(\Br,\svec) \ , \\
\svec\cdot\grad I_0 + \mu_s I_1 = \mu_s \int p(\svec,\svec')I_1(\Br,\svec') d\svec' \ , \\
\frac{1}{c} \frac{\partial I_0}{\partial t} + \svec\cdot\grad I_1
+  \mu_a I_0 + \mu_s I_2 = \mu_s \int p(\svec,\svec')I_2(\Br,\svec') d\svec'  \ .
\label{order_epsilon_squared}
\end{eqnarray}
The normalization condition~\eref{normalization_p} forces $I_0$ to depend only upon the spatial coordinate $\Br$. 
If the phase function $p(\svec,\svec')$ depends only upon the quantity $\svec\cdot\svec'$, it can be seen that 
\begin{equation}
\label{def_I_1}
I_1(\Br,\svec)= -\frac{1}{1-g}\svec\cdot\grad I_0(\Br) \ ,
\end{equation}
where the anisotropy $g$ is given by 
\begin{equation}
g=\int \svec\cdot\svec' p(\svec\cdot\svec') d\svec' \ ,
\label{eq:gdef}
\end{equation}
with $-1 < g < 1$. Note that $g=0$ corresponds to isotropic scattering and $g=1$ to extreme forward scattering. If we insert the above expression for $I_1$ into \eref{order_epsilon_squared} and integrate over $\svec$, we obtain the diffusion equation for the energy density $\U$:
\begin{equation}
\frac{\partial}{\partial t}\U(\Br,t) = \div\left[D(\Br)\grad \U(\Br,t)\right] - c\mu_a(\Br) u(\Br,t)  \ ,
\label{eq:da_time}
\end{equation}
where $I_0 = c\U/(4\pi)$.
Here the diffusion coefficient is defined by 
\begin{equation}
\label{def_D}
D= \frac{1}{3} c \ell^* \ , \ \ \  \ell^* = \frac{1}{(1-g)\mu_t} \ ,
\end{equation}
where $\ell^*$ is known as the transport mean free path. The usual expression for $\ell^*$ obtained from the $P_N$ method
\begin{equation}
\ell^* = \frac{1}{(1-g)\mu_s + \mu_a} \ ,
\end{equation}
is asymptotically equivalent to~\eref{def_D} since $\mu_a = \epsilon^2 \mu_s$.

The above derivation of the DA holds in an infinite medium and at long times. In a bounded domain, it is necessary to account for both boundary and initial layers, since the boundary conditions for the diffusion equation and the RTE are not compatible~\cite{Larsen_1974_1}. For the remainder of this paper, we will make use of the Robin boundary condition for a bounded domain $\domain$, which is of the form
\begin{equation}
\sbc\U:= c\U +  2 \bccons \diffcoef\ \n \cdot \grad \U = \Jin  \quad {\rm on} \quad \partial\domain \ ,
\label{eq:da_time_bc}
\end{equation}
where $\n$ is the outward unit normal and the extrapolation length $\ell_{\mathrm{ext}} = 2 \bccons \diffcoef/c$ can be computed from radiative transport theory~\cite{Case_book}. We note that $\bccons=0$ corresponds to an absorbing boundary and $\bccons \to \infty$ to a reflecting boundary. In the case of a boundary between diffusing media of differing refractive index $n_0$ and $n_1$ then $\bccons$ and $\ell_{\mathrm{ext}}$ take the form 
\begin{equation}
\bccons = \frac{1+R}{1-R} \ , \quad  \ell_{\mathrm{ext}} = \frac{2}{3}\ell^{\ast}\dfrac{1+R}{1-R}
\ .
\end{equation}
where $R$ is the unpolarized Fresnel reflection coefficient.
We also define the outgoing current
\begin{equation}
\mbc \U := \Jout = -\diffcoef\ \n\cdot \nabla\U = \frac{1}{2 \bccons}\left(c\U - \Jin\right) 
\label{eq:diffusion_measureable}
\end{equation}

If the phase function is not rotationally invariant, then the diffusion coefficient in~\eref{eq:da_time} becomes a second rank tensor, whose elements are given by angular moments of $p(\svec,\svec')$. The diffusion process is then anisotropic. We return to anisotropic diffusion briefly in \sref{sect:anisotropy}. We also note that if the speed of light is not constant, then the above derivation of the DA can be modified beginning from~\eref{RTE_variable_n}~\cite{Bal_2006_1}. The resulting diffusion equation depends explicitly upon the index of refraction $n$ and is of the form
\begin{equation}
\label{diff_eq_n}
\frac{\partial}{\partial t}\U(\Br,t) = \div\left[\frac{D(\Br)}{n^2(\Br)} \grad \left(n^2(\Br) \U(\Br,t)\right)\right] - c \mu_a(\Br) \U(\Br,t) \ .
\end{equation}

We now consider the scattering theory of time-harmonic diffuse waves. Assuming an $\e^{i\omega t}$ time-dependence with modulation frequency $\omega$, the energy density obeys the equation
\begin{equation}
\label{diff}
-\nabla \cdot[ D({\bm r})\nabla \Phi({\bm r})] + (c\mu_a({\bm r})+i\omega)\Phi({\bm r}) = 0 \ ,
\end{equation}
In addition to (\ref{diff}), the energy density must satisfy the boundary condition~\eref{eq:da_time_bc}. The solution to \eref{diff} obeys the Lippmann-Schwinger equation
\begin{equation}
\label{LS_diff}
\Phi = \Phi_i - GV\Phi
\end{equation}
where $\Phi_i$ is the energy density of the incident diffuse wave and $G$ is the Green's function for a homogeneous medium with absorption $\bar\mu_a$ and diffusion coefficient $D_0$. We have also
made use of the specific form of the generalised potential operator introduced in \eref{eq:ptlop_def}, which for the diffusion equation is given by
\begin{equation}
\label{V_def}
{V} = c\delta\mu_a - \nabla
  \cdot (\delta D  \nabla)  \ ,
\end{equation}
where $\delta\mu_a=\mu_a-\bar\mu_a$ and
$\delta D=D-D_0$. The unperturbed Green's
function $G({\bm r},{\bm r}^{\prime})$ 
satisfies
\begin{equation}
\label{G0_def}
\left (\nabla^2 - k^2 \right) G({\bm r},{\bm r}^{\prime}) =
-\frac{1}{D_0}\delta({\bm r}-{\bm r}^{\prime}) \ ,
\end{equation}
where the diffuse wave number $k$ is given by
\begin{equation}
k^2 = \frac{c\bar\mu_a + i\omega}{D_0} \ .
\label{k_def}
\end{equation}
We note here that the fundamental solution to the diffusion equation is given by
\begin{equation}
\label{diff_fund}
G(\Br,\Br') = \frac{1}{4\pi D} \frac{e^{-k|\Br-\Br'|}}{|\Br-\Br'|} \ .
\end{equation}
By iterating (\ref{LS_diff}) beginning with $\Phi=\Phi_i$, we obtain (compare to \eref{eq:born_general})
\begin{eqnarray}
\label{pert_V}
\Phi = \Phi_i - GV\Phi_i + GVGV\Phi_i + \cdots \ ,
\end{eqnarray}
which is the analog of the Born series for diffuse waves.

\label{sec:physics}

\section{Optical Tomography Modalities}
\label{sect:modalities}
%

We may now define some of the different problems of optical tomography
\subsection{Optical Tomography Based on the Diffusion Equation\label{sect:OT-DA}}

In diffuse optical tomography the  space $\sspace$ consists of pairs of functions $\param \equiv (\mua,\diffcoef)$  and the governing PDE $\PDE$ is the diffusion equation. We consider the domain as a simply connected set in $\mathbb{R}^n$
\begin{equation}
\gdomain \equiv \domain\,,\qquad \bgdomain \equiv \bdomain
\end{equation}
The source functions $\Jin$ and their resultant outgoing currents $\Jout$ 
belong to the space $H^{-1/2}(\bdomain)$.
For the time-dependent case the PDE $\PDE$ is given by \eref{eq:da_time} 
and in the frequency domain case by \eref{diff}. In both cases $\U$ nust satisfy the Robin boundary condition $\sbc$ given by \eref{eq:da_time_bc} and Neumann boundary condition $\mbc$ given by \eref{eq:diffusion_measureable},
therefore the mapping 
\begin{equation}
\tfn : H^{-1/2}(\bdomain) \rightarrow H^{-1/2}(\bdomain)
\end{equation}
is the Robin-to-Neumann map

\subsubsection{Diffuse Optical Tomography (DOT)\label{sect:DOT}}
The non-linear inverse problem in DOT is to recover  functions $\mua \in \absspace, \diffcoef \in \kapspace$ from measurements  $\y \in \dspace$. 
The data may be collected either in the time domain or in the frequency domain 
either at one or several modulation frequencies $\omega$ \cite{hebden97a}. If 
the former is used, it is typical to Fourier Transform the data and to develop 
the analysis in the frequcny domain. At frequency $\omega = 0$ ('DC' 
measurements), the data simply represents a total photon count without any 
phase information, and the recovery of both
$\mua$ and $\diffcoef$ suffers from non-uniqueness\cite{arridge98b} (but see\cite{harrach2009}). In this case a simpler problem is commonly defined in which one parameter (usually scattering) is assumed known and the forward mapping is restricted to recovery of the second function.

For the DOT problem the \Frechet\ derivative \eref{eq:gfrechet}  has a kernel based on the density functions for absorption and diffusion 
\begin{eqnarray}
\mtrx{c}{\pmdf_{\mua}(\r;\omega) \\ \pmdf_{\diffcoef}(\r;\omega)} &=&  \mtrx{c}{\aU(\r;\omega) \U(\r;\omega)\\ \nabla\aU(\r;\omega)\cdot\nabla\U(\r;\omega)}\\
\mtrx{c}{\pmdf_{\mua}(\r,t) \\ \pmdf_{\diffcoef}(\r,t)} &=& \mtrx{c}{\int_0^t \aU(\r,-t')  \U(\r,t-t')\rmd t'\\ \int_0^t \nabla\aU(\r,-t')\cdot\nabla\U(\r,t-t')\rmd t'}
\end{eqnarray}
where we note that the adjoint problem leads to a backward-time equation.

Although DOT is a non-linear inverse problem, a  linearised version is often considered, which we refer to as \emph{Difference Diffuse Optical Tomography} (DDOT). 
Here we take differences in measurements $\y^{\delta} = y_2 - y_1$  and formulate a linear mapping as in \eref{eq:lin_discrete}.
Clearly the linear mapping $A$ is given by the linear \Frechet\ derivative operator, evaluated at $\param_1 = (\submua{1},\diffcoef_1)$. As in nonlinear DOT, the simplified case of absorption only imaging is frequently considered. The success of this approach depends on how sensitive the reconstruction is to the correct choice of linearisation point $\param_1$. This typically involves careful calibration and consideration of the modelling errors. We return to this point in 
\sref{sect:errors_priors}.

\subsubsection{Fluorescence Diffuse Optical Tomography (FDOT)\label{sect:FDOT}}
In fluorescence optical tomography sources are introduced at an \emph{excitation} wavelength $\lambda^{e}$ giving rise to an excitation field $\U^{e}$.  Fluorescence is regarded as a function $\fx(\r)$  which governs the absorption of radiation at wavelength $\lambda^{e}$ and (partial) re-emission as a source at the longer wavelength (lower energy) $\lambda^{f}$. Measurements are taken at $\lambda^{f}$ and $\lambda^{e}$.
Using the frequency domain formulation, we  consider two coupled PDEs 
\begin{eqnarray}
\PDE(\param^{e})\Ue &\equiv&  -\nabla \cdot \diffcoef^{e}(\r) \nabla \Ue(\r; \omega) 
+\left( c \mua^{e}(\r) + {\myi \omega}\right)\Ue(\r;\omega)= 0 \label{eq:FDOT1}\\
\PDE(\param^{f})\Uef &\equiv&  -\nabla \cdot \diffcoef^{f}(\r) \nabla \Uef (\r; \omega) 
+\left(c \mua^{f}(\r) + {\myi \omega}\right)\Uef(\r;\omega) = \nonumber \\
&& \hspace{1cm} \fx(\r;\omega) \Ue(\r;\omega)\label{eq:FDOT2}\\
\sbc  \Ue &\equiv& c \Ue(\rd;\omega) + 2\bccons\diffcoef^{e}(\rd)
\mneumop{\Ue(\rd;\omega)} =  \eJin(\rs;\omega)\\
\sbc \Uef &\equiv& c \Uef(\rd;\omega) + 2\bccons\diffcoef^{f}(\rd)
\mneumop{\Uef(\rd;\omega)} = 0 \\
\mbc \Ue  &\equiv& \eJout(\rd;\omega) = -\diffcoef^{e}(\rd)\mneumop{\Ue(\rd;\omega)}\label{eq:FDOT-m1}\\
\mbc \Uef &\equiv& \efJout(\rd;\omega) = -\diffcoef^{f}(\rd)\mneumop{\Uef(\rd;\omega)}\label{eq:FDOT-m2}
\end{eqnarray}
The quantity $\fx$ is a product of the concentration of fluorescent material and its conversion efficiency. Since re-emission
is a Poisson process it is characterised by a lifetime $\tau$. In frequency domain this leads to a complex valued parameter
\begin{equation}
\fx(\r,\omega) = \fx_0(\r) \frac{1}{1 + \rmi \omega \tau(\r)}
\label{eq:fltf}
\end{equation}

The simplest problem considered is to recover $\fx \in \fluspace$ from measurements $\y \in \dspace^{\ef}$ and the linear operator
\begin{equation}
\mfmap_{\fx}^{\ef} : \sspace_{\fx} \rightarrow \dspace^{\ef}
\end{equation}
Consider the frequency domain version, the kernel of the \Frechet\ derivative \eref{eq:gfrechet} in this case is the density function with respect to $\fx$ given by
\begin{equation}
{\pmdf_{\fx}(\r;\omega)}= {\aUef(\r;\omega) \Ue(\r;\omega)}
\label{eq:frechet_FDOT}
\end{equation}
giving a complex valued reconstruction from which $\fx_0$, $\tau$ can be recovered using \eref{eq:fltf}. Clearly, for DC measurements, only the fluorescence $\fx_0$ can be recovered, not the lifetime; \emph{Fluorescence Lifetime Imaging Tomography} (FLIM tomography) is a term used to emphasise that both  $\fx_0$ and $\tau$ are being recovered. Time domain measurements may of course be used to solve this problem. If the data are Fourier transformed into the frequency domain the techniques for image reconstruction are unchanged. It may also be expressed directly in the time domain. In this case the kernel is given by a double convolution
\begin{equation}
{\pmdf_{\fx}(\r;t)}= \int_0^t \int_0^{t'}\aUf(\r;-t'') \Ue(\r;t-t') \fx(\r,t''-t')\rmd t' \rmd t''
\label{eq:flim_kernel}
\end{equation}
where
\begin{equation}
\fx(\r,t''-t') = \fx_0(\r) \rme^{-\frac{|t''-t'|}{\tau(\r)}}
\label{eq:fltt}
\end{equation}
gives explicitly the exponential decay of the Poisson process. In practice the half-life does not have a discrete value but a distribution, and the linear integral equation with \eref{eq:flim_kernel} as its kernel has to also be convolved with the finite time spread function of the measurement system~\cite{soloviev2007a}.
\subsubsection{Diffuse Optical Tomography with Fluorescence (DOT-FDOT)\label{sect:DOT-FDOT}}

The discussion of the linear  FDOT problem  in \sref{sect:FDOT} is predictated on knowing the parameters  $\param^{e}$ and $\param^{f}$ used to construct the kernel in \eref{eq:frechet_FDOT}. These could themselves be recovered using uncoupled DOT reconstructions at the excitation and fluorescence wavelengths.
Therefore it is natural to consider the joint problem
(using the notation of \eref{eq:FDOT1}-\eref{eq:FDOT-m2})

\begin{eqnarray}
\PDE(\param^{e})\Ue = 0\, \quad & \sbc \Ue = \eJin \,\quad & \mbc \Ue = \eJout \label{eq:DFDOTe}\\
\PDE(\param^{f})\Uf = 0\, \quad & \sbc \Uf = \fJin \,\quad & \mbc \Uf = \fJout\label{eq:DFDOTf}\\
\PDE(\param^{f})\Uef = \fx \Ue\,\quad  & \sbc \Uef = 0\,\quad  & \mbc \Uef = \efJout\label{eq:DFDOTef}
\end{eqnarray}

with the joint forward operator
\begin{equation}
\mtrx{l}{\mfmap_{\param}^e \\ \mfmap_{\param}^f \\ \mfmap_{\fx}^{\ef}} : \mtrx{c}{\absspacel{e},\kapspacel{e}\\\absspacel{f},\kapspacel{f} \\ \fluspace } \rightarrow \mtrx{l}{\dspace^e\\\dspace^f\\\dspace^{\ef}}
\end{equation}
The solution of the joint problem is computationally more demanding than the 
three seperate problems. 
This problem has hardly been addressed in the frequency domain in \cite{milstein2002a} and in the time domain in\cite{soloviev2008a}.


\subsubsection{Multispectral Diffuse Optical Tomography}

Whereas each of the above problems could be considered at a range of spectral samples and the spectral variation of the recovered solutions determined, the idea in multispectral DOT is to reformulate the problem into the recovery of a set of images of known chromophores with well characterised spectral dependence :
\begin{equation}
\mua(\lambda_j) = \sum_i \epsilon_i(\lambda_j) c_i \rightarrow \mua(\msset) = {\boldsymbol \epsilon}{\bf c}
\end{equation}
where ${\boldsymbol \epsilon}$ is a known matrix.
Similarly a spectral dependence of scattering can be written as
\begin{equation}
\mus(\lambda) = a \lambda ^{-b}
\end{equation}
This leads to a formulation
\begin{equation}
\mfmap^{\msset}_{\{\myvec{c},a,b\}} : \sspace_{\{\myvec{c},a,b\}} \rightarrow \dspace^{\msset}
\end{equation}

The spectrally decoupled linearisation is represented as
\begin{equation}
\begin{mtrx}{c}{
	\yv(\lambda_1)\\
	\yv(\lambda_2)\\
	\vdots\\
	\yv(\lambda_N)}
\end{mtrx} ~=~
\begin{mtrx}{cccccccccc}
{
\mA_{\mua}(\lambda_1) & \ldots & \mymat{0}  & 
\mA_{\mus}(\lambda_1) & \ldots & \mymat{0} \\
\vdots & \vdots & \vdots & \vdots & \vdots & \vdots \\
\mymat{0} & \ldots & \mA_{\mua}(\lambda_N) & 
\mymat{0} &\ldots & \mA_{\mus}(\lambda_N)}
\end{mtrx}
~
\begin{mtrx}{c}
{	\mua^{\delta}(\lambda_1)\\
	\vdots\\
	\mua^{\delta}(\lambda_N)\\
	\mus^{\delta}(\lambda_1)\\
	\vdots\\
	\mus^{\delta}(\lambda_N)
}
\end{mtrx}
\label{eq:ms1_linearsed}
\end{equation}

The spectral component at each pixel is independent so we introduce the 
diagonal matrices 
\begin{equation}
\mymat{C}_k(\lambda_j) := c_k(\lambda_j)\mymat{I} \qquad
\mymat{B}_{a}(\lambda_j) := \lambda_j^{-b}\mymat{I} \qquad
\mymat{B}_{b}(\lambda_n) := a \lambda_j^{-b}\ln \lambda_j \mymat{I} 
\end{equation}
and we can now represent  \eref{eq:ms1_linearsed} as
\begin{eqnarray}
\fl
\begin{mtrx}{c}{
	\yv(\lambda_1)\\
	\yv(\lambda_2)\\
	\vdots\\
	\yv(\lambda_N)
}\end{mtrx} ~=~
\begin{mtrx}{cccccccccc}
{
\mA_{\mua}(\lambda_1) & \ldots & \mymat{0}  & 
\mA_{\mus}(\lambda_1) & \ldots & \mymat{0} \\
\vdots & \vdots & \vdots & \vdots & \vdots & \vdots \\
\mymat{0} & \ldots & \mA_{\mua}(\lambda_N) & 
\mymat{0} &\ldots & \mA_{\mus}(\lambda_N)
}\end{mtrx} \times \nonumber \\
\begin{mtrx}{cccccc}
{
\mymat{C}_1(\lambda_1) & \ldots & \mymat{C}_K(\lambda_1) & \mymat{0} & \mymat{0}\\
\vdots & \vdots & \vdots & \vdots & \vdots &\vdots\\
\mymat{C}_1(\lambda_n) & \ldots & \mymat{C}_K(\lambda_n) & \mymat{0} & \mymat{0} \\
\vdots & \vdots & \vdots & \vdots &\vdots &\vdots\\
\mymat{C}_1(\lambda_N) & \ldots & \mymat{C}_K(\lambda_N) & \mymat{0}  & \mymat{0}\\
\mymat{0} & \ldots & \mymat{0}  & \mymat{B}_{a} (\lambda_1) & \mymat{B}_{b} (\lambda_1)\\
\vdots & \vdots & \vdots & \vdots & \vdots &\vdots\\
\mymat{0} & \ldots & \mymat{0}  & \mymat{B}_{a} (\lambda_n) & \mymat{B}_{b} (\lambda_n) \\
\vdots & \vdots & \vdots & \vdots & \vdots &\vdots \\
\mymat{0} & \ldots & \mymat{0}  & \mymat{B}_{a} (\lambda_N) & \mymat{B}_{b} (\lambda_1) 
}\end{mtrx}
\begin{mtrx}{c}{
\myvec{c}^{\delta}_1\\
\vdots\\
\myvec{c}^{\delta}_k\\
\vdots\\
\myvec{c}^{\delta}_K\\
\myvec{a}^{\delta}\\
\myvec{b}^{\delta}
}\end{mtrx}
\label{eq:ms2_linearsed}
\end{eqnarray}

See\cite{corlu2003,corlu2005}

Similarly to MSDOT in \emph{Multispectral fluorescence DOT} (MSFDOT) the fluorescence may be considered a linear combination of fluorophores emitting radiation in a known spectral pattern
\begin{equation}
\fx(\lambda_j) = \sum_i \epsilon_i(\lambda_j) \flph_i \rightarrow \fx(\msset) = {\boldsymbol \epsilon}^f\,{\bf \flph}
\end{equation}
In this case a linear forward operator may be constructed 
\begin{equation}
	\mfmap_{\mathrm{\bf \flph}}^{\msset^e\rightarrow\msset^f} : \mfpespace  \rightarrow \dspace^{\msset^e\rightarrow\msset^f}
\end{equation}
See\cite{zacharopoulos2009a}.
Finally, a very general problem could be considered in the recovery of all chromophores and flurophores from the operator

\begin{equation}
\mtrx{l}{\mfmap^{\msset^e,\msset^f}_{\{\myvec{c},a,b\}} \\ \mfmap_{\mathrm{\bf \flph}}^{\msset^e\rightarrow\msset^f}} : 
\mtrx{c}{ \sspace_{\{\myvec{c},a,b\}}\\ \mfpespace  } \rightarrow \mtrx{l}{\dspace^{\msset^e,\msset^f}\\ \dspace^{\msset^e\rightarrow\msset^f}}
\end{equation}
\subsubsection{Bioluminescence (BDOT)\label{sect:BDOT}}
Bioluminescence optical tomography 
detects the light emitted from within tissues
through the action of luciferase on its substrate, luciferin.
The
luciferase emission spectrum is quite wide\cite{kuo2005} (500-800 nm) and
overlaps with so-called near infrared window of biological tissues (650-850
nm), where the light scattering dominates over the light absorption. Red and
near infrared components of the emission spectrum penetrate biological tissue
appreciably well allowing to image bioluminescence objects imbedded deeply
inside tissue.


In terms of the general framework established in \sref{sect:general_framework}
the bioluminescence tomography forward problem is posed as
\begin{eqnarray}
	\PDE(\param) \U & =& q \label{eq:BDOT_gpde}\\
	\sbc \U &=& 0 \label{eq:BDOT_gin}\\
        \Jout &=& \mbc \U \label{eq:BDOT_gout}
\end{eqnarray}
There are no incoming radiation sources and spontaneously created light is subject to the same boundary conditions as in all other problems considered in this paper.
In the inverse problem, the unknown function $q$ is to be reconstructed from measured values of function $\Jout$ on the surface. This is a linear inverse problem which can be formulated in terms of a Fredholm integral equation of the first kind%
\begin{equation}
\Jout\left(\grd\right) =\gdomintvar {G(\grd,\gr)q(\gr)}{\gr},  \label{eq:BDOT_Fredholm}
\end{equation}%

In \cite{gu2004} the bioluminescent source distribution was recovered from a monochromatic data set and theoretical acpects of uniqueness of solutions were discussed in\cite{wang_G2004,cong_w2005}.
The idea of using multiple wavelengths in
bioluminescence imaging was introduced  by Chaudhari 
\textit{et al}\cite{chaudhari2005a} in which this approach is applied for 3D localization of
deep sources within a mouse phantom\textit{in vivo}. 
Combination of bioluminescence and Positron Emission Tomography (PET)
is discussed in\cite{alexandrakis2005a}.

\subsection{Optical Tomography based on the Radiative Transfer Equation}

Each of the problems in \sref{sect:OT-DA} could be considered based on the RTE instead 
of the diffusion equation. The parameters considered in the inverse problem 
are the same but the measurements are different. There are two possibilities :
angularly resolved measurements or angularly independent measurements. In the 
latter case the size of the measured dataset is the same as in the diffusion 
problems, but the construction of the linearisation is different. In the case 
of angularly dependent data the data set is much bigger. Issues regarding uniqueness depend critically on which measurements are made.

In the RTE-based problems, the generalised domain 
becomes $\gdomain = \domain \times S^{n-1}$ with boundary $\bgdomain = \rtebnd_- \cup \rtebnd_+$ where $\rtebnd_{\pm} = \bdomain \times S^{n-1}_{\pm} (\bnormal)$.
The governing PDE $\PDE$ is given by \eref{RTE_general} with parameters $\mua, \must$ and the boundary condition $\sbc$ specifies the specific intensity on $\rtebnd_-$ as given by \eref{RTE_bc_ihg}.
The transfer function \eref{eq:transfer_function} is known as the \emph{Albedo Operator}
\begin{equation}
\Lambda : L_1(\rtebnd_-) \rightarrow  L_1(\rtebnd_+)
\label{eq:albedo_op}
\end{equation}
which maps directional incoming radiation on $\bdomain$ to directional outgoing radiation 
\begin{equation}
\Jout_{\sind,\sind'}(\rd, \angvec)|_{\angvec \cdot \bnormvec > 0} = \Lambda \Jin_{\sind,\sind'}(\rs, \angvec)|_{\angvec \cdot \bnormvec < 0}
\end{equation}
where the set of source functions $\myset{S}$ is indexed by position index  $\sind$ and direction index $\sind'$.
 
As in the diffusion case, the specific problems may be time-dependent, time-independent, or frequency domain. The role of the aperture function is crucial.
In angular resolved meaurements
\begin{eqnarray}
	\y_{\apind,\apind',\sind,\sind'} &=& \left<w_{\apind,\apind'},\Jout_{\sind,\sind'}\right>_{\bgdomain} \nonumber \\
&=& \angint{
w_{\apind'}(\angvec) 
\boundintvar{\ap(\rd)\Jout_{\sind,\sind'}(\rd,\angvec)}{\rd}
}
\end{eqnarray}

The sensitivity functions are given by

\begin{eqnarray}
\hspace{-1.5cm}
  \rho_{\mua}(\r)  &=&
\angint {\; I^{\ast}(\r,\angvec)  I(\r,\angvec) } \\
\hspace{-1.5cm}
\rho_{\must}(\r) &=&
\angint {\; I^{\ast}(\r,\angvec)  I(\r,\angvec) } - \angint {\;\; \angintp {\; p (\angvec\cdot \angvecp) I^{\ast}(\r,\angvec) } I(\r,\angvec) }
\end{eqnarray}
where $I^{\ast}(\r,\angvec)$ is the solution to the adjoint transport equation with source boundary condition $w_{\apind'}(\angvec)\ap(\rd) $ on $\rtebnd_+$


Note that the equivalent right-hand side for the fluorescence problem in \eref{eq:FDOT2} is assumed isotropic. I.e.
\begin{equation}
q^f(\r) = \fx(\r) \angint{ I^e(\r,\angvec)}
\end{equation}

For detailed discussion of the uniqueness and stability of the different problems in inverse transport arising from angularly resolved vs angularly averaged and time-dependent vs time-independent measurements we refer to the recent review article\cite{bal2009a}.

\section{Forward Modelling Methods}
\label{sect:modelling_methods}
The inversion methods described in the sequel are dependent on the accuracy and efficiency with which Green's functions can be computed. The analytic form of such functions for the DA are known for several problems and geometries. The format for the RTE is generally limited to one-dimensional and other special cases (see \sref{sect:Physics_transport_regime}).


\subsection{Volume Discretisation Methods}


In the Finite Element Method (FEM) the volume  $\domain$ is discretised to a mesh
\begin{equation}
{\domain} \rightarrow \{\mathbb{T}_{\domain},\mathbb{N}_{\domain},\mathbb{U}_{\domain}\}
\label{eq:femmesh}
\end{equation}
where $\mathbb{T}_{\domain}$ is the set of $P$  elements $\elem_{e}; e=1, \ldots$,  $\mathbb{N}_{\domain}$ is the set of $N$ vertices $\myvec{N}_{p}; p=1,\ldots,N$, and 
 $\mathbb{U}_{\domain}$ is the set of $N$ 
locally supported basis functions $\fembasis_{p}(\r); p=1,\ldots,N$.
For the diffusion equation, writing 
\begin{equation}
\U(\r) \simeq \U^h(\r) =\sum_p \U_p \fembasis_{p}(\r)
\end{equation}
results to a discrete system
\begin{equation}
\mymat{K}(\param)\myvec{\U} = \myvec{q}
\end{equation}
where $\mymat{K}(\param)$ has matrix elements
\begin{eqnarray}
K_{lm} = \domintvar{ \left(\diffcoef(\r) \nabla \fembasis_l(\r) \cdot \nabla \fembasis_l(\r) + 
\left(c \mua(\r) + \myi \omega\right)\fembasis_l(\r) \fembasis_m(\r)\right)}{\r} 
+ \nonumber\\
\qquad \qquad \frac{1}{2\bccons}\boundintvar{\fembasis_l(\r) \fembasis_m(\r)}{\r}
\end{eqnarray}
and $\myvec{q}$ represents the discretisation of the boundary conditions. 
Making use of the basis expansion \eref{eq:sbs} allows the representation of 
$\mymat{K}(\param)$ as
\begin{equation}
\mymat{K} = \mymat{S} +  {\myi \omega}\mymat{B} + \sum_k \left(\diffcoef_k \mymat{K}_k^{\diffcoef} + c \mu_k  \mymat{K}_k^{\mu}\right) 
\label{eq:sysmat}
\end{equation}
where $\mymat{B}$ is the mass matrix, $\mymat{S}$ is the matrix of surface integrals with elements given by the last term on the right in \eref{eq:sysmat} and $\mymat{K}_k^{\diffcoef} \equiv \dbyd{\mymat{K}}{\diffcoef_k} $, $\mymat{K}_k^{\mu} \equiv \dbyd{\mymat{K}}{c \mu_k} $ are given by
\begin{eqnarray}
\mymat{K}_{k,lm}^{\diffcoef} &=& \domintvar{ b_k(\r) \nabla \fembasis_l(\r) \cdot \nabla \fembasis_m(\r)}{\r}\\
\mymat{K}_{k,lm}^{\mu} &=& \domintvar{ b_k(\r) \fembasis_l(\r)\fembasis_m(\r)}{\r}
\end{eqnarray}

In RTE FEM, a basis is also defined for the angular directions
\begin{equation}
S^{n-1} \rightarrow \mathbb{V}_{S^{n-1}} = \{\mathbb{T}_{S^{n-1}},\mathbb{N}_{S^{n-1}},\mathbb{U}_{S^{n-1}}\}
\end{equation}
Different schemes result by choosing different basis on the unit circle $S^1$ or sphere $S^2$
\begin{enumerate}
\item The discrete Ordinate method chooses a discrete set: $\mathbb{U}_{S^n} = \set{\delta(\angvec-\angvec_{p}}$ which are chosen to give exact quadrature points for a spherical harmonic expansion of the angular variable. The radiance is thus represented explicitly in a set of ray directions\cite{abdoulaev2003}.
\item The $P_N$ method chooses the spherical harmonics (rotated into real functions) directly:  $\mathbb{U}_{S^n} = \set{\tilde{Y}_p(\angvec)}$. This allows explicitly the representation as diffusion in the lowest order with higher terms representing higher order effects \cite{aydin2002}
\item The local basis function method chooses a set of locally supported basis functions on $\mathbb{T}_{S^{n-1}}$. This allows essentially the same machinary of sparse matrix manipulation as for the spatial variables\cite{richling01,tarvainen2005,tarvainen2005a}.
\item The wavelet basis chooses a heirarchical set of wavelets on the sphere to allow for variable and adaptive angular discrisation\cite{buchan2005}.
\end{enumerate} 
%

The RTE has been utilized as the forward model for optical 
tomography reconstruction in a few studies. 
In most of these papers, the forward solution of the RTE has been based either 
on the finite difference method (FDM) or the finite volume method (FVM) 
and the discrete ordinate  formulation of the RTE, see \cite{dorn98,klose1999a,klose2002b,abdoulaev2005,ren2006a}.  
In the FEM solution of 
the RTE in low-scattering or non-scattering regions, the ray-effect may 
become visible\cite{lathrop68,lathrop71}. 
In \cite{tarvainen2005b}, the FEM model of \cite{tarvainen2005,tarvainen2005a} was augmented with the streamline diffusion modification to stabilize the forward solution in this case. The streamline diffusion modification 
has been found to stabilise numerical solutions of the RTE in situations in 
which standard techniques produce oscillating results \cite{kanschat98,richling01,tarvainen2005b}.  
Due to the heavy computational and memory requirements of meshes, adaptive techniques have been used within an FVM approach\cite{soloviev2006a,soloviev2006b}, as well as the FEM scheme\cite{kwon2006,guven2007a,guven2007b,joshi2006b,lee_J2007}.
In order to overcome the difficulty of 3D meshing, meshless methods can be used\cite{qin2008}.

\subsection{Boundary Discretisation Methods}


Rather than meshing a volume we may consider the domain as the union of a number of subdomains
\begin{equation}
\domain = \cup_{\bdind}\domain_{\bdind}\,,\qquad \bdind = 1\ldots L
\label{eq:subdomains}
\end{equation}
together with  constant with subdomain optical parameters $\set{\param_{\bdind}}$ and a set of interfaces between domains 
\begin{equation}
\surf_j = \bdomain_{\bdind,\bdind'}\,,\quad j = 1\ldots J
\label{eq:domain_interfaces}
\end{equation}
We consider the diffusion approximation and introduce the following notation
\begin{eqnarray}
\I_{\sind} &:=& \Phi_{\sind-1}|_{\surf_{\sind}} =
\Phi_{\sind}|_{\surf_{\sind}}, \label{eq:Udef}\\ \J_{\sind} &:=&
\diffcoef_{\sind-1}
\left.\dbyd{\Phi_{\sind-1}}{\nu_{\sind-1}}\right|_{\surf_{\sind}} =
-\diffcoef_{\sind}
\left.\dbyd{\Phi_{\sind}}{\nu_{\sind}}\right|_{\surf_{\sind}}
\equiv \diffcoef_{\sind-1}
\left.\norgrad{\Phi_{\sind-1}}\right|_{\surf_{\sind}} =
\diffcoef_{\sind}
\left.\norgrad{\Phi_{\sind}}\right|_{\surf_{\sind}}\label{eq:Jdef}.
\end{eqnarray}
For simplicity we develop the discussion using only two regions.
For a more general implementation we refer to~\cite{sikora2006}.
For homogeneous region $\Omega_1, \Omega_2$, 
Green's second theorem provides the following \\
\begin{eqnarray}\label{eq:BIE_Gamma1,2}
\hspace{-1.5cm}\Phi_{1}(\r) +
\sintvar{\bdomain}{\left(\dbyd{G_{1}(\r,\rd')}{\myvec{\nu}} + 
\frac{G_{2}(\r,\rd')}{2\bccons \diffcoef_{1}}\right)
\I_{1}(\rd') }{\rd'} \nonumber\\ 
\hspace{-0.8cm}- \sintvar{\surf}{\left(\dbyd{G_{1}(\r,\rd')}{\myvec{\nu}}\I_{2}(\rd')-\frac{G_{1}(\r,\rd')}{\diffcoef_{1}}
\J_{2}(\rd')\right)}{\rd'} 
=\Q_{1}(\r) \\
\hspace{-1.5cm}\Phi_{2}(\r) +
\sintvar{\surf}{\left(\dbyd{G_{2}(\r,\rd')}{\myvec{\nu}}\I_{2}(\rd')-\frac{G_{2}(\r,\rd')}{\diffcoef_{2}}
\J_{2}(\rd')\right)}{\rd'} 
=\Q_{2}(\r)\label{eq:BIE_Gamma2}
\end{eqnarray}
where  
\begin{eqnarray}\label{bsourceterm}
\Q_{\ell}(\r)&=&\vintvar{\Omega_{\ell}}{{G_{\ell}(\r,\r')}
q_{\ell}(\r')}{\r' }
\end{eqnarray}
and the fundamental solutions $G_{\sind}$ are the three-dimensional Green's functions of the diffusion equation in an infinite medium as in \eref{diff_fund}.
Taking the limit as $\r \in \Omega_1 \rightarrow \partial \Omega$,  
 $\r \in \Omega_1 \rightarrow \surf$, and
 $\r \in \Omega_2 \rightarrow \surf$  results in three
coupled boundary integral equations (BIEs) 
in the three unknown functions 
\begin{equation}
\{ {\mathrm{f}}\}= \{\I_1,\I_2,\J_2\}\,.
\label{eq:BIEfunctions}
\end{equation}


In contrast to \eref{eq:femmesh} only the interfaces $\surf_{\sind}$ are discretised to meshes
\begin{equation}
\mathbb{V}_{\sind} = \{\mathbb{T}_{\sind},\mathbb{N}_{\sind},\mathbb{U}_{\sind}\}
\end{equation}
where $\mathbb{T}_{\sind}$ is the set of $P_{\sind}$ surface elements $\elem_{\sind,e}; e=1, \ldots, P_{\sind}$,  $\mathbb{N}_{\sind}$ is the set of $N_{\sind}$ vertices $\myvec{N}_{\sind,p}; p=1,\ldots,N_{\sind}$, and 
 $\mathbb{U}_{\sind}$ is the set of $N_{\sind}$ 
locally supported basis functions $\fembasis_{\sind,p}(\r); p=1,\ldots,N_{\sind}$.
Then the functions  \eref{eq:Udef}-\eref{eq:Jdef} are represented as
\begin{equation}\label{BEMapprox}
\I_{\sind}(\r)\approx\sum_{p=1}^{N_{\sind}}
\I_{\sind,p}\fembasis_{\sind,p}(\r), \qquad
\J_{\sind}(\r)\approx\sum_{p=1}^{N_{\sind}}
\J_{\sind,p}\fembasis_{\sind,p}(\r)\,.
\end{equation}
Using the representation ~\eref{BEMapprox} in the BIE system ~\eref{eq:BIE_Gamma1,2}-\eref{eq:BIE_Gamma2} followed by sampling at the nodal points, 
we obtain a discrete system for the  collocation Boundary Element Method
\cite{Aliabadi2002,becker92,beer01,sikora2004,sikora2002} 
as a linear matrix equation
\begin{equation}
\left(
\begin{array}{ccc}
\xi^{+}_{1}\mymat{I} + \mymat{A}_{11}^{1} + \frac{1}{2 \bccons}\mymat{B}_{11}^{1} & 
-\mymat{A}_{12}^{1} &  \mymat{B}_{12}^{1} \\
\mymat{A}_{11}^{2} + \frac{1}{2 \bccons}\mymat{B}_{11}^{2} & \xi^{-}_{1}\mymat{I} -\mymat{A}_{12}^{2} & \mymat{B}_{12}^{2} \\
0 &\xi^{+}_{2}\mymat{I} +  \mymat{A}_{22}^{2} & -\mymat{B}_{22}^{2}  
\end{array}
\right)
\mtrx{ccc}{\myvec{\I}_1\\\myvec{\I}_2\\\myvec{\J}_2}
=
\mtrx{ccc}{\myvec{Q}_1|_{\partial\Omega}\\\myvec{Q}_1|_{\surf}\\0}
\label{eq:BEM2layer}
\end{equation}
where we assumed that the source term is only in $\Omega_1$, and where  $\myvec{\I}_{\sind}$, $\myvec{\J}_{\sind}$ are the vectors of coefficients $\I_{\sind,p}$, $\J_{\sind,p}$ in \eref{BEMapprox}, and  the matrix elements are given by
\begin{eqnarray*}
A_{11}^{j}(p,p') &=& \sintvar{\bdomain} {\dbyd{G_{1}(\myvec{N}_{j,p},\rd')}{\myvec{\nu}}\fembasis_{1,p'}(\rd')}{\rd'} \\
A_{12}^{j}(p,p') &=& \sintvar{\surf}{ \dbyd{G_{1}(\myvec{N}_{j,p},\rd')}{\myvec{\nu}}\fembasis_{2,p'}(\rd')}{\rd'} \\
A_{22}^{j}(p,p') &=& \sintvar{\surf} {\dbyd{G_{2}(\myvec{N}_{j,p},\rd')}{\myvec{\nu}}\fembasis_{2,p'}(\rd')}{\rd'} \\
B_{11}^{j}(p,p') &=& \sintvar{\bdomain} {\frac{G_{1}(\myvec{N}_{j,p},\rd')}{\diffcoef_1}\fembasis_{1,p'}(\rd')}{\rd'} \\
B_{12}^{j}(p,p') &=& \sintvar{\surf} {\frac{G_{1}(\myvec{N}_{j,p},\rd')}{\diffcoef_1}\fembasis_{2,p'}(\rd')}{\rd'} \\
B_{22}^{j}(p,p') &=& \sintvar{\surf}{ \frac{G_{2}(\myvec{N}_{j,p},\rd')}{\diffcoef_2}\fembasis_{2,p'}(\rd')}{\rd'} \,.
\end{eqnarray*}
The extra function $\xi^{\pm}_{\sind}$ in~\eref{eq:BEM2layer}, arises due to singularities on the boundary. These terms can be
calculated by surrounding the point $\rd$, which lays on the
boundary, by a small hemisphere $\sigma_{\varepsilon}$ of radius
$\varepsilon$ and taking the limit
when $\varepsilon \rightarrow 0$.
However, as shown in \cite{becker92,bonnet95}, this term can be obtained indirectly by utilising some simple physical considerations. In particular, we have  $\xi^+_{i}(\r) = \xi^-_{i}(\r) = \frac 12$ when the observation point lies on a smooth surface, which is the case considered here. 

Note that appropriate numerical techniques for handling the singularity of the
kernels is required for the diagonal elements of these matrices~\cite{sikora2006}.
We represent \eref{eq:BEM2layer} as
\begin{equation}\label{THEBEM}
\mymat{T} \myvec{f }= \myvec{Q}\,,
\end{equation}
where $\myvec{f} = [\myvec{U}_1, \myvec{U}_2, \myvec{V}_2]$ is the discrete version of $\{\rm{f}\}$ in \eref{eq:BIEfunctions}.
The matrix $\mymat{T}$ takes the form of a dense un-symmetric block matrix. 

The normal flux is not always required. In this case we may take the Schur complement of \eref{eq:BEM2layer} to obtain
\begin{equation}
\mtrx{cc}{
\frac{1}{2}\mymat{I} + \mymat{A}_{11}^{1} + \frac{1}{2 \bccons}\mymat{B}_{11}^{1} & 
    -\mymat{A}_{12}^{1} +\mymat{B}_{12}^{1} \left(\mymat{B}_{22}^{2}\right)^{-1} \left(\frac{1}{2}\mymat{I} +  \mymat{A}_{22}^{2}\right) \\
\mymat{A}_{11}^{2} + \frac{1}{2 \bccons}\mymat{B}_{11}^{2} & 
    \frac{1}{2}\mymat{I} -\mymat{A}_{12}^{2} + \mymat{B}_{12}^{2} \left(\mymat{B}_{22}^{2}\right)^{-1}\left( \frac{1}{2}\mymat{I} +  \mymat{A}_{22}^{2}\right) 
}
\mtrx{c}{\myvec{\I}_1\\\myvec{\I}_2}
=
\mtrx{c}{\myvec{Q}_1|_{\partial\Omega}\\\myvec{Q}_1|_{\surf}}
\label{eq:BEM2layer_Schur}
\end{equation}
This is a smaller system than \eref{eq:BEM2layer} and therefore computationally cheaper. However, it is still usually lengthy to solve large dense matrix systems. Another approximation is obtained if  write \eref{eq:BEM2layer_Schur} in the form
\begin{equation}
\frac{1}{2}\left(\mymat{I} -2\mymat{G}\right)\myvec{\I}= \myvec{Q}\,,
\label{THEBEM_schur}
\end{equation}
We can then use a Neumann series approximation
\begin{equation}
\frac{1}{2} \myvec{\I} \simeq \left(\mymat{I} + 2 \mymat{G} + 4 \mymat{G}^2 +\ldots\right) \myvec{Q}
\end{equation}
truncation of this series at the $n^{\mathrm{th}}$ power in $\mymat{G}$ constitutes the Kirchoff approximation\cite{ripoll2001,ripoll2006a}.

\subsection{Monte Carlo Methods}

Monte Carlo methods are highly prevelant in many models of radiation in tissue. Within optical tomography, one of the early cited works include  Prahl \etal \cite{prahl89} dating from 1989 and Wang \etal \cite{Wang95} describing multi-layered tissues. For the more complex geometries involved in modelling light propogation in the head, \cite{boas02,heiskala09} describe 3D voxel based MC models including anisotropy\cite{Heiskala05}.

Use of the MC model within image reconstruction is usually limited to the construction of a linear perturbation model for reconstructing the difference in optical properties from changes in the data\cite{heiskala09,heiskalaPhD} with \cite{steinbrink2001} describing the use of  MC to recover optical absoprtion changes in a layered model

Because of the heavy computational cost of Monte Carlo methids, GPU based fast MC was investigated in\cite{alerstam:060504}. 
Another way to 'parallelize' computation of MC solutions with different optical parameters ($\mus$ and $\mua$) is to use white Monte Carlo. Absorption can always be scaled, but the problem of photon termination remains (if absorption is not known when photon packets are followed). In a homogeneous semi-infininite model also scattering can be scaled, but this is not applicable to heterogeneous models \cite{alerstam:041304}.

\subsection{Hybrid Methods}

\subsubsection{Coupled Radiative Transport and Diffusion}
To overcome the limitations of  diffusion theory in the proximity
of the light sources, hybrid methods which combine Monte Carlo
simulation with diffusion theory have been reported
\cite{wang93b,alexandrakis2000} and have been applied to 
media with low-scattering and non-scattering regions \cite{hayashi03}.
This approach still suffers from the time
consuming nature of the Monte Carlo methods and often require
iterative mapping between the models which increases computation times
even more. 

The Fokker--Planck equation, which has been found
to describe light propagation accurately when the scattering is
strongly peaked in the forward direction, can be utilized at the small
depths below highly collimated light sources \cite{kim2003}. 
However, it does not describe light propagation accurately 
at greater depths in biological tissues nor within low-scattering or
non-scattering regions.


A coupled radiative transfer equation and diffusion approximation
model was used in  \cite{bal2002a} for media
which contained strongly absorbing and low-scattering regions and
in\cite{tarvainen2005,tarvainen2005a} 
to describe light propagation in turbid
medium containing highly collimated light sources.
In\cite{tarvainen2005b} this model was extended to describe light
propagation in turbid medium with low-scattering and non-scattering
regions.  
In the coupled approach, the RTE is used as a forward model in sub-domains
in which the assumptions of the DA are not valid.
This includes the regions in the proximity of the source and
boundary and the low-scattering and non-scattering regions.
The DA is used as a forward model elsewhere in the domain.
The RTE and DA are coupled through boundary conditions between the RTE
and DA sub-domains and they are solved simultaneously using
the FEM.

\subsubsection{The Void Problem}
The void problem refers to the difficulty of modelling a non-scattering region within a highly scattering one. This has particular relevence in optical tomography applied to the brain where the presence of the cerebral spinal fluid (CSF) has just such a property~\cite{custo2006}.

The effect of void regions on particle transport was 
identified in early studies of
the Boltzmann equation for the propagation of neutrons 
in \cite{lathrop68,lathrop71,ackroyd86,ackroyd89}. 
In applications of optical tomography the effect of the voids has been studied more recently 
in \cite{firbank96,hielscher98,klose1999a,arridge2000a,dehghani2000a,bal2002a,bal2003,tarvainen2005b},
and it has been found that under these circumstances, 
%
most numerical solutions to the transport equation, apart from Monte Carlo, run into difficulties  so that special methods need to be developed.
One proposed solution is the \emph{radiosity-diffusion} model~\cite{firbank96,okada97,arridge2000a}. 
The principle here is to model propagation in voids through geometrical optics,
assuming that the distribution of directions is given by a \emph{non-local}
boundary condition at a diffuse/non-diffuse interface

\begin{eqnarray}
\hspace{-1.5cm}
\phidiff(\r) + 2\bccons\diffcoef\mneumop{\phidiff(\r)} = \frac{1}{\pi}\int_{\bvoid}\cos\theta\cos\theta'\left(\phidiff(\r')-
2\diffcoef\left(\frac{1+\refmm^{(1)}}{1-\refmm^{(0)}}\right)\mneumop{\phidiff(\r')}\right) \times \nonumber \\
(1-|\refmm(\theta)|^2)h(\r,\r')\frac{\exp[-(\mua + \rmi\frac{\omega}{c})|\r-\r'|]}{|\r-\r'|^2}\rmd\r' \\
\r,\r' \in \bvoid\,\quad
\cos\theta = \bnormvec(\r)\cdot \frac{\r-\r'}{|\r-\r'|}, \, \cos\theta'= \bnormvec(\r')\cdot \frac{\r-\r'}{|\r-\r'|} \nonumber
\label{eq:da_freq_bc_void}
\end{eqnarray}
where $\bvoid$ is the surface of the void region and $h(\r,\r')$ is
a \emph{visibility flag} that is unity if $\r,\r'$ are in line of sight across the void, and zero otherwise. See~\cite{ripollPhD} for
a detailed derivation. Both a FEM-radiosity \cite{arridge2000a} and BEM-radiosity
 \cite{ripoll2000a} implementation have been presented.
%
The radiosity-diffusion model was used in\cite{hyvonen2002,hyvonen2004,hyvonen2007a,hyvonen2007b} where a boundary recovery method was also developed for the inverse problem.

In \cite{bal2003, bal2006a} a different approach was proposed in which the void was replaced by an interface condition with an anisotropic Laplace-Beltrami operator. The tangential diffusion induced by this anisotropy was related to the void width and curvature.



\section{Direct Inversion Methods}
\label{sect:analytic_approaches}
By direct inversion we mean the use of inversion formulas and associated fast image reconstruction algorithms. In optical tomography, such formulas exist for particular experimental geometries, including those with planar, cylindrical and spherical boundaries. In this section, we consider several different inverse problems organized by length scale. We begin with the macroscopic case and proceed downward. 

\subsection{Diffuse optical tomography}

As described in section~\ref{sect:modalities}, the aim is to reconstruct the absorption and diffusion coefficients of a macroscopic medium from boundary measurements. We first consider the linearized inverse problem in one dimension and then discuss direct methods for linear and nonlinear inversion in higher dimensions. 

\subsubsection{One-dimensional problem}
\label{sect:1D}
We start by studying the time-dependent inverse problem in one dimension, which illustrates many features of the three-dimensional case. Let $\Omega$ be the half-line $x\ge 0$. The energy density $\Phi$ obeys the one-dimensional version of the time-dependent diffusion equation~\eref{eq:da_time}
\begin{equation}
\frac{\partial }{\partial t}\Phi(x,t) = D \frac{\partial^2}{\partial x^2}\Phi(x,t) - c\mu_a(x)\Phi(x,t) \ , \quad x\in\Omega \ ,
\end{equation}
where the diffusion coefficient $D$ is assumed to be constant, an assumption that will be relaxed later. The energy density is taken to obey the initial and boundary conditions
\begin{eqnarray}
\Phi(x,0)=\delta(x-x_1) \ , \\
\Phi(0,t) - \ell_{\rm ext} \frac{\partial \Phi}{\partial x}(0,t) = 0 \ .
\end{eqnarray}
Here the initial condition imposes the presence of a point source of unit-strength at $x_1$.
Since $\Phi$ decays exponentially, we consider for $k\ge0$ the Laplace transform
\begin{equation}
\Phi(x,k) = \int_0^\infty e^{-k^2 D t} \Phi(x,t) dt \ ,
\end{equation}
which obeys the equation
\begin{equation}
\label{1D_diff_eqn}
-\frac{d^2\Phi(x)}{dx^2} + k^2(1+\eta(x)) \Phi(x) =\frac{1}{D}\delta(x-x_1) \ ,
\end{equation}
where $\eta$ is the spatially-varying part of the absorption, which is defined by $\eta=c\mu_a/(Dk^2)-1$.
The solution to the forward problem is given by the integral equation
\begin{equation}
\label{1D_int_eq}
\Phi(x) = \Phi_i(x) - k^2 \int_\Omega G(x,y) \Phi(y)\eta(y) dy \ ,
\end{equation} 
where the Green's function is of the form
\begin{equation}
\label{1D_greens_function}
G(x,y) = \frac{1}{2Dk}\left(e^{-k|x-y|} + \frac{1-k\ell_{\rm ext}}{1+k\ell_{\rm ext}} e^{-k|x+y|}\right) \ ,
\end{equation}
and $\Phi_i$ is the incident field, which obeys~\eref{1D_diff_eqn} with
$\eta=0$. The above integral equation may be linearized with respect to
$\eta(x)$ by replacing $u$ on the right-hand side by $u_i$. This approximation is accurate when ${\rm supp}(\eta)$ and $\eta$ are small. If we introduce the scattering data
$\Phi_s=\Phi_i-\Phi$ and perform the above linearization, we obtain
\begin{equation}
\label{1D_int_eq_linearized}
\Phi_s(x_1,x_2) = k^2 \int_\Omega G(x_1,y)G(y,x_2) \eta(y) dy \ .
\end{equation}
Here $\Phi_s(x_1,x_2)$ is proportional to the change in intensity due to a point source at $x_1$ that is measured by a detector at $x_2$. 

In the backscattering geometry, the source and detector are placed at the origin ($x_1=x_2=0$) and \eref{1D_int_eq_linearized} becomes, upon using \eref{1D_greens_function} and omitting overall constants
\begin{equation}
\label{def_A}
\Phi_s(k) = \int_0^\infty e^{-kx} \eta(x) dx  \ ,
\end{equation}
where the dependence of $\Phi_s$ on $k$ has been made explicit. Thus, the linearized inverse problem can be seen to correspond to inverting the Laplace transform of $\eta$. Inversion of the Laplace transform is the paradigmatic exponentially ill-posed problem~which can be analyzed as follows~\cite{Epstein_2008_1} . Eq.~\eref{def_A} defines an operator $A: \eta \mapsto \Phi_s$ which is bounded and self-adjoint on $L^2([0,\infty])$. The singular functions $f$ and $g$ of $A$ satisfy
\begin{equation}
A^*A f = \sigma^2 f \ , \ \ \ AA^* g = \sigma^2 g \ ,
\end{equation}
where $\sigma$ is the corresponding singular value.
In addition, $f$ and $g$ are related by
\begin{equation}
Af=\sigma g \ , \ \ \ A^*g=\sigma f \ .
\end{equation}
If we observe that $A^*A(x,y) = 1/(x+y)$ and use the identity
\begin{equation}
\int_0^\infty \frac{y^{a}}{1+y} dy = \frac{\pi}{\sin(a + 1)\pi} \ , \ \ \  
-1 \le \mathrm{Re}(a) <0 \ ,
\end{equation}
we see that
\begin{equation}
f_s(x) = g_s^*(x) = {1\over \sqrt{2\pi}} x^{-{1\over 2} + i s} \ , 
\ \ \ s\in\mathbb R 
\end{equation}
and 
\begin{equation}
\sigma_s^2 = {\pi\over \cosh(\pi s)} \sim e^{-\pi |s|} \ .
\end{equation}
Note that the singular values of $A$ are exponentially small, which gives rise to severe ill-posedness. Using the above, we can write an inversion formula for 
\eref{1D_greens_function} in the form
\begin{eqnarray}
\eta(x) = \int_0^\infty dk \int_{-\infty}^\infty ds R\left(\frac{1}{\sigma_s}\right) f_s(x) g_s^*(k) \Phi_s(k)  \ ,
\end{eqnarray}
where the regularizer $R$ has been introduced to control the contribution of small singular values.

\subsubsection{Linearized inverse problem}
\label{sec:linearized3D}
We now consider the linearized inverse problem in three dimensions. To indicate the parallels with the one-dimensional case, we will find it convenient to work in the half-space geometry. Extensions to other geometries, including those with planar, cylindrical and spherical boundaries is also possible~\cite{Schotland_1997_1,Schotland_2001_1,Markel_2004_1}. In particular, we note that the slab geometry is often employed in optical mammography and small animal imaging. 
As before, we define the scattering data $\Phi_s=\Phi_i-\Phi$. Linearizing~\eref{LS_diff} with respect to $\delta\mu_a$ and $\delta D$ we find that, up to an overall constant, $\Phi_s$ obeys the integral equation
\begin{equation}
\label{phi_def}
\Phi_s({\bm r}_1,{\bm r}_2) = \int d\Br \left[G(\Br_1,\Br)G(\Br,\Br_2)c\delta\mu_a(\Br) +
\grad_{\Br}G(\Br_1,\Br)\cdot\grad_{\Br}G(\Br,\Br_2) \delta D(\Br) \right]  \ ,
\end{equation}
where $\Br_1$ is the position of the source, $\Br_2$ is the position of the detector and we have integrated by parts to evaluate the action of the operator $V$.
The Green's function in the half-space $z\ge 0$ is given by the plane-wave decomposition
\begin{equation}
\label{G_0}
G(\Br,\Br')=\int \frac{d\Bq }{(2\pi)^2} g(\Bq; z,z')
\exp[i\Bq\cdot(\Brho-\Brho')] \ ,
\end{equation}
where we have used the notation $\Br=(\Brho,z)$. If either $\Br$ or $\Br'$ lies in the plane $z=0$, then
\begin{equation}
\label{g_def}
g({\bm q}; z, z^{\prime}) = \frac{\ell_{\rm ext}}{D_0}{{\exp\left[ -Q({\bm q})\vert z - z^{\prime}\vert \right]} \over {Q({\bm q})\ell_{\rm ext} + 1}} \ ,
\end{equation}
where
\begin{equation}
Q(\Bq) = \sqrt{q^2 + k^2} \ .
\end{equation}
The inverse problem is to recover $\delta\alpha$ and $\delta D$ from boundary measurements. To proceed, we introduce the Fourier transform of $\Phi_s$
with respect to the source and detector coordinates according to
\begin{equation}
\label{Fourier_data}
\tilde\Phi_s(\Bq_1,\Bq_2) = \int d\Brho_1 d\Brho_2 e^{i(\Bq_1\cdot\Brho_1 +\Bq_2\cdot\Brho_2)}
\Phi_s(\Brho_1,0;\Brho_2,0) \ .
\end{equation}
If we define 
\begin{equation}
\psi(\Bq_1,\Bq_2)=(Q(\Bq_1)\ell_{\rm ext} +1)(Q(\Bq_2)\ell_{\rm ext} +1)\tilde\Phi_s(\Bq_1,\Bq_2) 
\end{equation}
and make use of \eref{G_0} and \eref{phi_def}, we find that
\begin{equation}
\label{psi_diff}
\psi(\Bq_1,\Bq_2) = \int d\Br e^{i(\Bq_1+\Bq_2)\cdot\Brho}\left[\kappa_A(\Bq_1,\Bq_2;z)c\delta\mu_a(\Br) + \kappa_D(\Bq_1,\Bq_2;z)\delta D(\Br)\right] \ .
\end{equation}
Here 
\begin{eqnarray}
\kappa_A(\Bq_1,\Bq_2;z)=c\exp[-\left(Q(\Bq_1)+Q(\Bq_2)\right)z] \ , \\
\kappa_D(\Bq_1,\Bq_2;z)=-\left(\Bq_1\cdot\Bq_2+Q(\Bq_1)Q(\Bq_2)\right)\exp[-\left(Q(\Bq_1)+Q(\Bq_2)\right)z] \ .
\end{eqnarray}
We now change variables according to
\begin{equation}
\label{FL_variables}
{\bm q}_1={\bm q}+{\bm p}/2 \ , \ \ \ 
{\bm q}_2={\bm q}-{\bm p}/2 \ , 
\end{equation}
where ${\bm q}$ and ${\bm p}$ are independent
two-dimensional vectors and rewrite (\ref{psi_diff}) as
\begin{equation}
\label{Fourier-Laplace_diff}
\psi({\bm q}+{\bm p}/2,{\bm q}-{\bm p}/2) = \int d\Br\exp(-i{\bm
  q}\cdot \Brho)\left[ \kappa_A({\bm p},{\bm q};z) \delta
  \mu_a({\bm r}) + \kappa_D({\bm p},{\bm q};z)\delta D(\bm r)\right]
\ , 
\end{equation}
The above result has the structure of a Fourier-Laplace transform by which $\delta\mu_a$ and $\delta D$ are related to $\psi$.  This relation 
can be used to obtain an inversion formula for the integral equation \eref{psi_diff}. To proceed, we note that the Fourier-transform in the transverse direction can be inverted separately from the Laplace transform in the longitudinal direction.  We thus arrive at the result
\begin{equation}
\label{main_1D}
\psi({\bm q}+{\bm p}/2,{\bm q}-{\bm p}/2) = \int \left[
  \kappa_A({\bm p},{\bm q};z) \tilde\delta\mu_a({\bm q},z) + \kappa_D({\bm p},{\bm q};z) \tilde\delta D({\bm q},z) \right]dz \ , 
\end{equation}
where $\tilde\delta\mu_a$ and $\tilde\delta D$ denote the two-dimensional Fourier-transform.
For fixed ${\bm q}$, Eq.~(\ref{main_1D}) defines an integral
equation for $\tilde\delta\mu_a({\bm q},z)$ and $\tilde\delta D({\bm q},z)$. It is readily seen that the minimum $L^2$ norm solution to (\ref{main_1D}) has the form
\begin{eqnarray}
\label{a_inv}
\fl
\tilde\delta\mu_a({\bm q},z) = \int d\Bp d\Bp^{\prime} \kappa_A^*({\bm p},{\bm q};z)
M^{-1}(\Bp,\Bp';\Bq) \psi({\bm p}^{\prime}+{\bm
  q}/2,{\bm p}^{\prime} -{\bm q}/2) \ , \\
\label{b_inv}
\fl
\tilde\delta D({\bm q},z) = \int d\Bp d\Bp^{\prime} \kappa_D^*({\bm p},{\bm q};z)
M^{-1}(\Bp,\Bp';\Bq)  \psi({\bm p}^{\prime}+{\bm
  q}/2,{\bm p}^{\prime} - {\bm q}/2) \ ,
\end{eqnarray}
where the matrix elements of $M$ are given by the integral
\begin{equation}
\label{M_def}
 M(\Bp,\Bp';\Bq) =
\int_0^L \left[ \kappa_A({\bm p},{\bm q};z)\kappa_A^*({\bm
  p}^{\prime},{\bm q};z)  +  \kappa_D({\bm p},{\bm q};z)\kappa_D^*({\bm
  p}^{\prime},{\bm q};z) \right] dz
\end{equation}
Finally, we apply the inverse Fourier transform in the transverse direction to arrive at the inversion formula
\begin{eqnarray}
\fl
\delta \mu_a({\bm r}) = \int{ d\Bq  
\over {(2\pi)^2}} e^{-i{\bm q}\cdot\Brho}
\int d\Bp d\Bp^{\prime} \kappa_A^*({\bm p},{\bm q};z)M^{-1}(\Bp,\Bp';\Bq)  \psi({\bm q}^{\prime}+{\bm
  p}/2,{\bm q}^{\prime} -{\bm p}/2) \ ,
\label{alpha_inv} \\
\fl
\delta D({\bm r}) = \int { d\Bq  
\over {(2\pi)^2}} e^{-i{\bm q}\cdot\Brho}
\int d\Bp d\Bp^{\prime} \kappa_D^*({\bm p},{\bm q};z) M^{-1}(\Bp,\Bp';\Bq) \psi({\bm q}^{\prime}+{\bm p}/2,{\bm q}^{\prime} - {\bm p}/2) \ .
\label{D_inv}
\end{eqnarray}

Several remarks on the above result are necessary. First, implementation of~\eref{alpha_inv} and \eref{D_inv} requires regularization to stabilize the computation of the inverse of the matrix $M$. Second, sampling of the data function $\Phi_s$ is easily incorporated. The Fourier transform in~\eref{Fourier_data} is replaced by a lattice Fourier transform. If the corresponding wave vectors $\Bq,\Bq_2$ are restricted to the first Brillouin zone of the lattice, then the inversion formula \eref{alpha_inv} and \eref{D_inv} recovers a bandlimited approximation  to the coefficients $\delta\mu_a$ and $\delta D$~\cite{Markel_2004_1,Markel_2002_1}. Third, the resolution of reconstructed images in the transverse and longitudinal directions is, in general, quite different. The transverse resolution is controlled by sampling and is determined by the highest spatial frequency that is present in the data. The longitudinal resolution is much lower due to the severe ill-posedness of the Laplace transform inversion which is implicit in~\eref{Fourier-Laplace_diff}, similar to the one-dimensional case discussed in section~\ref{sect:1D}. Finally, the inversion formula \eref{alpha_inv} and \eref{D_inv} can be used to develop a fast image reconstruction algorithm whose computational complexity scales as $O(NM\log M)$, where $M$ is the number of detectors, $N$ is the number of sources and $M\gg N$~\cite{Markel_2004_1}. The algorithm has recently been tested in noncontact optical tomography experiments. Quantitative reconstructions of complex phantoms with millimeter-scale features located centimeters within a highly-scattering medium have been reported~\cite{Wang_2005_1,Konecky_2008_1}. Data sets of order $10^8$ source-detector pairs can be reconstructed in approximately 1 minute of CPU time on a 1.5 GHz computer.

\subsubsection{Boundary removal}

The direct inversion method described in section~\ref{sec:linearized3D} is applicable to relatively simple geometries, including those with planar, cylindrical and spherical boundaries. An extension of the method can be employed to handle more complex geometries~\cite{Ripoll_2006_1}. Consider a bounded domain $\Omega$ in which the diffusion equation~\eref{diff} for the energy density $\Phi$ holds. Making use of the Kirchoff integral as described in section~5.2, we obtain
\begin{eqnarray}
\label{kirchoff}
\nonumber
\Phi(\Br) = && \Phi_0(\Br) - c\int_\Omega G(\Br,\Br')\delta\mu_a(\Br') \Phi(\Br')d\Br' \\ && + D \int_{\partial\Omega} \left[G(\Br,\Br')\grad \Phi(\Br')-\Phi(\Br')\grad_{\Br'} G(\Br,\Br') \right]
\cdot \n d\Br' \ ,
\end{eqnarray}
where, for simplicity, we have assumed that the diffusion coefficient is constant throughout $\Omega$.
Here
\begin{equation}
\Phi_0(\Br)=\int_\Omega G(\Br,\Br')J_{-}(\Br') d\Br' 
\end{equation} 
and $G$ is the fundamental solution to the diffusion equation defined in \eref{diff_fund}.
Using the boundary condition~\eref{eq:da_time_bc} and introducing the approximation
$\Phi=\Phi_0$ in $\Omega$, we find that
\eref{kirchoff} becomes
\begin{equation}
\label{boundary_data}
(\Phi_0-\Phi-M)(\Br) = \int_\Omega G(\Br,\Br')\Phi_0(\Br') c\delta\mu_a(\Br') d\Br' \ ,
\end{equation}
where 
\begin{equation}
M(\Br) = \int_{\partial\Omega}\left[\frac{1}{\ell_{\rm ext}} G(\Br,\Br') + \n\cdot\grad_{\Br'} G(\Br,\Br')  \right]\Phi(\Br')d\Br' \ .
\end{equation}
We observe that the left hand side of~\eref{boundary_data} is known from measurements. Thus, within the accuracy of the aforementioned approximation, the inverse problem becomes that of recovering $\delta\alpha$ in an infinite medium with measurements on the surface $\partial\Omega$. Let $B$ denote a ball that contains $\Omega$. Then the field 
on $\partial\Omega$ can be propagated to $\partial B$ by using the Rayleigh-Sommerfeld formula
\begin{equation}
\Phi(\Br) = D\int_{\partial\Omega} \n\cdot\grad_{\Br'} G(\Br,\Br') \Phi(\Br') d\Br' \ ,
\quad \Br \in \partial B \ .
\end{equation}
Thus, it is possible to reformulate the inverse problem for a bounded domain in terms of one for an infinite domain, with measurements prescribed on a spherical surface. This latter problem is one for which direct inversion can be performed. 

\subsubsection{Nonlinear inverse problem}

We now consider the nonlinear inverse problem of DOT. 
The Born series (\ref{pert_V}) can be rewritten in the form
\begin{equation}
\label{series_eta}
\Phi_s(\Br_1,\Br_2)=\int d\Br K_1^i(\Br_1,\Br_2;\Br)\eta_i(\Br) +
\int d\Br d\Br' K_2^{ij}(\Br_1,\Br_2;\Br,\Br')\eta_i(\Br)\eta_j(\Br') + \cdots \ ,
\end{equation}
where 
\begin{equation}
\label{eta_def}
\eta({\bm r}) = 
\left( \begin{array}{l}
\eta_1({\bm r}) \\
\eta_2({\bm r}) \\
\end{array}\right) = 
\left( \begin{array}{l}
c\delta\mu_a({\bm r}) \\
\delta D({\bm r}) \\
\end{array}\right) \ ,
\end{equation}
summation over repeated indices is implied with $i,j=1,2$. The components of 
the operators $K_1$ and $K_2$ are given by
\begin{eqnarray}
\label{K_1_1}
&& K_1^1({\bm r}_1,{\bm r}_2; {\bm r}) = G({\bm r}_1,{\bm
  r})G({\bm r},{\bm r}_2) \ , \\
\label{K_1_2}
&& K_1^2({\bm r}_1,{\bm r}_2; {\bm r}) = 
\nabla_{{\bm r}} G({\bm r}_1,{\bm r}) 
\cdot \nabla_{{\bm r}} G({\bm r},{\bm r}_2) \ , \\
\label{K_2_11}
&& K_2^{11}({\bm r}_1,{\bm r}_2; {\bm r},{\bm r}') = -G({\bm r}_1,{\bm
  r})G({\bm r},{\bm r}')G({\bm r}',{\bm r}_2) \ , \\
\label{K_2_12}
&& K_2^{12}({\bm r}_1,{\bm r}_2; {\bm r},{\bm r}') = - G({\bm r}_1,{\bm
  r})\nabla_{{\bm r}'}G({\bm r},{\bm r}') \cdot
\nabla_{{\bm r}'} G({\bm r}',{\bm r}_2)  \ , \\
\label{K_2_21}
&& K_2^{21}({\bm r}_1,{\bm r}_2; {\bm r},{\bm r}') = -\nabla_{{\bm r}} G({\bm r}_1,{\bm r}) \cdot\nabla_{{\bm r}}G({\bm r},{\bm r}') G({\bm r}',{\bm r}_2)  \ , \\ 
\label{K_2_22}
&& K_2^{22}({\bm r}_1,{\bm r}_2; {\bm r},{\bm r}') = -\nabla_{{\bm r}} G({\bm r}_1,{\bm r}) \cdot\nabla_{{\bm r}}
 \left[\nabla_{{\bm r}'}G({\bm r},{\bm r}') \cdot 
   \nabla_{{\bm r}'} G({\bm r}',{\bm
    r}_2)\right]   \ .
\end{eqnarray}
It will prove useful to express the Born series as a formal power series in tensor powers of $\eta$ of the form
\begin{equation}
\label{phi_def_again}
\Phi_s = K_1\eta + K_2 \eta\otimes\eta + K_3 \eta\otimes\eta\otimes\eta 
+ \cdots \ .
\end{equation}

The solution to the nonlinear inverse problem of DOT may be formulated from the ansatz that $\eta$ may be expressed as a series in tensor powers of $\Phi_s$ of the form
\begin{equation}
\label{inverse_series}
\eta = \K_1 \Phi_s + \K_2\Phi_s\otimes\Phi_s + \K_3\Phi_s\otimes\Phi_s\otimes\Phi_s + \cdots  \ , 
\end{equation}
where the $\K_j$'s are operators which are to be determined~\cite{Moses_1956_1,Prosser_1969_1,Markel_2003_1,Moskow_2008_1}.
To proceed, we substitute the expression { (\ref{phi_def_again}) for
$\Phi_s$  into  (\ref{inverse_series})}  and equate terms with the same tensor 
power of $\eta$. We thus obtain the relations
\begin{eqnarray}
\K_1K_1 & = & I \ , \label{ident}\\
\K_2K_1 \otimes K_1 + \K_1K_2 & = & 0 \ , \\
\K_3 K_1 \otimes K_1 \otimes K_1 + \K_2K_1\otimes K_2 + \K_2K_2 \otimes K_1 + \K_1K_3 
& = & 0  \ , \\
\sum_{m=1}^{j-1} \K_m \sum_{i_1+\cdots+i_m=j} K_{i_1} \otimes \cdots \otimes K_{i_m} + \K_j K_1 \otimes \cdots \otimes K_1 &=& 0 \  ,
\end{eqnarray}
which may be solved for the $\K_j$'s with the result
\begin{eqnarray}
\K_1 &=& K_1^+ \ , \\
\K_2 &=& -\K_1K_2\K_1\otimes\K_1 \ , \\
\K_3 &=&  -\left(\K_2K_1\otimes K_2 +
\K_2K_2\otimes K_1+\K_1K_3\right)\K_1\otimes\K_1\otimes\K_1 \ , \\
\label{def_K_j}
\K_j &=& - \left(\sum_{m=1}^{j-1} \K_m \sum_{i_1+\cdots+i_m=j} K_{i_1} \otimes \cdots \otimes
K_{i_m}\right) \K_1 \otimes \cdots \otimes \K_1 \ .
\label{def_kappa}
\end{eqnarray}

\begin{figure}[t]
\centering
\includegraphics[width=.7\textwidth,trim=-10 470 -10 -10,keepaspectratio=true]{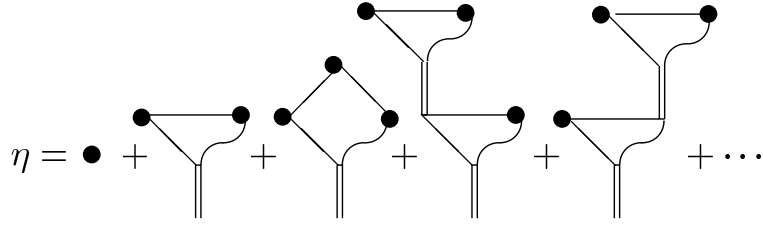} 
\caption{Diagrammatic representation of the inverse scattering series.}
\label{fig:diagrams}
\end{figure}

We will refer to (\ref{inverse_series}) as the inverse series. Here we note several of its properties. First, $K_1^+$ is the regularized pseudoinverse of the operator $K_1$.  The singular value decomposition of the operator $K_1^+$ can be computed analytically for particular geometries, as explained in section~3.5. Since the operator $\K_1$ is unbounded,  regularization of $K_1^+$ is required to control the ill-posedness of the inverse problem. Second, the coefficients in the inverse series have a recursive structure. The operator $\K_j$ is determined by the coefficients of the Born series $K_1, K_2,\ldots, K_j$. Third, the inverse scattering series can be represented in diagrammatic form as shown in Figure~\ref{fig:diagrams}. A solid line corresponds to a factor of $G$, a wavy line to the incident field, a solid vertex ($\bullet$) to $\K_1\Phi_s$, and the effect of the final application of $\K_1$ is to join the ends of the diagrams. Note that the recursive structure of the series is evident in the diagrammatic expansion which is shown to third order. 
Finally, inversion of only the linear term in the Born series is required to compute the inverse series to all orders. Thus an ill-posed nonlinear inverse problem is reduced to an ill-posed linear inverse problem plus a well-posed nonlinear problem, namely the computation of the higher order terms in the series.

We now characterize the convergence of the inverse series. We restrict our attention to 
the case of a uniformly scattering medium for which $\eta=\delta\alpha$. 
To proceed, we require an estimate on the $L^2$ norm of the operator $\K_j$. We define
the constants $\mu$ and $\nu$ by
\begin{equation}
\label{def_mu2}
\mu=\sup_{\Br\in B_a} k^2 \|G_0(\Br,\cdot)\|_{L^2(B_a)} \ .
\end{equation}
 
\begin{equation}
\label{def_nu2}
\nu = k^2 |B_a|^{1/2} \sup_{\Br\in B_a}\|G_0(\Br,\cdot)\|_{L^2(\partial\Omega)}  \ .
\end{equation}
Here $B_a$ denotes a ball of radius $a$ which contains the support of $\eta$.
It can be shown~\cite{Moskow_2008_1} that if $(\mu_p+\nu_p)\|\K_1\|_2 < 1$ then the operator
\begin{equation}
\K_j : L^2(\partial\Omega\times\cdots\times \partial\Omega) \longrightarrow L^2(B_a)
\end{equation}
defined by (\ref{def_kappa}) is bounded and 
\begin{equation}
\|\K_j \|_2 \le C(\mu+\nu)^j\|\K_1\|_2^j  \ ,
\end{equation}
where $C=C(\mu_p,\nu_p,\| \K_1 \|_p)$ is independent of $j$.

We can now state the main result on the convergence of the inverse series.

\begin{theorem}\cite{Moskow_2008_1}.
\label{thm1}
The inverse scattering series converges in the $L^2$ norm  if $\|\K_1\|_2 < 1/(\mu+\nu)$ and $\|\K_1 \Phi_s\|_{L^2(B_a)} < 1/(\mu+\nu)$. Furthermore, the following estimate for the series limit $\tilde\eta$ holds
\begin{equation*}
\Big\|\tilde\eta-\sum_{j=1}^N \K_j \Phi_s\otimes\cdots\otimes\Phi_s\Big\|_{L^p(B_a)} \le C \frac{\left[(\mu_p+\nu_p)\|\K_1\Phi_s \|_{L^2(B_a)}\right]^{N+1}}{1-(\mu+\nu)\|\K_1\Phi_s \|_{L^2(B_a)}} \ ,
\end{equation*}
where $C=C(\mu,\nu, \| \K_1\|_2 )$  does not depend on $N$ nor on the scattering data $\Phi_s$.
\end{theorem}

The stability of the limit of the inverse series under perturbations in the scattering data can be analyzed as follows: 

\begin{theorem}\cite{Moskow_2008_1}.
\label{thm2}
Let  $\|\K_1\|_2 < 1/(\mu+\nu)$ and let $\Phi_{s1}$ and $\Phi_{s2}$ be scattering data for which $M\|\K_1\|_2  < 1/(\mu+\nu)$, where
$M=\max{(\|\Phi_{s1}\|_2,\| \Phi_{s2}\|_2)}$.
Let $\eta_1$ and $\eta_2$ denote the corresponding limits of the inverse scattering series. Then the following estimate holds
\begin{equation*}
\|\eta_1-\eta_2\|_{L^2(B_a)} < \tilde{C}  \|\Phi_{s1}-\Phi_{s2}\|_{L^2(\partial\Omega\times\partial\Omega)}  \ ,
\end{equation*}
where $\tilde{C}=\tilde{C}(\mu,\nu,\| \K_1\|_2,M)$ is a constant that is otherwise independent of $\Phi_{s1}$ and $\Phi_{s2}$.
\end{theorem}

It is a consequence of the proof of Theorem~\ref{thm2} that $\tilde C$ is proportional to $\|\K_1\|_2$. Since regularization sets the scale of $\|\K_1\|_2$, it follows that the stability of the nonlinear inverse problem is controlled by the stability of the linear inverse problem.

The limit of the inverse scattering series does not, in general, coincide with $\eta$. We characterize the approximation error as follows.

\begin{theorem}\cite{Moskow_2008_1}
\label{error_estimate}
Suppose that $\| \K_1 \|_2< 1/(\mu +\nu)$, $\|\K_1\Phi_s\|_{L^2(B_a)} < 1/(\mu+\nu)$. Let $\mathcal{M}=\max{( \| \eta\|_{L^2(B_a)} , \| \K_1 K_1 \eta\|_{L^2(B_a)} )}$ and assume that $\mathcal{M} < 1/(\mu+\nu)$. Then the norm of the difference between the partial sum of the inverse series and the true absorption obeys the estimate
\begin{equation*}
\fl
\Big\|\eta-\sum_{j=1}^N\K_j \Phi_s\otimes\cdots\otimes\Phi_s\Big\|_{L^2(B_a)} \le
C \|(I-\K_1 K_1)\eta \|_{L^2(B_a)} + 
\tilde C{[(\mu_p+\nu_p)\|\K_1\|_2\|\Phi_s\|]^N \over{ 1-(\mu+\nu)\|\K_1\|_2\|\Phi_s\|_2 }}\ ,
\end{equation*}
where $C=C(\mu,\nu,\| \K_1\|_2, \mathcal{M})$ and $\tilde C=\tilde C (\mu,\nu\| \K_1\|_2)$ are independent of $N$ and $\Phi_s$.
\end{theorem}

We note that, as expected, the above result shows that regularization of $\K_1$ creates an error in the reconstruction of $\eta$.  For a fixed regularization, 
the relation $\K_1 K_1=I$ holds on a subspace of $L^2(B_a)$ which, in practice, is finite dimensional. By regularizing $\K_1$ more weakly, the subspace becomes larger, eventually approaching all of $L^2$. However, in this instance, the estimate in Theorem~\ref{error_estimate} would not hold since $\|\K_1\|_2$ is so large that the inverse scattering series would not converge.  Nevertheless, Theorem~\ref{error_estimate} does describe what can be reconstructed exactly, namely those $\eta$ for which $\K_1 K_1 \eta =I$. That is, if we know apriori that $\eta$ belongs to a particular finite-dimensional subspace of $L^2$, we can choose $\K_1$ to be a true inverse on this subspace. Then, if $\| \K_1\|_2 $ and $\|\K_1 \Phi_s\|_{L^2} $ are sufficiently small, the inverse series will recover $\eta$ exactly. 

It is straightforward to compute the constants $\mu$ and $\nu$ in dimension three. We have
\begin{eqnarray}
\label{mu}
\mu = k^2 e^{-ka/2} \left(\frac{\sinh(ka)}{4\pi k}\right)^{1/2} \ .
\end{eqnarray}
\begin{eqnarray}
\label{nu}
\nu\le  k^2 |\partial\Omega| |B_a|^{1/2}  \frac{e^{-2k \dist(\partial\Omega, B_a)}}{(4\pi\dist(\partial\Omega, B_a))^2} \ .
\end{eqnarray}
Note that $\nu$ is exponentially small.
It can be seen that the radius of convergence of the inverse series $R=1/(\mu+\nu) \lesssim O(1/(ka)^{3/2})$ when $ka \gg 1$.

Numerical studies of the inverse scattering series have been performed. For inhomogeneities with radial symmetry, exact solutions to the forward problem were used as scattering data and reconstructions were computed to fifth order in the inverse series. It was found that the series appears to converge quite rapidly for low contrast objects. As the contrast is increased, the higher order terms systematically improve the reconstructions until, at sufficiently large contrast, the series diverges.

\subsection{Inverse transport}

We now turn our attention to the inverse problem for the RTE. This is a very large subject in its own right. Indeed, a recent topical review has already covered the key mathematical issues regarding existence, uniqueness and stability of the inverse transport problem~\cite{bal2009a}. Here we aim to discuss the inverse problem in some special cases which lead to direct inversion procedures and fast algorithms. See~\cite{Schotland_2007_1} for further details.
As in section~\ref{sec:linearized3D}, we will work in the $z\ge0$ half-space with the source and detector located on the $z=0$ plane. The source is assumed to be pointlike and oriented in the inward normal direction. The light exiting the medium is further assumed to pass through a normally-oriented angularly-selective aperture which collects all photons with intensity
\begin{equation}
\label{avg_I}
{\mathcal I}(\Br)=\int_{\n\cdot\svec>0} \n\cdot\svec A(\svec) I(\Br,\svec) d\svec\ ,
\end{equation}
where $A$ accounts for the effect of the aperture and the integration is carried out over all outgoing directions. When the aperture selects only photons traveling in the normal direction, then $A(\svec)=\delta(\svec-\n)$ and ${\mathcal I}(\Br) = I(\Br,\n)$. This case is relevant to noncontact measurements in which the lens is focused at infinity. The case of complete angularly-averaged data corresponds to $A(\svec)=1$. 

We now consider the linearized inverse problem. If the medium is inhomgeneously absorbing, it follows from (\ref{phi}) and (\ref{avg_I}) that the change in intensity measured relative to a homogeneous reference medium with absorption $\bar\mu_a$ is proportional to the data function $\Phi_s(\Brho_1,\Brho_2)$ which obeys the integral equation
\begin{equation}
\Phi_s(\Brho_1,\Brho_2)=\int_{\n\cdot\svec>0} \n\cdot\svec A(\svec) \Phi_s(\Brho_1,0,\z;\Brho_2,0,-\svec) d\svec  \ .
\end{equation}
Following the development in section~\ref{sec:linearized3D}, we consider the Fourier transform of $\Phi_s$ with respect to the source and detector coordinates. Upon inserting the plane-wave decomposition for $G$ given by (\ref{Green_RTE}) into (\ref{phi}) and carrying out the Fourier transform, we find that
\begin{eqnarray}
\label{FL_RTE}
\tilde\Phi_s(\Bq_1,\Bq_2)=\sum_{\mu_1,\mu_2}M_{\mu_1\mu_2}(\Bq_1,\Bq_2)\int d\Br && \exp\big[i(\Bq_1+\Bq_2)\cdot\Brho \\ && - (Q_{\mu_1}(\Bq_1) + Q_{\mu_2}(\Bq_2)z\big]\delta\mu_a(\Br) \ ,
\nonumber
\end{eqnarray}
where
\begin{eqnarray}
\label{def_M}
\fl
M_{\mu_1\mu_2}(\Bq_1,\Bq_2)&=&\sum_{l_1m_1,l_1'm_1'}\sum_{l_2m_2}{\mathcal B}_{l_1'm_1'}^{l_1m_1}(\Bq_1,\mu_1) {\mathcal B}_{l_2m_2}^{l_1'm_1'}(\Bq_2,\mu_2)\int_{\n\cdot\svec>0} \n\cdot\svec A(\svec) Y_{l_2m_2}(\svec) d\svec .
\end{eqnarray}
Eq.~(\ref{FL_RTE}) is a generalization of the Fourier-Laplace transform which holds for the diffusion approximation, as discussed in section~\ref{sec:linearized3D}. It can be seen that (\ref{FL_RTE}) reduces to the appropriate form in the diffuse limit, since only the smallest discrete eigenvalue contributes.

The inverse problem now consists of recovering $\delta\mu_a$ from $\tilde\Phi_s$. To proceed, we make use of the change of variables~\eref{FL_variables} and rewrite (\ref{FL_RTE}) as
\begin{equation}
\label{Phi}
\tilde\Phi_s(\Bq+\Bp/2,\Bq-\Bp/2)=\int dz K(\Bq,\Bp;z)\tilde\delta\mu_a(\Bq,z) \ ,
\end{equation}
where $\tilde\delta\mu_a(\Bq,z)$ denotes the two-dimensional Fourier transform of $\delta\mu_a$ with respect to its transverse argument and
\begin{eqnarray}
K(\Bq,\Bp;z)=&&\sum_{\mu_1,\mu_2}M_{\mu_1\mu_2}(\Bq + \Bp/2,\Bq - \Bp/2) \\ 
&&\times\exp\left[-\left(Q_{\mu_1}(\Bq + \Bp/2)+Q_{\mu_2}(\Bq - \Bp/2)\right)z\right] \ .
\end{eqnarray}
This change of variables can be used to separately invert the transverse and longitudinal functional dependences of $\delta\mu_a$ since for fixed $\Bq$, (\ref{Phi}) defines a one-dimensional integral equation for $\tilde\delta\mu_a(\Bq,z)$ whose pseudoinverse solution can in principle be computed. We thus obtain a solution to the inverse problem in the form
\begin{equation}
\delta\mu_a(\Br)=\int\frac{d\Bq}{(2\pi)^2}e^{-i\Bq\cdot\Brho}\int d\Bp K^+(z;\Bq,\Bp)\tilde\Phi_s(\Bq+\Bp/2,\Bq-\Bp/2)\ ,
\end{equation}
where $K^+$ denotes the pseudoinverse of $K$. It is important to note the ill-posedness due to the exponential decay of the evanescent modes (\ref{modes_general}) for large $z$. Therefore, we expect that the resolution in the $z$ direction will degrade with depth, but that sufficiently close to the surface the transverse resolution will be controlled by sampling.

\subsection{Single-scattering tomography}
\label{sec:SSOT}

Consider an experiment in which a narrow collimated beam is normally incident on a highly-scattering medium which has the shape of a slab.
Suppose that the slab is sufficiently thin that the incident beam undergoes predominantly single-scattering and that the intensity of transmitted light is measured by an angularly-selective detector. If the detector is collinear with the incident beam and its aperture is set to collect photons traveling in the normal direction, then only unscattered  photons will be measured. Now, if the aperture is set away from the normal direction, then the detector will not register any photons. However, if the detector is no longer collinear and only photons which exit the slab at a fixed angle are collected, then it is possible to selectively measure single-scattered photons. Note that the contribution of single-scattered photons is described by~\eref{delta_final}.

The inverse problem for single-scattered light is to recover $\mu_a$ from measurements of $\delta I$ as given by~\eref{delta_final}, assuming $\mu_s$ and $p$ are known. The more general problem of simultaneously reconstructing $\mu_a$, $\mu_s$ and $p$ can also be considered. In either case, what must be investigated is the inversion of the broken-ray Radon transform which is defined as follows. Let $f$ be a sufficiently smooth function. The broken-ray Radon transform is defined by
\begin{equation}
\label{def_Rb}
R_bf(\Br_1,\svec_1;\Br_2,\svec_2)=\int_{BR(\Br_1,\svec_1;\Br_2,\svec_2)} f(\Br) d\Br \ .
\end{equation}
Here $BR(\Br_1,\svec_1;\Br_2,\svec_2)$ denotes the broken ray which begins at $\Br_1$, travels in the direction $\svec_1$ and ends at $\Br_2$ in the direction $\svec_2$. Note that if $\Br_1,\Br_2,\svec_1$ and $\svec_2$ all lie in the same plane and $\svec_1$ and $\svec_2$ point into and out of the slab, then the point of intersection $\BR$ is uniquely determined. Thus it will suffice to consider the inverse problem in the plane and to reconstruct the function $f$ from two-dimensional slices. 

Evidently, the problem of inverting (\ref{def_Rb}) is overdetermined. However, if the directions $\svec_1$ and $\svec_2$ are taken to be fixed, then the inverse problem is formally determined. To this end, we consider transmission measurements in a slab of width $L$ and choose a coordinate system in which $\svec_1$ points in the $\z$ direction. The sources are chosen to be located on the line $z=0$ in the $yz$-plane and are taken to point in the $\z$ direction. The detectors are located on the line $z=L$ and we assume that the angle $\theta$ between $\svec_1$ and $\svec_2$ is fixed. Under these conditions, it can be seen that the solution to (\ref{def_Rb}) is given by
\begin{eqnarray}
\label{soln_int_eq}
f(y,z) = &&\int \frac{dk}{2\pi}e^{-iky}\bigg[ \cot\left(\frac{\theta}{2}\right)F(k,z) - i \cot^2\left(\frac{\theta}{2}\right) k e^{-i \cot(\theta/2)k\xi} \\
&& \times\int_0^z e^{i\cot(\theta/2)k \xi} F(k,\xi) d\xi \bigg]\ ,
\end{eqnarray}
where 
\begin{equation}
\label{def_F}
F(k,\xi)=\lim_{y_2\to \Delta y(\xi)}\left(\frac{\partial}{\partial y_2} + i k \right)\int e^{iky_1}R_bf(y_1,y_1+y_2)dy_1 
\end{equation}
and $\Delta y(\xi)=(L-\xi)\tan\theta$.

Eq.~\eref{soln_int_eq} is the inversion formula for the broken-ray Radon transform. We note that in contrast to x-ray computed tomography, it is unnecessary to collect projections along rays which are rotated about the sample. This considerably reduces the complexity of an imaging experiment. We also note that the presence of a derivative in (\ref{def_F}), as in the Radon inversion formula, means that regularization is essential. Finally, by making use of measurements from multiple detector orientations it is possible to simultaneously reconstruct $\mu_a$ and $\mu_s$. See~\cite{Florescu_2009_1} for further details.

\section{Numerical Inversion Methods}
\label{sect:numerical_inversion}
The starting point for numerical inversion methods is the definition of a variational form whose minimum represents the solution. The most commonly used is a regularised weighted least squares form
\begin{equation}
\Obj(\y,\mfmap(\param)) = \frac{1}{2} || \y - \mfmap(\param)||_{\cvy}^2 +  \frac{1}{2} ||\param||^2_{\cvx}
\label{eq:WLS_reg}
\end{equation}
From the Bayesian point of view \eref{eq:WLS_reg} represents the negative log
of the posterior probability density function
\begin{equation}
\Prob(\param | \y) = \exp(-\Obj(\y,\mfmap(\param))) = \Prob(\y|\param) \Prob(\param)
\label{eq:posterior}
\end{equation}
where $ \Prob(\y|\param)$ is the Gaussian (normal) probability density function of the  distribution of the noise in the data with mean zero and covariance $\icvy$
\begin{equation}
\Prob(\y | \param)=\Prob_{\rm noise}(\y-\mfmap(\param)) = {\mathcal N}(0,\icvy)
\end{equation}
and $\Prob(\param)$ is the Gaussian (normal) prior probability density function of 
the  distribution of $\param$ with mean zero and covariance $\icvx$. 
Thus minimisation of \eref{eq:WLS_reg} represents the \emph{maximum a posteriori} (MAP) estimate of \eref{eq:posterior}. 

%

The forms in \eref{eq:WLS_reg} and \eref{eq:posterior} can be generalised in two main ways. Firstly the prior may be assumed not normal but of the form
\begin{equation}
\Prob(\param) = \exp(-\reglr(\param))
\end{equation}
Secondly the noise need not be Gaussian distributed but could follow an alternative probability model such as Poisson. This consideration leads to alternative \emph{data fit  functionals} of which the commonest is the cross-entropy or KL  divergence
\begin{equation}
KL(\y,\mfmap(\param)) = \int_{\dspace} \left(\y \log\left[\frac{\y}{\mfmap(\param)}\right] - \y + \mfmap(\param) \right)\rmd \y
\end{equation}

From an optimisation point of view the minimisation of \eref{eq:WLS_reg} is 
known as the regularised output least squares solution and can be achieved 
with classical methods if $\Obj(\y,\mfmap(\param))$ is convex. Non-convexity does not 
normally arise from the forward mappings but may arise if the prior is 
non-convex. 
An alternative
formulation of the problem which has its origins in control theory is known as the ``all-at-once'' or ``PDE-constrained method''. In this formulation the solution of the forward problem is treated as a constraint
\begin{eqnarray}
\textnormal{minimise:~}\quad & ||\y - \operator{M}\field||_{\cvy}^2 + \regln\reglr(\param)&\hspace{1cm} \\
\textnormal{subject to:}\quad& \eref{eq:gpde}-\eref{eq:gout}&\nonumber
\label{eq:constrained_PDE_abstract}
\end{eqnarray}

Several of the problems in \sref{sect:modalities} are linear. In addition it is frequently the case that a linearised problem is considered even when the forward problem is non-linear. Suppose we assume a known linearisation point $\param_0$ and we consider the true solution to be a perturbation from this point
$\param = \param_0 + \pparam$ then we consider the minimisation 
\begin{eqnarray}
\pparam_{\ast} 
&=& \mystack{\textnormal{arg min}}{\pparam}\left[\frac{1}{2}||\y - \mfmap(\param_0) - \dfmfmap(\param_0)\pparam||^2_{\cvy} + \regln\frac{1}{2}||\pparam||^2_{\cvx} \right] \label{eq:linearisation1}\\
&=&\mystack{\textnormal{arg min}}{\pparam}\left[\frac{1}{2}||\y^{\delta} - \dfmfmap(\param_0)\pparam||^2_{\cvy} + \regln\frac{1}{2}||\pparam||^2_{\cvx} \right]
\end{eqnarray}
where $\y^{\delta} = \y - \mfmap(\param_0) $
is the change in measurement assumed to be linearly related to $\pparam$. Discretisation of this problem as discussed in \sref{sect:general_framework} leads to the form \eref{eq:lin_discrete}. In the next section we consider several generic linear solvers as applied in optical tomography.

\subsection{Linear Methods\label{sect:linear_methods}}

%

Consider the solution of a linear discrete problem
\begin{equation}
\yv = \mA \fv \quad \yv \in \R^M, \fv \in \R^N, \mA \in \R^M \times \R^N\,.
\label{eq:lin1}
\end{equation}

%
The weighted least-squares term to be minimised is
\begin{eqnarray}
\fv_{\ast} 
&=& \mystack{\textnormal{arg min}}{\fv}\left[\frac{1}{2}||\yv - \mA \fv||^2_{\cvy} + \regln\frac{1}{2}||\fv||^2_{\cvx} \right] \nonumber \\
&=& \mystack{\textnormal{arg min}}{\fv}\left[\frac{1}{2}(\yv - \mA \fv)\tr\cvy(\yv - \mA \fv) + \regln\frac{1}{2}\fv\,\tr \cvx \fv \right]
\label{eq:linlsq3}
\end{eqnarray}

We define factorisations
\begin{equation}
\lcvy\tr\lcvy = \cvy \,\qquad \lcvx\tr\lcvx = \cvx
\end{equation}
and the canonical (dimensionless) variables
\begin{equation}
\tyv = \lcvy\yv\,\qquad \tilde{\mA} = \lcvy\mA\lcvx^{-1}\, \qquad \tfv = \lcvx\fv
\label{eq:var_trans}
\end{equation}
Then the solution to \eref{eq:linlsq3} is given by
\begin{equation}
\tfv_{\ast} = \mystack{\textnormal{arg min}}{\tfv}\left[\tilde{\Obj}_{\regln}(\tyv,\tilde{\mA}\tfv) := \frac{1}{2}||\tyv - \tilde{\mA} \tfv||^2 + \regln\frac{1}{2}||\tfv||^2 \right] 
\label{eq:linlsq_trans}
\end{equation}
In purely algebraic terms, the transformation $\mA \rightarrow \tilde{\mA}$ given 
in \eref{eq:var_trans} is a conditioning step, with $\lcvy, \lcvx^{-1}$ the 
left and right preconditioners respectively\cite{calvetti2007}. 

We now consider some linear inversion methods for \eref{eq:linlsq_trans}.
\subsubsection{Newton methods (overdetermined system)}
A Newton solution for an over determined system ($M > N$) is written
\begin{eqnarray}
\tfv_{\ast} &=& \left(\tilde{\mA}\tr\tilde{\mA} + \regln\mymat{I}\right)^{-1}\tilde{\mA}\tr\tyv\label{eq:Newton_overdet_trans}\\
 \lcvx\fv_{\ast} &=& \left(\lcvx^{-\mathrm{T}}{\mA\tr}\cvy \mA \lcvx^{-1} +
\regln\mymat{I}\right)^{-1}\lcvx^{-\mathrm{T}}\mA\tr\cvy \yv \nonumber\\
\fv_{\ast} &=& \left({\mA\tr}\cvy \mA  +\regln \cvx \right)^{-1}\mA\tr\cvy \yv
\label{eq:Newton_overdet}
\end{eqnarray}
\subsubsection{Newton methods (underdetermined)}
A Newton solution for an underdetermined system ($M < N$) is written
\begin{eqnarray}
\tfv_{\ast} &=& \tilde{\mA}\tr\left(\tilde{\mA}\tilde{\mA}\tr + \regln\mymat{I}\right)^{-1}\tyv\label{eq:Newton_underdet_trans}\\
 \lcvx\fv_{\ast} &=& \lcvx^{-\mathrm{T}}\mA\tr\lcvy\tr \left(\lcvy \mA \icvx{\mA\tr}\lcvy\tr +
\regln\mymat{I}\right)^{-1}\lcvy \yv \nonumber\\
\fv_{\ast} &=& \icvx\mA\tr \left(\mA\icvx \mA\tr  +\regln \icvy\right)^{-1} \yv\label{eq:Newton_underdet}
\end{eqnarray}

\subsubsection{Landweber Method and Steepest Descent}

The steepest descent  method is an iterative reconstruction scheme. 
Beginning with an arbitrary initial estimate $\tfv^{0}$ (usually all 
zeros), we iterate solutions using
\begin{equation}
\tfvk{k+1} = \tfvk{k} + \step_{k} \tsdsk{k}
\label{eq:sd_iter}
\end{equation}
where the steepest direction for \eref{eq:linlsq_trans} is given by
\begin{equation}
\tsdsk{k}= \tilde{\mA}\tr \left( \tyv - \tilde{\mA}\tfvk{k}\right) - \regln \tfvk{k}
\label{eq:sd_trans}
\end{equation}
and
\begin{equation}
\step_{k} = \frac{||\tsdsk{k}||^2}{||\tilde{\mA}\tsdsk{k}||^2 + \regln ||\sdsk{k}||^2 }\,,
\label{eq:sdlinesearchmin}
\end{equation}
is the step length in direction $\tilde{\sds}^{(k)}$ that minimises the one dimensional error function
\[
\tilde{\Obj}(\step) := \tilde{\Obj}_{\regln}(\tyv, \tilde{\mA}(\tfv+\step \tsdsk{k} ))
\] 
The implication is that the vectors in the sequence $\set{\tsdsk{1}, \tsdsk{2},\ldots }$ have the properties \cite{bertero98}
that they can be obtained iteratively
\begin{equation} 
\tsdsk{k+1} = \tsdsk{k} - \step_k \tilde{\mA} \tr \tilde{\mA} \sdsk{k}\,,
\end{equation}
and  satisfy an orthogonality condition
\begin{equation}
	(\tsdsk{k+1}, \tsdsk{k} ) =0\,.
\end{equation}

The Landweber method replaces \eref{eq:sd_iter} with
\begin{equation}
\tfvk{k+1} = \tfvk{k} + \step \tilde{\mA}\tr\left( \tyv - \tilde{\mA}\tfvk{k}\right)
\label{eq:landweber_trans}
\end{equation}
where $\step$ is a relaxation parameter.
Rather than take the exact step that would move the solution to the minimum of 
the error function in the direction $\tsdsk{k}$ the Landweber method 
takes a fixed step.
Thus its convergence is slower than the steepest descent method.

Considering the form of \eref{eq:landweber_trans} we have
\begin{equation}
\fvk{k+1} = \fvk{k} + \step \icvx{\mA}\tr\cvy\left( {\yv} - {\mA}\fvk{k}\right)
\label{eq:landweber_with_priors}
\end{equation}

\subsubsection{Krylov Methods\label{sect:KrylovMethods}}

Consider the general, overdetermined quasi-Newton scheme,
\begin{equation}
\fv = \mymat{H}_{\param}^{-1}\mA\tr \yv\,,
\end{equation}
or the underdetermined scheme
\begin{equation}
\fv = \mA\tr\mymat{H}_{e}^{-1}\yv
\end{equation}
where $\mymat{H}_{\param} \in \R^{N}\times \R^{N}$,  $\mymat{H}_{e} \in \R^{M}\times \R^{M}$ are symmetric positive definite.
Then the Krylov spaces may be defined as the spaces spanned by the set of vectors
\begin{eqnarray}
\Krylsp_{\param} &:=& \textnormal{span}\set{\mA\tr \yv,\mymat{H}_{\param}\mA\tr \yv\ldots \mymat{H}_{\param}^{j}\mA\tr \yv \ldots}\\
\Krylsp_{e}      &:=& \textnormal{span}\set{\mA\tr \yv,\mA\tr\mymat{H}_{e} \yv\ldots \mA\tr \mymat{H}_{e}^{j}\yv,\ldots}
\end{eqnarray}
Deriving these Krylov sequences for the problems \eref{eq:Newton_overdet_trans}, or  \eref{eq:Newton_underdet_trans} both give rise to the same set
\begin{eqnarray}
\tilde{\Krylset}_{\alpha} &=& \set{\tmA\tr\tyv, \tmA\tr\tmA\tmA\tr\tyv + \regln\tmA\tr\tyv, \ldots , \sum_{j'=0}^j \regln^{j'} \left(\tmA\tr\tmA\right)^{j-j'}\tmA\tr\tyv\ldots}\\
&=:& \set{\tKrylvec_{\regln}^{(0)},\tKrylvec_{\regln}^{(1)},\ldots,\tKrylvec_{\regln}^{(j)}\,\ldots } \nonumber
\end{eqnarray}
Defining the unregularised Krylov space of dimension J, via the basis set
\begin{eqnarray}
\tilde{\Krylset} \equiv \set{\tKrylvec^{(j)}} := \set{\left(\tmA\tr\tmA\right)^{j}\tmA\tr\tyv \,,\qquad j = 0\ldots J-1}
\end{eqnarray}
we may derive the Krylov space  for any regularisation parameter $\alpha$ via the basis set
\begin{equation}
\tilde{\Krylset}_{\regln} = \left\{
\begin{array}{lcl}
\tKrylvec_{\regln}^{(0)} &=& \tKrylvec^{(0)}\\
\tKrylvec_{\regln}^{(1)} &=& \tKrylvec^{(1)} + \regln \tKrylvec_{\regln}^{(0)}\\
\vdots &&\\
\tKrylvec_{\regln}^{(j)} &=& \tKrylvec^{(j)} + \regln \tKrylvec_{\regln}^{(j-1)}\\
\end{array}\right.
\end{equation}
Finally we may construct the Krylov space for the original parameters:
\begin{equation}
\Krylset_{\regln} = \left\{
\begin{array}{lcl}
\Krylvec_{\regln}^{(0)} &=& \icvx \mA\tr \cvy \yv\\
\Krylvec_{\regln}^{(1)} &=& \left(\icvx \mA \tr \cvy \mA + \regln\right) \Krylvec_{\regln}^{(0)}\\
\vdots &&\\
\Krylvec_{\regln}^{(j)} &=& \left(\icvx \mA \tr \cvy \mA + \regln\right)\Krylvec_{\regln}^{(j-1)}\\
\end{array}\right.
\end{equation}

Contrast this with the basis defined from \eref{eq:Newton_overdet} 
\begin{equation}
\Krylset^{\mathrm{(ON)}}_{\regln} = \left\{
\begin{array}{lcl}
\Krylvec_{\regln}^{(0)} &=& \mA\tr \cvy \yv\\
\Krylvec_{\regln}^{(1)} &=& \mA \tr \cvy \mA \mA\tr \cvy \yv + \regln\cvx\mA\tr \cvy \yv\\
\vdots &&\\
\Krylvec_{\regln}^{(j)} &=& \left(\mA \tr \cvy \mA + \regln\cvx\right)^{j-1} \mA\tr \cvy \yv \\
\end{array}\right.
\end{equation}
or that from \eref{eq:Newton_underdet}
\begin{equation}
\Krylset^{\mathrm{(UN)}}_{\regln} = \left\{
\begin{array}{lcl}
\Krylvec_{\regln}^{(0)} &=& \icvx \mA\tr \yv\\
\Krylvec_{\regln}^{(1)} &=& \icvx \mA\tr \mA \icvx \mA\tr \yv + \regln \icvx \mA\tr \icvy \yv\\
\vdots &&\\
\Krylvec_{\regln}^{(j)} &=& \icvx \mA\tr \sum_{j'}^{j} \regln^{j'}\left(\mA \icvx \mA \tr\right)^{j-j'} \icvy^{j'} \yv \\
\end{array}\right.
\end{equation}

%

In the conjugate gradient method a second sequence of $\tmA$-conjugate directions $\set{\tcgsk{1},\tcgsk{2},\ldots}$ is constructed and \eref{eq:sd_trans} becomes
\begin{equation}
\tfvk{k+1} = \tfvk{k} + \stepk{k} \tcgsk{k}\,.
\label{eq:sd_cg}
\end{equation}
The direction $\tcgsk{k}$ is given by 
\begin{equation}
\tcgsk{k+1} = \tsdsk{k+1} + \cgpk{k}\tcgsk{k}\,
\end{equation}
where the Gram-Schmidt constants $\cgpk{k}$ are given by one of the methods
\begin{equation}
\cgpk{k} = \left\{
\begin{array}{ll}
\textnormal{Fletcher-Reeves:}
& \frac{||\tsdsk{k+1}||^2}{|| \tsdsk{k}||^2}       \\
\textnormal{Polak-Ribiere:}
& \frac{\left<\tsdsk{k+1}, \tsdsk{k+1}-\tsdsk{k}\right>}{\left<\tcgsk{k+1}, \tsdsk{k+1}-\tsdsk{k}\right>}
\end{array}
\right.
\end{equation}
and the step length (c.f. \eref{eq:sdlinesearchmin}) is given by
\begin{equation}
\stepk{k} = \frac{||\sdsk{k}||^2}{||\tilde{\mA}\tcgsk{k}||^2 + \regln ||\tcgsk{k}||^2 }\,,
\label{eq:cglinesearchmin}
\end{equation}

We note that given a Krylov vector (an image in parameter space) the next vector is formed by the steps
\begin{equation}
\textnormal{forward project} \rightarrow \textnormal{filter in data space} \rightarrow \textnormal {back project} \rightarrow \textnormal {filter in image space}
\label{eq:FFBFscheme}
\end{equation}
Thus the sequence can be formed without building the matrix $\mA$ and constitutes a \emph{matrix free} approach.

\subsubsection{Row and Column Normalisation\label{sect:row_column_normalisation}}

Let us define $\myvec{a}^{R}_i \in \R^N$ as the $i^\mathrm{th}$ row of $\mA$ and
$\myvec{a}^{C}_j \in \R^M$ as the $j\mathrm{th}$ column. Then we can define a
number of purely algebraic conditioning matrices
\begin{eqnarray}
\begin{array}{ccc}
  \mymat{R}_1 &=& \diag\left[\set{||\myvec{a}^{R}_i||_1; i=1\ldots M}\right]\\
  \mymat{C}_1 &=& \diag\left[\set{||\myvec{a}^{C}_j||_1; j=1\ldots N}\right]
\end{array}\label{eq:L1rowcolsum}\\
\begin{array}{ccccc}
  \mymat{R}^2_2 &=& \diag\left[\set{||\myvec{a}^{R}_i||^2_2; i=1\ldots M}\right] &=&\ \diag\left[\mA \mA\tr \right]\\
  \mymat{C}^2_2 &=& \diag\left[\set{||\myvec{a}^{C}_j||^2_2; j=1\ldots N}\right] &=& \diag\left[\mA\tr \mA \right]
\end{array}\label{eq:L2rowcolsum}
\end{eqnarray}
where the notation $||.||_p := \left(\sum .^p\right)^{1/p}$ defines the p-norm of a vector.


For the special case where $\mA$ contains purely nonnegative elements the forms
$\mymat{R}_1$, $\mymat{C}_1$ can be constructed as row and columns sums:
\begin{equation}
 \mymat{R}_1 = \mA \myvec{1}\,,\qquad \mymat{C}_1 = \mA\tr \myvec{1}
\end{equation}

The choice
\[
    \icvy \leftarrow \mymat{R}_1 \,\qquad \cvx  \leftarrow \mymat{C}_1
\]
leads to the \emph{Simultaneous Algebraic Reconstruction Technique} (SART)
\begin{equation}
\fvk{k+1} = \fvk{k} + \step \mymat{C}_1^{-1}\mA\tr\mymat{R}_1^{-1}\left( {\yv} - {\mA}\fvk{k}\right)
\label{eq:SART}
\end{equation}



This is suitable for emission tomography and also for fluorescence optical tomography, where the sought for parameters are positive and the matrix $\mA$ represents a probability of a photon emitted at voxel $i$ arriving at detector $j$.


The choice
\[
    \cvy \leftarrow \mymat{R}^2_2 \,\qquad \cvx  \leftarrow \mymat{I}
\]
leads to the \emph{Simultaneous Iterative Reconstruction Technique} (SIRT) 
\begin{equation}
\fvk{k+1} = \fvk{k} + \step \mA\tr\mymat{R}_2^{-2}\left( \yv - \mA\fvk{k}\right)
\label{eq:SIRT}
\end{equation}

A particular preconditioner that is used in conjugate gradient solvers is given by (Golub and van Loan 1989)
\begin{equation}
\cvx = \mymat{M} = \diag\left[\mA \tr \cvy \mA \right]
\end{equation}
which seeks to equalise the curvature of the objective function along the coordinate axes (``sphering'').

The SIRT and SART algorithms derived from the Landweber method can also be
thoguht of as derived from the Algebraic reconstruction technique (ART), or Kaczmarz method~\cite{kaczmarz37}. This is a row-action method that computes an update to the solution by processing one row of the linearised system at a time. The inner loop generates an update from row $i$ by
\begin{equation}
\fvk{k+1} = \fvk{k} + \varpi_k \frac{\yv_i - \langle\myvec{a}^{R}_i,\fvk{k} \rangle}{||\myvec{a}^{R}_i||^2 } \myvec{a}^{R}_i
\end{equation}
where $\varpi_k$ is a relaxation parameter. The loop over rows of $\mA$ is repeated $n$ times. The rows $i$ are processed in randomised order, which has been shown to improve the convergence rate of the method~\cite{strohmer06}.
\subsection{Multiplicative Methods}

Consider the subsitution
\begin{equation}
\y = \rme^{\logy} \leftrightarrow \logy = \log \y \label{eq:logy}
\end{equation}

An optimisation problem in terms of \eref{eq:logy} may be defined
\begin{equation}
\param_{\ast} 
= \mystack{\textnormal{arg min}}{\param}\left[\frac{1}{2}||\log \y - \log \mfmap(\param)||^2_{\cvy} + \regln\reglr(\param) \right]
\end{equation}

Taking the linearisation of this as in \eref{eq:linearisation1} leads to a 
problem
\begin{eqnarray}
\pparam_{\ast} 
&=& \mystack{\textnormal{arg min}}{\pparam}\left[\frac{1}{2}||\log \y - \log \mfmap(\param_0) - \frac{\dfmfmap(\param_0) }{\mfmap(\param_0)}\pparam||^2_{\cvy} + \regln\reglr(\param) \right] \\
&=& \mystack{\textnormal{arg min}}{\pparam}\left[\frac{1}{2}||\logy^{\delta}  - \hat{A} \pparam||^2_{\cvy} + \regln\reglr(\param) \right]
\end{eqnarray}
where $\logy^{\delta} = \log \y - \log \mfmap(\param_0)$ and the linear operator
\begin{equation}
\hat{A} = \frac{\dfmfmap(\param_0) }{\y_0}
\end{equation}
with $\y_0 = \mfmap(\param_0)$. Similarly, considering
\begin{equation}
\param = \rme^{\logparam} \leftrightarrow \logparam = \log \param \label{eq:logparam}
\end{equation}
leads to a linearised problem
\begin{eqnarray}
\logparam^{\delta}_{\ast} 
&=& \mystack{\textnormal{arg min}}{\logparam^{\delta}}\left[\frac{1}{2}||\log \y
- \log \mfmap(\param_0) - \frac{\dfmfmap(\param_0)\param_0 }{\mfmap(\param_0)}\logparam^{\delta}||^2_{\cvy} + \regln\reglr(\param) \right] \\
&=& \mystack{\textnormal{arg min}}{\logparam^{\delta}}\left[\frac{1}{2}||\logy^{\delta}  - \breve{A} \logparam^{\delta}||^2_{\cvy} + \regln\reglr(\param) \right] \label{eq:log_multiplicative}
\end{eqnarray}
with
\begin{equation}
\breve{A} = \frac{\dfmfmap(\param_0) \param_0}{\y_0}
\end{equation}

Taking the discrete version we arrive at

\begin{equation}
\tmA = \diag[1/\yyv_0] \lcvy\, \mA \, \lcvx^{-1}\diag[\xv_0]
\end{equation}
which means that the covariances have been transformed :
\begin{equation}
\icvy \rightarrow \diag[\yyv_0]\icvy \diag[\yyv_0]\,,\qquad 
\icvx \rightarrow \diag[\xv_0]\icvx \diag[\xv_0]\,
\end{equation}


The steepest descent scheme for \eref{eq:log_multiplicative} is given by
\begin{equation}
\xvk{k+1} = \xvk{k}\odot\exp\left[\tmA \tr \log \left(
\frac{\yyv}{\mfmap(\xvk{k})} 
\right)-\regln \reglr'(\xvk{k})\right]
\end{equation}


when the mapping $\mfmap$ is linear (as in fluorescence DOT for example) and 
when regularisation is only taken implicitly, we arrive at the
\emph{Multiplicative ART} (MART) scheme 
\begin{equation}
\xvk{k+1} = \xvk{k} \odot \exp\left[\mA\tr \log \left(\frac{\yv}{ \mA \xvk{k}}\right)\right]
\end{equation}

If we start instead from the KL divergence
\begin{equation}
KL(\yyv,\mfmap(\xv)) = \sum \yyv \log\left[\frac{\yyv}{\mfmap( \xv)}\right] - \yyv + \mfmap(\xv)
\end{equation}
we get to the \emph{ Maximum Likelihood Expectation Maximisation} (MLEM)
 algorithm
\begin{equation}
\xvk{k+1} = \frac{\xvk{k}}{\mA\tr \myvec{1}}\odot \left[\mA\tr \left(\frac{\yyv}{ \mfmap(\xvk{k})}  \right)\right]
\end{equation}
In this case the explicit regularisation leads to the MAP-EM algorithm
\begin{equation}
\xvk{k+1} = \frac{\xvk{k}}{\mA\tr \myvec{1} + \regln \reglr'(\xvk{k})}\odot \left[\mA\tr \left(\frac{\yyv}{ \mfmap(\xvk{k})}  \right)\right]
\end{equation}


\subsection{Non-Linear Methods\label{sect:Optimisation}}
%
%
%


\subsubsection{Gauss-Newton Method}
The Gauss-Newton method can be considered as the iterative 
minimisation of the quadratic 
Taylor-series approximation to $\Obj(\y,\mfmap(\param))$
\begin{eqnarray}
\Obj(\yv,\mfmap(\xvk{k} + \fv)) &\simeq& \left.
\frac{1}{2}||\yyv - \mfmap(\xvk{k}) - \dfmfmap(\xvk{k})\fv||^2_{\cvy} +\right.\nonumber \\
&& \quad \left. \regln \left(\reglr(\xvk{k}) + 2\left<\fv,\rPDE \xvk{k}\right> + \left<\fv,\rPDE\fv\right> \right)\right.
\label{eq:GNK1}
\end{eqnarray}
where $\rPDE(\xvk{k}) : \sspace \rightarrow \sspace$ represents the mapping induced by the linearisation of the functional $\reglr$. Minimisation of \eref{eq:GNK1} is given by the solution to
\begin{eqnarray}
\qquad \left(\mA\tr\cvy\mA + \regln \rPDE\right) \fv &=& \mA\tr\cvy\left(\yyv - \mfmap(\xvk{k})\right) - \rPDE\xvk{k} \ , \\
\hspace{3.1cm}\mymat{H}\,\fv &=&  -\myvec{g}\label{eq:generic_Newton}
\end{eqnarray}
with Hessian $\mymat{H} = \left(\mA\tr\cvy\mA + \regln \rPDE\right)$ and 
gradient $\myvec{g} = \rPDE\xvk{k} - \mA\tr\cvy\left(\yyv - \mfmap(\xvk{k})\right) $. 
From the form of the Hessian we can see it is guaranteed to be symmetric 
non-negative definite provided that $\reglr$ is convex, and therefore $\fv$ 
is guaranteed to be in a descent direction for $\Obj$.
Solution of \eref{eq:generic_Newton} can be carried out with any of the 
methods in \sref{sect:linear_methods}. In particular 
the \emph{Gauss-Newton-Krylov} method uses a Krylov solver for 
\eref{eq:generic_Newton}. For a badly ill-posed problem the size of the 
requisite Krylov space can be quite small, and may be constructed entirely 
through forward and adjoint solutions and image and data filtering operations as in \eref{eq:FFBFscheme}, and therefore also as a matrix-free approach.

Since the approximation in \eref{eq:GNK1} is only locally quadratic, the update 
given by solving \eref{eq:generic_Newton} may i) not be optimal or may ii) not 
be a descent step (if the Hessian is not symmetric positive definite). 
There are two strategies for resolving these problems.
The \emph{Damped Gauss-Newton} method is a globalisation strategy that 
addresses the first problem. In this approach the update direction is used in a 
one-dimensional line search to find an update step $\step_k$  that 
minimises $\Obj$ along this direction:
\begin{equation}
\step_k = \textnormal{arg min}_{\step} \Obj(\yyv,\mfmap(\xvk{k} + \step \fv))
\label{eq:linsearch1}
\end{equation}
Note that the full non-linear mapping $\mfmap$ is used in \eref{eq:linsearch1},
not its linearisation. The \emph{\LM} method can also address the second problem.
In this approach a control parameter $\LMp$ is used to modify the Hessian
\[
\mymat{H} \rightarrow \mymat{H} +\LMp \Id \,.
\]
When $\LMp$ is large the update $\fv$ tends towards the steepest descent 
direction with increasingly shorter steps. When $\LMp$ is small the update 
tends towards the Newton direction. The idea in the algorithm is to decrease 
$\LMp$ whenever the preceding step reduced $\Obj$ but to increase it and 
re-solve \eref{eq:GNK1} if the preceding step increased $\Obj$. The role 
of $\LMp$ also serves to modify the eigenspectrum of $\mymat{H}$ and force it 
to be symmetric positive definite. The role of column normalisation of $\mA$ 
(see \sref{sect:row_column_normalisation}) is often crucial in the use of 
the \LM algorithm since it can ``sphere'' the level sets of the local quadratic
approximation. 
Typically the \LM method is required for highly non-linear problems,
which is not the case in optical tomography. A comparison of \LM and 
damped Gauss-Newton methods can be found in \cite{schweiger2005a}. The Gauss-Newton method was used for the inverse RTE problem in optical tomography in\cite{tarvainen2008a}

\subsubsection{Nonlinear conjugate gradient method}

The algorithm is effectively the one in \sref{sect:KrylovMethods} with line 
search as in \eref{eq:linsearch1} replacing the calculation of the step length 
give in \eref{eq:cglinesearchmin}. Rather than present it in the canonical 
variables we can present it in the normal variables (Fletcher-Reeves version)

%

\begin{algorithmic}
\STATE Set $\mvk{g}{0} = \regln \reglr'(\xvk{0})- \afmfmap\cvy\left(\yyv - \mfmap(\xvk{0})\right)  $
\STATE $\cgsk{0} = -\icvx \mvk{g}{0}$ \\
\FOR{$k=1,...$: }
\STATE $\step_k  = \textnormal{arg min}_{\step}\Obj(\xvk{k-1} + \step \cgsk{k-1})$ \\
\STATE $\xvk{k} = \xvk{k-1} + \step_k  \cgsk{k-1}$ \\
\STATE $\mvk{g}{k} = \regln \reglr'(\xvk{k})-\afmfmap\cvy\left(\yyv - \mfmap(\xvk{k})\right) $ \\
\STATE $\cgp_k= \frac{\left<\mvk{g}{k},\icvx \mvk{g}{k}\right>}{\left<\mvk{g}{k-1}, \icvx \mvk{g}{k-1}\right>} $ 
\STATE $\cgsk{k} = -\icvx\mvk{g}{k} + \cgp_{k} \cgsk{k-1}$ \\
\ENDFOR
\end{algorithmic}

Note that $\tilde{\tilde{\sds}}_{k} = -\icvx\mvk{g}{k} $ represents the 
application of the preconditioner $\myvec{M}^{-1}=\icvx$ to the gradient. 
If we assume $\reglr'(\xv) \equiv \rPDE \xv = \lcvx\tr \lcvx \xv = \cvx \xv$ 
then this conditioned gradient will be
\[
\tilde{\tilde{\sds}}_{k} = \icvx \afmfmap\cvy\left(\yyv - \mfmap(\xvk{k})\right) - \alpha\xvk{k}
\]
Note the similarity to the general scheme in \eref{eq:FFBFscheme}. Also note 
that the variables $\tilde{\tilde{\sds}} $ 
are no longer {dimensionless}. They have the reciprical 
dimensions of the original parameters themselves.

Non-linear Conjugate Gradients is usually restarted after a certain number of 
iterations, as for any descent method, may be combined with projection onto 
convex sets for enforcement of constraints such as lower and upper bounds, 
but at the expense of loss of its conjugacy properties.
Nonlinear CG was used in optical tomography for the diffusion based problem 
in\cite{arridge98c} and for the RTE based problem in~\cite{kim_h2007}

\subsubsection{Limited-memory BFGS method}
 The BFGS (Broyden-Fletcher-Goldfarb-Shanno) algorithm is a quasi-Newton approach that builds up estimates $\tilde{\mymat{H}}^{-1}$ of the inverse of the Hessian 
matrix $\mymat{H}$~\cite{Nocedal,liu89}. The update rule
\begin{equation}
\xvk{k+1} = \xvk{k} - \lambda_k \tilde{\mymat{H}}^{-1}_k \mvk{g}{k}
\end{equation}
is employed, where $\tilde{\mymat{H}}_k$ is updated at each iteration by the 
formula 
\begin{eqnarray}
\tilde{\mymat{H}}^{-1}_{k+1} &=& \mymat{V}_k\tr \tilde{\mymat{H}}^{-1}_k \mymat{V}_k + \rho_k \mvk{d}{k} \mvk{d}{k}\tr, \\
\rho_k &=& 1/\left<\mvk{z}{k}, \mvk{d}{k}\right>, \\
\mymat{V}_k &=& \mymat{I} - \rho_k \mvk{z}{k} \mvk{d}{k} \tr, \\
\mvk{d}{k} &=& \xvk{k+1} - \xvk{k}, \\
\mvk{z}{k} &=& \mvk{g}{k+1} - \mvk{g}{k}
\end{eqnarray}
In the limited memory version of the algorithm (L-BFGS), the approximate matrices
$\tilde{\mymat{H}}_k^{-1}$ are not stored explicitly, but described implicitly 
by a limited number of pairs of vectors $\lbrace \mvk{d}{i}, \mvk{z}{i} \rbrace$,
where in each iteration, a new vector pair is added, and the oldest pair is 
discarded.

\subsubsection{Nonlinear \KZ\, Method}
The non-linear \KZ\, method is widely used in non-linear tomography \cite{natterer2000}. It is applicable for problems such as DOT which use multiple sources.
Using one source at a time, the subset of the data from this source is back-projected and added to the solution

\begin{equation}
\xvk{k+1} = \xvk{k} + \operator{C}\afmfmap_{i(k)} \cvy \left( \yyv_{i(k)} - \mfmap_{i(k)} (\xvk{k}) \right)
\end{equation}

where the index $i(k)$ refers to the subset of the data accessed on the $k^{\mathrm{th}}$ iteration of the algorithm. The operator $\operator{C}$ is a simple operator playing the role of a right-preconditioner. A natural choice would be
\begin{equation}
\operator{C} = \left( \mfmap_{i(k)}\afmfmap_{i(k)} + \alpha \mymat{I}\right)^{-1}
\end{equation}
but this may be hard to compute. Applications in optical tomography can be found
in \cite{dorn98,rodriquez2009}


\subsubsection{Iterative Coordinate Descent}

A relatively simple approach to generating step directions is to update one pixel
at a time. It is equivalent to taking a descent direction $\sds$, projecting to 
one unit axis, and  iterating through the dimensions of $\sspace$.  This is Gibbs
sampling procedure for the posterior distribution within the Bayesian framework.

Writing $\mvk{e}{k}$ for a unit with 1 in the $k\th$ entry and 0 otherwise, the
ICD method considers
\begin{eqnarray}
\tau_k &=& \textnormal{arg min}_k \left[\Obj(\mvk{x}{k} + \tau \mvk{e}{k})\right]\\
\mvk{\param}{k+1} &=& \mvk{\param}{k} + \tau_k\mvk{e}{k}
\end{eqnarray}

Taking a local linearisation around $\mvk{\param}{k}$ the one dimensional objective functional to be minimised is
\begin{eqnarray}
\Obj(\mvk{x}{k} + \tau \mvk{e}{k}) &=& \frac{1}{2}||\myvec{\y}-\mfmap(\mvk{\param}{k})||^2_{\cvy} - \tau \left<\myvec{\y}-\mfmap(\mvk{\param}{k}),\cvy \mvk{a}{k}^{C}\right> + \nonumber \\
&&\qquad \frac{1}{2}\tau^2 ||\mvk{a}{k}^{C}||^2_{\cvy} + \regln\reglr(\mvk{x}{k} + \tau \mvk{e}{k})
\label{eq:ICD1}
\end{eqnarray}
If the prior is Gaussian then \eref{eq:ICD1} is a weighted least-squares problem with minimum
\begin{equation}
\tau_k = \frac{\left<\myvec{\y}-\mfmap(\mvk{\param}{k}),\cvy \mvk{a}{k}^{C}\right> -  \left<\mvk{e}{k},\cvx \xvk{k}\right> }{ ||\mvk{a}{k}^{C}||^2_{\cvy} + \left<\mvk{e}{k},\cvx \mvk{e}{k}\right> }
\end{equation}
but it is relatively simple also for a non-Gaussian prior since evaluation of the prior usually only involves  neighbouring pixels. Projection onto constraint sets is also efficient since once the constraints are imposed for one pixel it is not revsited. However the evaluation of the likelihood is only efficient for an explicit matrix method and not for the matrix-free approach.

The ICD method was used in PET by Fessler\cite{fessler94} and in Optical Tomography by Bouman and Webb~\cite{ye99a,milstein2003a}. Acceleration was achieved by a 
multiresolution strategy in\cite{ye2001a}.

\subsubsection{PDE Constrained Method\label{sect:PDE_constrained_method}}

%
The PDE constrained method considers the approach presented in \eref{eq:constrained_PDE_abstract}. To implement the method we consider Lagrangian (dual) fields
$\lgnfield$ and define an objective function 
\begin{equation}
\Lgn_{\sind}(\param,\field_{\sind},\lgnfield_{\sind}) = ||\y_{\sind} - \operator{M}\field_{\sind}||_{\cvy}^2 + \regln\reglr(\param) + \left<\lgnfield_{\sind} ,(\PDE(\param)\field_{\sind} -q_{\sind})\right>_{\domain}
\label{eq:PDE_Lagrangian}
\end{equation}
where $q_{\sind}$ represents an equivalent source for the boundary condition 
\eref{eq:gin}. When considering all sources the full Lagrangian is
\begin{equation}
\Lgn(\param,\field,\lgnfield) = \sum_i \Lgn_{\sind}(\param,\field_{\sind},\lgnfield_{\sind})
\end{equation}

The minimum of \eref{eq:PDE_Lagrangian} occurs where the first variation
\begin{eqnarray}
\Lgn_{\sind,\param} &=& \regln\reglr'(\param) + \langle\lgnfield_{\sind}, \PDE_{\param}\field_{\sind}\rangle_{\domain} \label{eq:Lgn_x}\\
\Lgn_{\sind,\field} &=& \aPDE\lgnfield_{\sind} - \myop{M}^{\ast} \cvy \left(\y_{\sind} - \operator{M}\field_{\sind}\right) \label{eq:Lgn_U}\\
\Lgn_{\sind,\lgnfield} &=& \PDE\field_{\sind} - q_{\sind} \label{eq:Lgn_Z}
\end{eqnarray}
becomes zero. We recognise \eref{eq:Lgn_U} as the equation for the backprojected 
field of the residual difference between the data $\y_{\sind}$ and the 
measurement of the direct field $\field_i$. We 
recognise $\PDE_{\param} = \ptlop(\pparam)$ as the potential operator
introduced in \eref{eq:ptlop_def}. Therefore if \eref{eq:Lgn_Z} and \eref{eq:Lgn_Z} are both equated to zero,  \eref{eq:Lgn_x} is the gradient of the negative log posterior of a MAP estimation scheme. If instead $\set{-\Lgn_{\sind,\param} ,-\Lgn_{\sind,\field}, -\Lgn_{\sind,\lgnfield}}$ is used as an update direction for 
$\set{\param,\field,\lgnfield}$ then we may solve for all variables simultaneously
without requiring convergence for any one variable until termination. 

Taking the  second variation of $\Lgn$ results to a Newton system
\begin{equation}
\mtrx{ccc}{\regln\reglr''(\param) & \aPDE_{\param}\lgnfield_{\sind}  & \PDE_{\param}\field_{\sind} \\
 \aPDE_{\param}\lgnfield_{\sind} & \myop{M}^{\ast}\cvy\myop{M} & \PDE\\
\PDE_{\param}\field_{\sind} & \PDE & 0
}
\mtrx{c}{\param^{\delta}\\\field_{\sind}^{\delta}\\\lgnfield_{\sind}^{\delta}}
= -\mtrx{c}{\regln\reglr'(\param) + \langle\lgnfield_{\sind}, \PDE_{\param}\field_{\sind}\rangle_{\domain}\\ \aPDE\lgnfield_{\sind} - \myop{M}^{\ast}\cvy \left(\y - \operator{M}\field_{\sind}\right)\\ \PDE\field_{\sind} - q}
\label{eq:Lgn_KKT}
\end{equation}

%

In many applications the system in \eref{eq:Lgn_KKT} is simplified by taking the Schur complement of the complete Hessian and by using the Gauss-Newton or other quasi-Newton approximations for the solution scheme~\cite{abdoulaev2005,bangareth2008a,kim_h2009}. For the time-domain problem, even with such Hessian reduction techniques, Newton methods are infeasible but first-order descent schemes can still be used\cite{soloviev2008a}.

%




\subsection{Error and Prior Modelling\label{sect:errors_priors}}

So far we did not discuss the forms of the data covariances $\icvy$ or the prior
$\reglr(\param)$. The format of noise usually be predicted on physical grounds.
It is usually takem to be zero-mean Gaussian noise with a possibly non-white 
covariance. When considering photon counting measurements the implication is 
that the noise should be Poisson with variance $\sigma^2_j = \y_j$. Assuming 
sufficient signal to approach the Central Limit Theorem would lead to an 
equivalent Gaussian model of additive noise with
\begin{equation}
\icvy = \diag[\y] \quad \Leftrightarrow \quad \lcvy = \diag\left[\frac{1}{\y^{1/2}}\right]
\end{equation}
A more commonly used model assumes that equal numbers of photons are collected at each detector leading to a constant relative error and the formal covariance structure 
\begin{equation}
\icvy = \diag[\y^2] \quad \Leftrightarrow \quad \lcvy = \diag\left[\frac{1}{\y}\right]
\end{equation}
which also corresponds to the error model implicit in using the Rytov series in 
place of the Born series for the linearised model. Within this photon counting 
paradigm there is no correlation between errors on different detectors. However, 
when considering the actual errors between measured and modelled data the 
correlation is far from being negligible, and their mean  is far from being zero.
This discrepency can be understood by formally considering the modelling error as a random variable\cite{kaipio05}.
\begin{equation}
\Prob(\meas - \mfmap_h(\param_{\mathrm{true}})) \sim \myop{N}(\overline{\epsilon},\icvy+\mymat{\Gamma}_{\epsilon})
\end{equation}
where $\mfmap_h$ is an approximate model with the accuracy of the computational method employed, $\overline{\epsilon}$ is the model bias representing the discrepency between this model and the real Physics, and $\mymat{\Gamma}_{\epsilon}$ is the model error covariance. In practical applications of optical tomography an estimate of the bias is made by measureing a reference problem and comparing it to a reference model
\begin{equation}
\overline{\epsilon} = \y^{\mathrm{ref}} - \mfmap_h(\param_{\mathrm{ref}})
\end{equation}
which leads to a corrected model
\begin{equation}
\mfmap_{\mathrm{corrected}}(\param) = \mfmap_h(\param) + \y^{\mathrm{ref}} - \mfmap_h(\param_{\mathrm{ref}})
\end{equation}
Clearly the corrected model and the measured data agree exactly for the value
$\param = \param_{\mathrm{ref}}$ and the assumption is that they will also agree at nearby values. This is questionable.

A more principled approach is the {\em approximation error method} \cite{kaipio05}.
In this approach the statistical properties of the modelling error are estimated 
by by sampling over a plausible distribution of solutions and comparing the 
modeled data from each sample with ``measured'' data from this sample. In 
practice the measured data could be also be modelled but with  a much more 
accurate technique. 
In \cite{arridge2006a} this technique was shown to
result in reconstructed images using a relatively inaccurate forward
model that were of almost equal quality to those using a more
accurate forward model; the increase in computational efficiency was
an order of magnitude.

Choice of regularisation $\reglr$ is critical and difficult to justify unequivocally. Assuming a form
\begin{equation}
\Psi(\param) = \domintvar{\psi\left(|\nabla \param|)\right)}{\r},
\label{eq:Psi_variational}
\end{equation}
where  $\psi : \sspace \rightarrow \sspace$ is an image to image mapping leads to the linearisation
\begin{equation}
\cvx = \rPDE(\param) = \nabla^{\mathrm{T}}\diffusivity (\r) \nabla
\end{equation}
where the \emph{diffusivity} is given by
\begin{equation}
\diffusivity := \frac{\psi'\left(|\nabla\param|\right)}{|\nabla \param|} 
\end{equation}
Note that the matrix $\icvx$ is a covariances whose entries represent a correlation between  pixels. On the other hand $\cvx$, is sparse, representing the local relationship between neighbours. We could therefore specify $\icvx$ in a number of ways
\begin{enumerate}
\item From a database of representative images 
\item By specifying a Markov Random field $\cvx \equiv \mymat{\MRF}$
\item By specifying  a PDE $\cvx \equiv \nabla\tr { \diffusivity \nabla}$
\end{enumerate}



Methods such as first-order Tikhonov (Phillips-Twomey),
Total Variation (TV), or Generalised Gaussian Markov Random Field
methods (GGMRF)\,\cite{bouman93} make statistical assumptions about
the distribution of edges. They pose assumptions about the regularity (e.g., smoothness) 
of the solution, the regularity being measured in terms of some norm of the solution and its differentials\cite{douiri2005}.
In addition the ``lumpy background''
prior provides an effective method of
estimating information content of different imaging
systems\,\cite{pineda2006}.

Image processing
methods for computer vision offer a large set of techniques for
de-noising and segmentation based on anisotropic diffusion
processes\,\cite{acar94,geman84,rudin92a,vogel96}.
These techniques formulate directly a PDE for image flow, rather than
as the Euler equation of a variational form and can be considered more general.
%
%
%
More specific prior information can be incorporated if we assume that some approximate knowledge of the objects such as shape topology and intra-region parameter regions is available.
\cite{kolehmainen99a} 

\section{Shape Based Methods}
\label{sect:shape_methods}

Optical Tomography as discussed so far in this article is an example of a \emph{parameter
identification} problem because the reconstructed images represent
the parameters of a model of light propagation. In common with many
medical imaging modalities the reconstructed images are not an end in
themselves; they need to be analyzed for structural and functional
information, including classification of regions, segmentation and
cross validation with other modalities.

In the particular case of OT, this post-processing step
has some drawbacks.  In reconstruction, we use small sets of measurements to
reconstruct a large image with many parameters.  In the subsequent
analysis stage we analyze these parameters to categorise the image
into only a few discrete categories
e.g. to estimate a classification into discrete regions, or to determine
the parameters of a low-dimensional shape model. In
each case, the final result is of considerably smaller
dimensionality than the intermediate reconstruction.
An alternative approach would be to integrate the classification or
segmentation with the reconstruction.  This problem may be better
posed than the two-stage approach, as we are only moving from the
sparse measurements to another low-dimensional space.

In many biomedical applications, the parameter distribution which is
sought in the inverse problem, contains some sorts of interfaces. 
These can for example be boundaries of some tumor or hematoma, 
or the interfaces between different organs or tissue types, or
the boundaries of regions filled with some tracer or marker 
substance. These interfaces are typically not recovered well in
classical reconstruction schemes due to the above mentioned need
for relatively strong regularization tools. Most of these 
regularization tools penalize variations or gradients in the
parameter distribution, which yields oversmoothed reconstructions. 
Therefore, interfaces and boundaries are blurred and region
boundaries cannot easily be detected. However, in many applications
it is important to be able to find the boundaries of certain
inclusions as accurately as possible. Whereas inside and outside these
subregions 
the parameters might not vary much, often across the
interfaces significant jumps occur in the tissue 
parameters. This motivates the possibility of direct \emph{shape-based} reconstruction methods, which generally can be categorised into \emph{explicit} and \emph{implicit} methods.

Shape reconstruction
techniques for Optical Tomography using the DA as the forward model 
are presented in \cite{kolehmainen99a,kolehmainen2000a,kolehmainen2000b,kilmer00,kilmer03,dorn2004a,dorn2004b,bal2006a,schweiger2006a,zacharopoulos2006a} and using the RTE for the forward model in\cite{arridge2006b}.

\subsection{Explicit Shape Method}

In shape based methods it is assumed that the domain is represented 
as in \eref{eq:subdomains},\eref{eq:domain_interfaces}. The assumption that the optical parameters are constant in subdomains is usually also imposed, but can in principle be relaxed.
In the explicit shape method we represent the boundaries by an explicit parameterisation and construct a forward problem in terms of the coefficients of the parameterisation and possibly the optical parameters inside each domain.

\subsubsection{Parametric representation of surfaces} \label{sect:surface_parameterisation}

A typical parametisation of a closed surface is in terms trigonometric functions (2D) or spherical harmonics (3D) 
\begin{equation}
b_k =  \mtrx{c}{\Re \\ \Im}\left\{ 
\begin{array}{ll}
       \rme^{\myi k \vartheta} & \textnormal{in 2D} \\
        Y_{l}^{m}(\vartheta,\varphi) , k = (l+1)^2 +m& \textnormal{in 3D}
\end{array}
\right.
\end{equation}
 where a purely real basis is constructed from the real and imaginary parts of the complex functions.

These bases may be used to parameterise the radial distance of the surface from an internal point, restricting the analysis to star-shaped objects, or to develop a  a harmonic mapping of the surface onto a circle or sphere, which is made possible by the specific representation of each cartesian component seperately in the basis.

In general, the higher order basis functions are roughly  assumed to represent more detailed characteristics of the surface, whereas the lower order ones describe more the overall features like volume, orientation etc. 

For simplicity we introduce the notation
\begin{equation}\label{gamdef}
\shape= \{\gamma_k\},\quad k = 1,\cdots K
\end{equation}
to describes the finite set of basis coefficients for the surface $\surf$  up to degree $K$.

\subsubsection{The Forward Problem} \label{sect:shape_forward_model}

As well as the set of surface coefficients  $\{\gamma_k\}$ we define the set 
of parameters $\set{\param_{\bdind} }$ for each domain $\domain_{\bdind}$.
We now define the nonlinear forward operator 
  as the mapping from the optical properties $\param=\{\param_{\bdind}\}\, {\bdind=1,\ldots L}$ of the individual regions
$\domain_{\bdind}$, $\bdind =1,\ldots L$, 
and the geometric parameters $\shape=\{\gamma_k\}\, {k=1,\ldots,K}$ 
to the measurements on the surface of $\bdomain$. Combining the measurements from the $S$ independent sources, the forward mapping takes the form
\begin{equation}
\y = \mfmap(\shape,\param)
 \,.
\label{eq:forwaoper}
\end{equation}


\subsubsection{The inverse problem}\label{sect:shapeinverse}

The most general inverse problem to be considered in our framework would be the 
simultaneous reconstruction of shape parameters $\gamma$ and optical parameters $\param$ from the data $\meas$. However, in most practical applications which we 
have in mind, good estimates for the optical parameters in most of the 
regions $\domj$ are available, and only those optical parameters which 
correspond to  unknown anomalies need to be recovered. 
In particular the outer boundary shape $\partial\Omega$ is usually considered known as well as the optical parameters in the outermost region. Relaxation of these assumptions, in particular the first, leads to more complex considerations that have so far not been addressed in optical tomography, but for which methods have been developed in related fields~\cite{kolehmainen2007a,kolehmainen2008a}

As in the parameter identification problem, the key aspect of using the shape based method is to discover a formula for calculating shape sensitivities functions
which form the kernel of the \Frechet\, derivative of \eref{eq:forwaoper}.
Roughly speaking, the shape sensitivity function for shape basis function $b_k$ is given by the spatial sensitivity function $\pmdf$ defined in \eref{eq:pmdf1}, projected onto $b_k$ on $\surf$
\begin{equation}
\y^{\delta}_{\apind,\sind,k} = \left<\pmdf_{\apind,\sind},b_k \right>_{\surf}
\label{eq:shape_sensitivity}
\end{equation}
The precise form of the sensitivity function depends on the forward model. 
In \cite{kolehmainen99a} a FEM model was used and in \cite{zacharopoulos2006a} 
a BEM model. In the latter the model is naturally represented as the union of 
subdomains as in \eref{eq:subdomains} and the domain boundaries are explictly 
meshed. The fields $\U$ and the normal currents $J_n = \bnormvec \cdot \nabla \U$ are explicitly represented, and sensitivity function $\pmdf$ is defined in 
terms of both the field and the normal currents. Shape updating requires only 
movement of the mesh points of the internal surfaces. By contrast, in the FEM 
model the surface $\surf$ intersects to a subset of elements of the mesh. The 
appropriate form of \eref{eq:shape_sensitivity} calls for the integration of products of forward and adjoint fields along split elements.
If the simplified BEM model is used, where the normal currents are removed by 
taking the Schur complement of the discrete system\eref{eq:BEM2layer_Schur}, then the sensitivity functions form FEM and BEM are the same. 

To develop a shape reconstruction scheme, starting from a geometric 
configuration defined by the initial set of shape coefficients $\shapeiter{0}$, we will search for the set $\shaperecon$ that minimises the distance between computed $\mfmap(\gamma,\param)$ and measured data $\meas$ 
\begin{equation}
\shaperecon\; =~\textnormal{arg min}_{\shape}\left[\Obj(\shape):=
\frac{1}{2}\|\meas-\mfmap(\shape,\param)\|^{2}\right]\,.
\label{eq:shapeinvprob}
\end{equation}
In \cite{kolehmainen99a,kolehmainen2000a,kolehmainen2000b,zacharopoulos2006a} a \LM scheme was used for the inversion step
\begin{equation}
\shapeiter{\itind+1}=\shapeiter{\itind}+(\mA_{\itind}\tr\mA_{\itind}+
\lambda\mymat{I})^{-1}\mA_{\itind}\tr(\meas-\mfmap(\shapeiter{\itind},\param)).
\label{eq:shapeLM}
\end{equation}

\subsection{Implicit Shape Method}

The level set technique represents the boundary of a domain as the zero-set of a
 function $\varphi$. Rather than updating the boundary of the domain an update 
rule for $\varphi$ is developed and the optical parameters are defined
\begin{equation}
\param = \left\{ \begin{array}{ll} 
\param_{\mathrm{int}} & \varphi \leq 0 \\
\param_{\mathrm{ext}} & \varphi > 0 
\end{array}\right.
\end{equation}
providing a tool which is able to reconstruct interfaces 
in the region of interest together with certain characteristics
of the interior and exterior subregions. The key difference from the explicit methods is that topological changes in the domains are easily accomodated.
A description of this approach for the DA with examples in 2D is provided in\cite{schweiger2006a} and the 3D case in\cite{schweiger2008a}. 
Level set methods based on the RTE model 
in optical tomography were  presented in
\cite{dorn2004a,dorn2004b}. 
In these papers the scattering function $\must$ 
was assumed known {\em a priori} and only the subdomain boundaries of the
absorption function $\mua$ were reconstructed. 
So far, these methods have not yet been applied to the reconstruction of the 
low scattering void regions.

Due to the presence of two parameters $\mua$ and $\diffcoef$ in the DA, the
application to Optical Tomography uses two level set function $\varphi_{\mua}, \varphi_{\diffcoef}$.
Assuming that the background absorption and diffusion
parameters $\mu_{\mathrm{ext}}$ and $\diffcoef_{\mathrm{ext}}$ are
constant and known, the unknowns of the inverse problem
are therefore $\varphi_{\mua}, \varphi_{\diffcoef}$ and the interior absorption and diffusion
parameters $\mu_{\mathrm{int}}$ and $\diffcoef_{\mathrm{int}}$. The forward mapping is defined
\begin{equation}
\y = \mfmap\left(\diffcoef(\varphi_{\diffcoef},\diffcoef_{\mathrm{int}}),\mua(\varphi_{\mua},\mu_{\mathrm{int}})\right)
\end{equation}
and we formulate the output least squares cost functional
\begin{equation}\label{Met3}
\Obj\left(\varphi_{\diffcoef}, \varphi_{\mua}, \diffcoef_{\mathrm{int}}, \mu_{\mathrm{int}} \right)\,=\,
\frac{1}{2}\left\|
\meas - \mfmap\left(\diffcoef(\varphi_{\diffcoef},\diffcoef_{\mathrm{int}}),\mua(\varphi_{\mua},\mu_{\mathrm{int}})\right)
\right\|^2.
\end{equation}
The goal is to find a mimimizer $(\varphi,\psi,\diffcoef_{_{obj}},
\mu_{_{obj}})$
 of this cost functional.

For solving the inverse problem an evolution
approach is adopted in which  the unknowns of the inverse 
problem are assumed to depend on an artificial evolution
time $t$ and evolve during the reconstruction 
\begin{equation} 
\begin{array}{lcllcl}\dbyd{\varphi_{\diffcoef}}{t}&=& f_{\diffcoef}(t),\qquad &
\dbyd{\varphi_{\mua}}{t}\,&=&\, f_{\mua}(t), \\ 
\dbyd{\diffcoef_{\mathrm{int}}}{t}&=& h_{\diffcoef}(t),\qquad&
\dbyd{\mu_{\mathrm{int}}}{t}\,&=&\, h_{{\mua}}(t).
\end{array}
\end{equation}
where  $f_{\diffcoef}$, $f_{\mua}$, $ h_{\diffcoef}$, $ h_{\mua}$  are forcing terms which point into a descent direction
of the cost (\ref{Met3}). We use  $\surf := \varphi = 0$ to denote the zero set of $\varphi$ and 
\begin{equation}
s = \afmfmap (\meas - \mfmap\left(\diffcoef(\varphi_{\diffcoef},\diffcoef_{\mathrm{int}}),\mua(\varphi_{\mua},\mu_{\mathrm{int}})\right)
\end{equation}
for the ``classical'' descent direction. Then the forcing terms are defined
\begin{eqnarray}
f_{\param} &=& \rPDE^{-1}(\param_{\mathrm{int}} - \param_{\mathrm{ext}})s_{\param} \\
h_{\param} &=& (\param_{\mathrm{int}} - \param_{\mathrm{ext}})s_{\param}|_{\surf}
\end{eqnarray}
where $\rPDE$ is a regularisation term as in \sref{sect:Optimisation}.
Note that the update for the level set $\varphi$ is defined on the whole domain, whereas the update for the parameters is only defined on the zero set interface.
In practice the level set update may be resticted to a narrow band within a specified distance of this interface.
For details we refer to~\cite{schweiger2006a,schweiger2008a}.

The level set approach can easily be generalized 
to more complex scenarios such as smoothly varying or parameterized interior and
exterior parameter profiles.  

\section{Conclusions and Further Topics
\label{sect:OtherTopics}}
Before drawing this review to a conclusion we briefly mention some important topics which space precludes us from developing in further detail.

\subsection{Anisotropy}\label{sect:anisotropy}

As mention in \sref{sec:diffuse}, if we relax 
the assumption that the RTE phase function is angularly invariant, we can develop an anisotropic model. The derivation of this model in terms of a 
spherical harmonic expansion (i.e. the $P_N$ method) leads to a diffusion equation with a tensor diffusion coefficient\cite{heino02,heino2003}
\begin{equation}
\dbyd{}{t}\phidiff (\posvec,t)
-\nabla \cdot \difftensor (\posvec) \nabla \phidiff (\posvec,t) 
+ c \mua(\posvec)\phidiff (\posvec,t) =0
\label{eq:da_time_anisotropic}
\end{equation}

In terms of the inverse problem we are lead to similar difficulties with regard to non-uniquness as occur in impedence imaging. In particular a diffeomorphism of $\domain$ that preserves the Robin-to-Neumann map will be equivalent to a change in the local orientation of the eigenvectors of $\difftensor$. As a mechanism for compensating for such non-uniqueness, the modelling error approach (see \sref{sect:errors_priors}) has been employed,~\cite{heino2004,heino2005}.

\subsection{Refractive Index Variation}\label{sect:refindex}

All the inverse problems we have discussed in this article were based on the 
transport equation or diffusion approximation with constant speed (i.e. 
constant refractive index). However as discussed in \sref{sec:physics}, 
there exist a number of models for variable refractive index including  
\eref{RTE_variable_n} and \eref{diff_eq_n}. In addition there have been several treatments of the radiative transport equation in piecewise homogeneous (i.e. layered) domains with constant absorption and scattering but differing refractive index\cite{cassell2006,cassell2007,garcia2008,elaloufi2007a} wherein semi-analytic methods can be derived in terms of Fresnel boundary conditions at the interfaces between sub-domains.

In terms of an inverse problem ~\cite{khan2005,khan2006} has considered the diffusion approximation version of the problem as in \eref{diff_eq_n} as the basis for an inversion to recover refractive index changes. 

\subsection{Time-varying Optical Tomography}

In real biological applications, the subject being studied is not static but 
time varying. In this case we have the choice to reconstruct time-point by 
time-point or to consider a spatial temporal model.There are two main 
considerations. Firstly the time scale on which $\param$ is changing may be 
comparable to the time scale on which the experiments are performed. In other 
words the data acquisition process may be corrupted by movement. In this case 
there are methods from time-series analysis in which the dynamic changes 
in $\param$ are included in the reconstruction model. A well known example 
is the Kalman filter method which assumes a random walk model in some state 
space variables whose statistics is learnt over the course of time-series 
experiment. This idea was explored in diffuse optical tomography by \cite{kolehmainen2003}.  In \cite{prince2003} the Kalman filter approach was extended to utilise a state space model in which the underlying temporal variation was taken to be a mixture of pseudo-sinusoidal signals with random amplitude and phase, which allowed seperable reconstruction of functional components in brain images.
In \cite{zhang2005b} a fully space-time model was used with data reduction acheived via a Kronecker product representation of forward operators.

The second consideration it that physiological models may be used to predict temporal variation which may be used to constrain the reconstruction. These methods used coupled systems of non-linear equations which are much more complex than those used in a state-space model\cite{diamond2006}. A review can be found in~\cite{huppert2009}
\subsection{Multimodality\label{sect:multimodality}}

By Multimodality imaging we mean the combination
of the varying capacities of two or more medical imaging technologies
in extracting physiological and anatomical information of organisms
in a complementary fashion, in order to enable the consistent retrieval
of accurate and content rich biological information \cite{ardekani96,comtat2002,gindi93a,lu98,Baillet1999, rangarajan2000,sastry97,som98,zaidi2003, Ahlfors04, Babiloni2005}. In optical tomography, the aim is to enhance the low resolutionfunctional image information with high resolution complementary structural information using for example ultrasound \cite{zhu2003,zhu2005uti}, MRI \cite{ntziachristos2002, Brooksby2003}, or  X-rays \cite{Li2003, Zhang2005}.

One method for utilising the prior information is to construct priors that constrain the reconstruction of the optical information to be commensurate with the auxiliary information in a well-defined way. One approach is the use of structural priors which smooth the solution within regions but not across boundaries of regions inferred from the axiliary modality \cite{kaipio99a, brooksby2005, douiri2006a}. An alternative is to define an information theoretic measure of similarity between the reconstructed image and the prior image~\cite{panagiotou2009a}.

The topic of multimodality and the appropriate use of cross-information is rapidly expanding within optical tomography in particular and medical imaging in general.

\section*{References}

\bibliography{abbrsrc,ref,JCS}

\end{document}